\documentclass[12pt]{amsart}

\usepackage[leqno]{amsmath}

\usepackage[mathscr]{euscript}
\usepackage{amssymb}

\usepackage[margin=1.2in]{geometry}

\usepackage{pb-diagram,pb-xy} 
\usepackage[all,cmtip,arrow]{xy} 
\usepackage{quiver}
\usepackage{tikz-cd}

\usepackage{hyperref} 
\hypersetup{colorlinks=true,linkcolor=blue,citecolor=magenta}

\numberwithin{equation}{section}

\newtheorem{theorem}[equation]{Theorem}
\newtheorem{lemma}[equation]{Lemma}
\newtheorem{corollary}[equation]{Corollary}
\newtheorem{proposition}[equation]{Proposition}
\newtheorem{conjecture}[equation]{Conjecture}

\theoremstyle{definition}
\newtheorem{definition}[equation]{Definition}
\newtheorem{construction}[equation]{Construction}
\newtheorem{remark}[equation]{Remark}
\newtheorem{example}[equation]{Example}
\newtheorem{goal}[equation]{Goal}

\newcommand{\isom}{\cong}                       
\newcommand{\homeq}{\simeq}                     
\newcommand{\smsh}{\wedge}                      
\newcommand{\wdge}{\vee}                        
\newcommand{\Smsh}{\bigwedge}                   

\newcommand{\cat}[1]{\mathscr{#1}}              

\providecommand{\hocolim}{\mathop{\rm hocolim}}
\providecommand{\colim}{\mathop{\rm colim}}
\providecommand{\holim}{\mathop{\rm holim}}
\providecommand{\lim}{\mathop{\rm lim}}

\providecommand{\hofib}{\mathop{\rm hofib}}
\providecommand{\thofib}{\mathop{\rm thofib}}

\providecommand{\Map}{\mathop{\rm Map}\nolimits}
\providecommand{\Fun}{\mathop{\rm Fun}\nolimits}
\providecommand{\Nat}{\mathop{\rm Nat}\nolimits}

\providecommand{\Hom}{\mathop{\rm Hom}\nolimits}
\providecommand{\creff}{\mathop{\rm cr}\nolimits}
\providecommand{\Tot}{\mathop{\rm Tot}\nolimits}

\newcommand{\un}[1]{\underline{#1}}

\newcommand{\spaces}{\mathscr{T}op}
\newcommand{\based}{\mathscr{T}op_*}
\newcommand{\spectra}{{\mathscr{S}p}}

\newcommand{\Fin}{\cat{F}\mathrm{in}_*}
\newcommand{\Surj}{\cat{S}\mathrm{urj}_*}

\newcommand{\weq}{\; \tilde{\longrightarrow} \;}      
\newcommand{\epi}{\twoheadrightarrow}           
\newcommand{\into}{\hookrightarrow}

\newcommand{\der}{\partial}                     

\newcommand{\ord}[1]{#1\textsuperscript{th}}

\newcommand{\Alg}{\cat{A}\mathrm{lg}}
\newcommand{\CAlg}{\cat{CA}\mathrm{lg}}

\providecommand{\Pro}{\mathop{\rm Pro}}
\providecommand{\Ind}{\mathop{\rm Ind}}

\providecommand{\Bimod}{\mathop{\rm Bimod}}

\newcommand{\ev}{\mathrm{ev}}
\newcommand{\evb}{\ev_\bullet}
\newcommand{\ex}{\mathrm{ex}}

\newcommand{\Catinf}{\mathscr{C}\mathrm{at}_\infty}

\begin{document}

\title[Infinity-operads and Day convolution]{Infinity-operads and Day convolution in Goodwillie calculus}
\author{Michael Ching}
\address{
	Department of Mathematics and Statistics \\
	Amherst College \\
	PO Box 5000 \\
	Amherst, MA 01002 \\
	USA}
\email{mching@amherst.edu}


\begin{abstract}
We prove two theorems about Goodwillie calculus and use those theorems to describe new models for Goodwillie derivatives of functors between pointed compactly-generated $\infty$-categories. The first theorem says that the construction of higher derivatives for spectrum-valued functors is a Day convolution of copies of the first derivative construction. The second theorem says that the derivatives of any functor can be realized as natural transformation objects for derivatives of spectrum-valued functors. Together these results allow us to construct an $\infty$-operad that models the derivatives of the identity functor on any pointed compactly-generated $\infty$-category.

Our main example is the $\infty$-category of algebras over a stable $\infty$-operad, in which case we show that the derivatives of the identity essentially recover the same $\infty$-operad, making precise a well-known slogan in Goodwillie calculus.

We also describe a bimodule structure on the derivatives of an arbitrary functor, over the $\infty$-operads given by the derivatives of the identity on the source and target, and we conjecture a chain rule that generalizes previous work of Arone and the author in the case of functors of pointed spaces and spectra.
\end{abstract}

\maketitle

The fundamental construction of Goodwillie calculus is, for a functor $F: \cat{C} \to \cat{D}$, a tower of approximations to $F$ that mimics the Taylor series in ordinary calculus. One of the basic principles of this theory is that the fibres of the maps in that tower can be described relatively simply in terms of stable homotopy theory. Goodwillie showed that when $\cat{C}$ and $\cat{D}$ are either the categories of pointed spaces or spectra, the \ord{$n$} homogeneous piece of a functor $F: \cat{C} \to \cat{D}$ is determined by a single spectrum $\der_nF$ together with an action of the \ord{$n$} symmetric group $\Sigma_n$.

A central question in calculus is how to reconstruct the Taylor tower of the functor $F$ (and hence, in cases where the tower converges, the functor $F$ itself) from these homogeneous pieces, i.e. from the symmetric sequence $\der_*F = (\der_nF)_{n \geq 1}$. In the cases where $\cat{C}$ and $\cat{D}$ are each either pointed spaces or spectra, this question was answered in a pair of papers by Greg Arone and the author \cite{arone/ching:2011,arone/ching:2015}. We first showed that, for $F: \cat{C} \to \cat{D}$, the symmetric sequence $\der_*F$ has the structure of a bimodule over the two operads $\der_*I_{\cat{C}}$ and $\der_*I_{\cat{D}}$ formed by the derivatives of the identity functor on the categories $\cat{C}$ and $\cat{D}$. We then showed that the resulting adjunction, between the categories of ($n$-excisive) functors $F: \cat{C} \to \cat{D}$ and ($n$-truncated) bimodules, is comonadic, so that the Taylor tower of $F$ can be recovered from the action of a certain comonad on the bimodule $\der_*F$.

In this paper, we extend the first part of that previous work to a broad class of $\infty$-categories. In particular, we show that the derivatives of the identity functor on any pointed compactly-generated $\infty$-category form a stable $\infty$-operad in a natural way, and that the derivatives of any functor form a bimodule over the appropriate $\infty$-operads.

Note that the approach taken in this paper is significantly different from that of \cite{arone/ching:2011}, and even in the cases of pointed spaces and spectra it gives a new perspective on how the operad structures arise. In particular, this paper provides a new construction of the spectral Lie operad, as an $\infty$-operad, distinct from the cobar construction described in~\cite{ching:2005}.

One of the differences we encounter in the general case is that the \ord{$n$} layer of the Taylor tower is no longer determined by a single spectrum with $\Sigma_n$-action. Our definition of derivatives (given in Section~\ref{sec:der}) is therefore necessarily more involved. For us the \ord{$n$} derivative of a functor $F: \cat{C} \to \cat{D}$ is a diagram of spectra of the form
\[ \der_nF: \spectra(\cat{C})^n \times \spectra(\cat{D})^{op} \to \spectra \]
that is symmetric in the $n$ copies of $\spectra(\cat{C})$, and is linear in each variable. Here $\spectra(\cat{C})$ denotes the stabilization of the $\infty$-category $\cat{C}$ as described by Lurie in \cite[1.4]{lurie:2017}.

When $\cat{C}$ and $\cat{D}$ are each the $\infty$-category of pointed spaces or spectra, these stabilizations are both $\spectra$, the $\infty$-category of spectra. If $F$ preserves filtered colimits, the resulting symmetric multilinear functor $\spectra^n \times \spectra^{op} \to \spectra$ is determined by its value on the sphere spectrum in each variable. This value recovers the spectrum with $\Sigma_n$-action that is usually referred to as the \ord{$n$} derivative of the functor $F$. To simplify this introduction we suppress the dependence of the derivative on other variables in what follows. More explicit statements in the case of general $\cat{C}$ and $\cat{D}$ can be found in the main body of the paper.

Our philosophy is to start by focusing on functors $F: \cat{C} \to \spectra$. Let $\cat{F}_{\cat{C}}$ denote the $\infty$-category of those functors of this type that are reduced (i.e. $F(*) \homeq *$) and finitary (preserve filtered colimits). The construction of the \ord{$n$} derivative can then be viewed as a functor
\[ \der_n : \cat{F}_{\cat{C}} \to \spectra. \]
Now the $\infty$-category $\cat{F}_{\cat{C}}$ has a (non-unital) symmetric monoidal product given by the objectwise smash product of functors, and therefore the category $\Fun(\cat{F}_{\cat{C}},\spectra)$ of functors $\cat{F}_{\cat{C}} \to \spectra$ has a symmetric monoidal product $\otimes$ given by the \emph{Day convolution} of the objectwise smash product on $\cat{F}_{\cat{C}}$ and the ordinary smash product on $\smsh$.

Our first main theorem (proved in Section~\ref{sec:conv}) gives a relationship between the functors $\der_n$, for different $n$, in terms of this Day convolution structure.

\begin{theorem} \label{thm:A}
Let $\cat{C}$ be a pointed compactly-generated $\infty$-category. Then there is a $\Sigma_n$-equivariant equivalence, in the $\infty$-category $\Fun(\cat{F}_{\cat{C}},\spectra)$, of the form
\[ \der_n \homeq \der_1^{\otimes n}. \]
\end{theorem}

Next we turn to (reduced, finitary) functors $F: \cat{C} \to \cat{D}$ between two arbitrary pointed compactly-generated $\infty$-categories. Our second main theorem (proved in Section~\ref{sec:generalF}) allows us to identify the derivatives of such a functor $F$ in terms of the derivatives of spectrum-valued functors on $\cat{C}$ and $\cat{D}$.

\begin{theorem} \label{thm:B}
Let $F: \cat{C} \to \cat{D}$ be a reduced functor that preserves filtered colimits. Then there is a natural equivalence
\[ \der_nF \homeq \Nat(\der_1(-), \der_n(- \circ F)) \]
where the right-hand side is the spectrum of natural transformations between two functors of type $\cat{F}_{\cat{D}} \to \spectra$.
\end{theorem}

Combining Theorems~\ref{thm:A} and \ref{thm:B}, we get new models for $\der_nF$ that can be defined entirely in terms of the first derivative construction for spectrum-valued functors, and Day convolution:
\begin{equation} \label{eq:der} \der_nF \homeq \Nat(\der_1(-),\der_1^{\otimes n}(- \circ F)). \end{equation}
One unanswered question of \cite{arone/ching:2011} was whether such models can admit (unital and associative) composition maps of the form
\begin{equation} \label{eq:comp} \der_*(G) \circ \der_*(F) \to \der_*(GF) \end{equation}
which, in particular, provide the derivatives of the identity functor (or, indeed, any monad) with an operad structure. The models given in (\ref{eq:der}) permit the construction of composition maps of the form (\ref{eq:comp}) for a wide range of $\infty$-categories. (In the case of pointed spaces and spectra a simpler approach is due to Yeakel~\cite{yeakel:2020}.)

When $F$ is the identity functor $I_{\cat{C}}$ on a pointed compactly-generated $\infty$-category $\cat{C}$, equation (\ref{eq:der}) takes the form
\[ \der_nI_{\cat{C}} \homeq \Nat(\der_1,\der_1^{\otimes n}). \]
The symmetric sequence $\der_*I_\cat{C}$ has an operad structure given by composition of natural transformations, a so-called `coendomorphism operad' for the object $\der_1$ with respect to Day convolution. In a similar way, for $F: \cat{C} \to \cat{D}$, the derivatives of $F$ form a bimodule over the operads $\der_*I_{\cat{C}}$ and $\der_*I_{\cat{D}}$.

To be more precise, what we get are $\infty$-operads (in the sense of Lurie~\cite[2.1]{lurie:2017}) and bimodules over those $\infty$-operads. We give explicit constructions of these objects in Sections~\ref{sec:operad} and~\ref{sec:bimodule}. Those constructions rely heavily on some technical constructions with symmetric monoidal $\infty$-categories: work of Lurie~\cite[2.2.6]{lurie:2017} (on Day convolution) and of Barwick-Glasman-Nardin~\cite{barwick/glasman/nardin:2018} (on opposite symmetric monoidal structures). Combining these two pieces of work, we get a (non-unital) symmetric monoidal $\infty$-category
\[ \Fun(\cat{F}_{\cat{C}},\spectra)^{op,\otimes} \]
whose underlying $\infty$-category is the opposite of the $\infty$-category of functors $\cat{F}_{\cat{C}} \to \spectra$. (In fact, we have to take care over size issues at this point, and replace $\cat{F}_{\cat{C}}$ with a small symmetric monoidal subcategory, but we ignore that issue for the remainder of this introduction.)

Taking the suboperad of $\Fun(\cat{F}_{\cat{C}},\spectra)^{op,\otimes}$ generated by $\der_1$ then produces an $\infty$-operad $\mathbb{I}_{\cat{C}}^{\otimes}$ that encodes the coendomorphism operad structure on $\der_*I_{\cat{C}}$ described above.

In Section~\ref{sec:ex-AlgO}, we focus on one principal example of our general theory: the case where $\cat{C}$ is the $\infty$-category of (non-unital) stable algebras over a stable $\infty$-operad $\cat{O}^{\otimes}$. This example includes $\infty$-categories of structured ring spectra, such as $E_n$-ring spectra and spectral Lie algebras. Our general result is a calculation of the $\infty$-operad $\mathbb{I}_{\cat{C}}^{\otimes}$ in this case:

\begin{theorem} \label{thm:alg}
Let $\cat{O}^{\otimes}$ be a stable $\infty$-operad and let $\Alg_{\cat{O}^{\otimes}}$ be the $\infty$-category of non-unital stable $\cat{O}^{\otimes}$-algebras. Then $\mathbb{I}^{\otimes}_{\Alg_{\cat{O}^{\otimes}}}$ is (Morita-)equivalent to $\cat{O}^{\otimes}$ itself.
\end{theorem}

By a Morita equivalence of stable $\infty$-operads, we mean an equivalence between the $\infty$-categories of stable algebras over those $\infty$-operads. We actually show a stronger result: that $\mathbb{I}^{\otimes}_{\Alg_{\cat{O}^{\otimes}}}$ contains $\cat{O}^{\otimes}$ as a full suboperad, and it is the inclusion of that suboperad that induces an equivalence between $\infty$-categories of algebras.

Theorem~\ref{thm:alg} verifies a longstanding principle in Goodwillie calculus: that the derivatives of the identity functor on a category of operadic algebras recovers the original operad. For example, this principle can be seen on an arity-wise basis in the work of Harper and Hess~\cite[1.14]{harper/hess:2013}. Here we promote that principle to a full equivalence of ($\infty-$)operads. An alternative approach to this calculation was recently introduced by Clark~\cite{clark:2020}.

In Section~\ref{sec:bimodule} we turn back to the derivatives of an arbitrary (reduced, finitary) functor $F: \cat{C} \to \cat{D}$ between pointed compactly-generated $\infty$-categories, and provide a precise construction of the derivatives of $F$ as a bimodule over the $\infty$-operads $\der_*I_{\cat{C}}$ and $\der_*I_{\cat{D}}$. By a \emph{bimodule} over two $\infty$-operads, we mean a $\Delta^1$-family of $\infty$-operads (in the sense of Lurie~\cite[2.3.2.10]{lurie:2017}) that restricts to the given $\infty$-operads over the endpoints of $\Delta^1$.

For a pair of functors $F: \cat{C} \to \cat{D}$ and $G: \cat{D} \to \cat{E}$, our construction gives rise to a `composition' map
\[ \mathbb{D}_G^{\otimes} \circ_{\mathbb{I}_{\cat{D}}^{\otimes}} \mathbb{D}_{F}^{\otimes} \to \mathbb{D}_{GF}^{\otimes} \]
We conjecture a Chain Rule, generalizing that of~\cite{arone/ching:2011}, which says that the above map is an equivalence of $(\mathbb{I}_{\cat{E}}^{\otimes},\mathbb{I}_{\cat{C}}^{\otimes})$-bimodules.

\subsection*{Technical background}

We use $\infty$-categories, or quasi-categories, as our basic model for $(\infty,1)$-categories, yet very little technical knowledge of this theory is required in sections 1-3 of the paper. Our two main results about Goodwillie calculus depend only on basic homotopy theory such as properties of homotopy limits and colimits. These results could be stated, and proved, in more-or-less exactly the same way in the context of simplicial model categories instead.

In later sections, the theory of $\infty$-categories, and in particular that of $\infty$-operads, as developed by Lurie, plays a much more concerted role. We rely heavily on \cite{lurie:2009} and \cite{lurie:2017} as references, though we do recall the basic principles of the theory of $\infty$-operads in Section~\ref{sec:operad}. We require two particularly technical constructions on symmetric monoidal $\infty$-categories: the Day convolution structure of~\cite[2.2.6]{lurie:2017}; and the opposite monoidal $\infty$-category construction due to Barwick-Glasman-Nardin~\cite{barwick/glasman/nardin:2018}.

For the initial development of Goodwillie calculus in the context of $\infty$-categories, we rely on \cite[Ch. 6]{lurie:2017}, though the reader will not need any of the technical details of that work.


\subsection*{Notation}

We use letters such as $\cat{C},\cat{D}$ to stand for pointed compactly-generated $\infty$-categories, and the symbol $*$ to denote a null object in such. In particular, we have $\based$, the $\infty$-category of pointed (small) Kan complexes, and $\spectra$, the $\infty$-category of spectra from \cite[1.4.3]{lurie:2017}. We also make use of the standard adjunction
\[ \Sigma^\infty : \based \rightleftarrows \spectra : \Omega^\infty. \]

For a pointed $\infty$-category $\cat{C}$, we write $\Hom_{\cat{C}}(-,-)$ for (some model of) the pointed simplicial set of maps between two objects of $\cat{C}$. If $\cat{C}$ is stable, it admits mapping spectra which we denote $\Map_{\cat{C}}(-,-)$, so that $\Hom_{\cat{C}}(-,-) \homeq \Omega^\infty \Map_{\cat{C}}(-,-)$.

A pointed compactly-generated $\infty$-category $\cat{C}$ admits tensors by pointed simplicial sets, which we write with a smash product symbol. In particular, we have suspensions $\Sigma^L x \homeq S^L \smsh x$ for $x \in \cat{C}$. Note that we then have natural equivalences of pointed simplicial sets
\[ \Omega^L\Hom_{\cat{C}}(x,-) \homeq \Hom_{\cat{C}}(\Sigma^L x, -). \]

We often consider the $\infty$-category of functors between two other $\infty$-categories, which we denote in the form $\Fun(\cat{C},\cat{D})$. When $\cat{D} = \spectra$, the $\infty$-category $\Fun(\cat{C},\cat{D})$ is stable in which case we write
\[ \Nat_{\cat{C}}(-,-) := \Map_{\Fun(\cat{C},\spectra)}(-,-). \]
This plays the role of a spectrum of natural transformations between two such functors.

We omit notation for the nerve construction: for example, $\Surj$ is the nerve of the category of finite pointed sets and pointed surjections.

When we say limit or colimit, we almost always mean \emph{homotopy} limit or colimit (and we denote these as $\holim$ or $\hocolim$). The exception is when constructing an $\infty$-category, for example as a pullback of other $\infty$-categories, in which case we intend a strict pullback in the category of simplicial sets.

\subsection*{Acknowledgements}

None of the work in this paper would exist without considerable support from, and helpful conversations with, Greg Arone. Feedback from Jacob Lurie was useful for understanding stable $\infty$-operads and the Chain Rule. Questions and comments from Daniel Fuentes-Keuthan greatly improved the exposition, as did many suggestions and corrections from an anonymous referee which led to improvements to Section~\ref{sec:ex-AlgO} in particular.

The research in this paper is supported by the National Science Foundation through grant DMS-1709032. Some of this work was undertaken during the programme \emph{Homotopy Harnessing Higher Structures} in 2018 at the Isaac Newton Institute for Mathematical Sciences, Cambridge, supported by EPSRC grant number EP/R014604/1.

This is the accepted version of an article to be published in the Journal of the London Mathematical Society.

\section{Goodwillie derivatives in \texorpdfstring{$\infty$}{infinity}-categories} \label{sec:der}

Let $F: \cat{C} \to \cat{D}$ be a reduced functor between pointed compactly-generated $\infty$-categories. Such $F$ has a \emph{Taylor tower} constructed in this generality by Lurie~\cite[6.1]{lurie:2017} following Goodwillie's original approach \cite{goodwillie:2003}. This tower is a sequence of functors of the form
\[ F \to \dots \to P_nF \to P_{n-1}F \to \dots \to P_1F \to P_0F \homeq * \]
where $F \to P_nF$ is initial (up to homotopy) among natural transformations from $F$ to an $n$-excisive functor. The \emph{\ord{$n$} layer} in the Taylor tower is the fibre
\[ D_nF := \hofib(P_nF \to P_{n-1}F) \]
and $D_nF: \cat{C} \to \cat{D}$ is an $n$-homogenous functor.

One of Goodwillie's main results provides a classification of homogeneous functors, which shows that the \ord{$n$} layer $D_nF$ can be recovered from a symmetric multilinear functor $\Delta_nF: \spectra(\cat{C})^n \to \spectra(\cat{D})$ by the formula
\[ D_nF(X) \homeq \Omega^\infty_{\cat{D}}[\Delta_nF(\Sigma^\infty_{\cat{C}}X,\dots,\Sigma^\infty_{\cat{C}}X)_{h\Sigma_n}] \]
where
\[ \Sigma^\infty_{\cat{C}} : \cat{C} \rightleftarrows \spectra(\cat{C}) : \Omega^\infty_{\cat{C}} \]
is the stabilization adjunction for $\cat{C}$; see \cite[6.1.4.7 and 6.2.3.22]{lurie:2017}.

\begin{definition} \label{def:deriv}
Let $F: \cat{C} \to \cat{D}$ be a reduced functor between pointed compactly-generated $\infty$-categories, and let $\Delta_nF: \spectra(\cat{C})^n \to \spectra(\cat{D})$ be the symmetric multilinear functor described above. The \emph{\ord{$n$} derivative of $F$} is the functor
\[ \der_nF: \spectra(\cat{C})^n \times \spectra(\cat{D})^{op} \to \spectra \]
defined by
\[ \der_nF(X_1,\dots,X_n;Y) := \Map_{\spectra(\cat{D})}(Y,\Delta_nF(X_1,\dots,X_n)) \]
where $\Map_{\spectra(\cat{D})}(-,-)$ denotes a mapping spectrum construction for the stable $\infty$-category $\spectra(\cat{D})$. In other words, we can think of $\der_nF$ as the composite of $\Delta_nF$ with the stable Yoneda embedding for the stable $\infty$-category $\spectra(\cat{D})$.

Note that $\der_nF$ is symmetric multilinear in the $\spectra(\cat{C})$ variables, and preserves all limits in $\spectra(\cat{D})^{op}$ (that is, takes colimits in $\spectra(\cat{D})$ to limits in $\spectra$).
\end{definition}

\begin{example} \label{ex:topsp}
When $\cat{C}$ and $\cat{D}$ are both either $\based$ or $\spectra$, and $F$ preserves filtered colimits, the functor $\der_nF$ of Definition~\ref{def:deriv} is determined by the single spectrum (with $\Sigma_n$-action)
\[ \der_nF(S^0,\dots,S^0;S^0) \]
where $S^0$ is the sphere spectrum. We write $\der_nF$ also for this individual spectrum, which is the object typically referred to as \emph{the} \ord{$n$} derivative of $F$ in this case.
\end{example}

\begin{example} \label{ex:Csp}
When $\cat{D}$ is $\based$ or $\spectra$, there is an equivalence
\[ \der_nF(X_1,\dots,X_n;S^0) \homeq \Delta_nF(X_1,\dots,X_n). \]
We also write this object as $\der_nF(X_1,\dots,X_n)$. More generally, whenever either $\cat{C}$ or $\cat{D}$ is $\based$ or $\spectra$, we may omit the corresponding arguments of the functor $\der_nF$, in which case those arguments are assumed to be the sphere spectrum $S^0$.
\end{example}

\section{Derivatives of spectrum-valued functors} \label{sec:conv}

We now turn to our first main result, which concerns the derivatives of spectrum-valued functors.

\begin{definition}
Fix a pointed compactly-generated $\infty$-category $\cat{C}$ and let $\cat{F}_{\cat{C}}$ be the full subcategory of $\Fun(\cat{C},\spectra)$ whose objects are the reduced functors $\cat{C} \to \spectra$ that preserve filtered colimits. For objects $X_1,\dots,X_n \in \spectra(\cat{C})$, Example~\ref{ex:Csp} says that we have a functor
\[ \der_n(-)(X_1,\dots,X_n) : \cat{F}_{\cat{C}} \to \spectra. \]
The goal of this section is to understand how these functors are related to one another for varying $n$.
\end{definition}

The relationship we are looking for is via a version of Day convolution (see \cite{day:1970}) for functors $\cat{F}_{\cat{C}} \to \spectra$ with respect to the following symmetric monoidal structures: on $\cat{F}_{\cat{C}}$ the objectwise smash product of functors; and on $\spectra$ the ordinary smash product. Later in the paper, we work with a symmetric monoidal $\infty$-category that represents this Day convolution, but for now it is sufficient to describe convolution by its universal property.

\begin{definition} \label{def:Day}
The \emph{Day convolution} of ${A},{B}: \cat{F}_{\cat{C}} \to \spectra$, if it exists, consists of a functor
\[ {A} \otimes {B} : \cat{F}_{\cat{C}} \to \spectra \]
and a natural transformation of functors $\cat{F}_{\cat{C}} \times \cat{F}_{\cat{C}} \to \spectra$ of the form
\[ \alpha: {A}(-) \smsh {B}(-) \to ({A} \otimes {B})(- \smsh -) \]
that induces equivalences of mapping spaces
\[ \Hom_{\Fun(\cat{F}_\cat{C},\spectra)}({A} \otimes {B}, {C}) \weq \Hom_{\Fun(\cat{F}_{\cat{C}} \times \cat{F}_{\cat{C}},\spectra)}({A}(-) \smsh {B}(-), {C}(- \smsh -)) \]
for an arbitrary functor ${C}: \cat{F}_{\cat{C}} \to \spectra$. Note that we use the symbol $\smsh$ on the right-hand side to denote both the smash product of spectra and the objectwise smash product on $\cat{F}_{\cat{C}}$. We define convolution of more than two functors in a similar way.
\end{definition}

\begin{remark} \label{rem:Day}
Definition~\ref{def:Day} says that Day convolution is a left Kan extension and implies that it is unique up to equivalence. In the cases we care about, we will prove existence directly, primarily via Lemma~\ref{lem:Day} below.
\end{remark}

The main result of this section is the following relationship between the \ord{$n$} and $1^{\text{st}}$ derivative constructions for functors from $\cat{C}$ to $\spectra$.

\begin{theorem} \label{thm:dn}
Let $\cat{C}$ be a pointed compactly-generated $\infty$-category, and consider objects $X_1,\dots,X_n \in \spectra(\cat{C})$. Then there is a natural equivalence
\[ \der_n(-)(X_1,\dots,X_n) \homeq \der_1(-)(X_1) \otimes \dots \otimes \der_1(-)(X_n) \]
where $\der_n$ denotes the \ord{$n$} derivative construction for functors $\cat{C} \to \spectra$, and $\otimes$ denotes the Day convolution of Definition~\ref{def:Day}.
\end{theorem}

\begin{corollary}
When $\cat{C} = \based$ or $\spectra$, taking $X_1 = \dots = X_n = S^0$ in Theorem~\ref{thm:dn} gives the formula
\[ \der_n \homeq \der_1^{\otimes n}. \]
\end{corollary}

\begin{remark}
The Day convolution cannot be calculated objectwise. In particular, Theorem~\ref{thm:dn} does \emph{not} imply that the \ord{$n$} derivative of a \emph{particular} functor $F: \cat{C} \to \spectra$ can be calculated from the first derivative of $F$ (which would clearly be false). Rather it says that $\der_nF$ can be calculated as a homotopy colimit of the form
\[ \der_nF \homeq \hocolim_{G_1 \smsh \dots \smsh G_n \to F} \der_1G_1 \smsh \dots \smsh \der_1G_n \]
calculated over the $\infty$-category of $n$-tuples of functors $G_1,\dots,G_n$ with a map $G_1 \smsh \dots \smsh G_n \to F$.
\end{remark}

\begin{remark} \label{rem:I_C}
Since Day convolution is, in fact, a symmetric monoidal structure, Theorem~\ref{thm:dn} allows us to see that the collection of functors $(\der_n)_{n \geq 1}$ possesses additional structure. Suppose we define a coloured operad $\mathbb{I}_{\cat{C}}$, enriched in $\spectra$, with colours given by the objects of $\spectra(\cat{C})$ and terms
\begin{equation} \mathbb{I}_{\cat{C}}(X_1,\dots,X_n;Y) = \Nat_{\cat{F}_{\cat{C}}}(\der_1(-)(Y), \der_1(-)(X_1) \otimes \dots \otimes \der_1(-)(X_n)). \label{eq:P} \end{equation}
where $\Nat_{\cat{F}_{\cat{C}}}(-,-)$ denotes a mapping \emph{spectrum} construction for the stable $\infty$-category $\Fun(\cat{F}_{\cat{C}},\spectra)$. The operad structure is given by composition of natural transformations. It is an easy consequence of Theorem~\ref{thm:dn} that the derivatives of any functor $\cat{C} \to \spectra$ form a right module over the operad $\mathbb{I}_{\cat{C}}$. As stated, these operad and module structures are only associative up to homotopy; a more precise definition of $\mathbb{I}_{\cat{C}}$ as an $\infty$-operad is given in Section~\ref{sec:operad}.
\end{remark}

The remainder of this section consists of the proof of Theorem~\ref{thm:dn}. This proof relies largely on Goodwillie's identification of the derivative as a multilinearized cross-effect. That is, we have, for $x_1,\dots,x_n \in \cat{C}$ and $F: \cat{C} \to \spectra$:
\begin{equation} \label{eq:goodwillie} \Delta_nF(\Sigma^\infty_{\cat{C}}x_1,\dots,\Sigma^\infty_{\cat{C}}x_n) \homeq \hocolim_{L \to \infty} \Omega^{nL} \creff_nF(\Sigma^L x_1,\dots,\Sigma^L x_n). \end{equation}
where $\Delta_nF$ is the symmetric multilinear functor that classifies the homogeneous functor $D_nF$. (An $\infty$-categorical version of this result follows from \cite[6.1.3.23 and 6.1.1.28]{lurie:2017}.)

We also use the fact, extending \cite[3.13]{arone/ching:2016}, that cross-effects of spectrum-valued functors can be represented as natural transformation objects. To see this fact we first need a version of the Yoneda Lemma in this context.

\begin{lemma} \label{lem:yoneda}
Let $\cat{C}$ be a pointed $\infty$-category, $x$ an object of $\cat{C}$, and $F: \cat{C} \to \spectra$ a reduced functor. Then there is a natural equivalence of spectra
\[ \Nat_{\cat{C}}(\Sigma^\infty \Hom_{\cat{C}}(x,-), F(-)) \homeq F(x) \]
where $\Nat_{\cat{C}}(-,-)$ denotes a mapping spectrum for the stable $\infty$-category $\Fun(\cat{C},\spectra)$.
\end{lemma}
\begin{proof}
Any functor $F: \cat{C} \to \spectra$ admits a natural map
\[ \Hom_{\cat{C}}(x,-) \to \Hom_{\spectra}(F(x),F(-)) \homeq \Omega^\infty \Map_{\spectra}(F(x),F(-)) \]
which is basepoint-preserving when $F$ is reduced and which corresponds by adjunction to the desired map
\[ F(x) \to \Nat_{\cat{C}}(\Sigma^\infty \Hom_{\cat{C}}(x,-),F(-)). \]
To prove this map is an equivalence of spectra, it is sufficient to show that the induced map
\[ \Omega^\infty \Sigma^{-k} F(x) \to \Omega^\infty \Sigma^{-k} \Nat_{\cat{C}}(\Sigma^\infty \Hom_{\cat{C}}(x,-),F(-)) \]
is an equivalence of simplicial sets for all $k \in \mathbb{Z}$. We can identify the right-hand side with
\[ \Hom_{\Fun(\cat{C},\spectra)}(\Hom_{\cat{C}}(x,-),\Omega^\infty \Sigma^{-k}F(-)), \]
and the claim follows from the ordinary Yoneda Lemma.
\end{proof}

We then have the following description of the cross-effects of a functor $F: \cat{C} \to \spectra$.

\begin{lemma} \label{lem:creff}
Let $\cat{C}$ be a pointed $\infty$-category. For any $F: \cat{C} \to \spectra$ and objects $x_1,\dots,x_n \in \cat{C}$, we have a natural equivalence
\[ \creff_nF(x_1,\dots,x_n) \homeq \Nat_{\cat{C}}(\Sigma^\infty \Hom_{\cat{C}}(x_1,-) \smsh \dots \smsh \Sigma^\infty \Hom_{\cat{C}}(x_n,-), F(-)). \]
\end{lemma}
\begin{proof}
The first cross-effect is the fibre
\[ \creff_1F(x) := \hofib(F(x) \to F(*)). \]
When $G: \cat{C} \to \spectra$ is reduced, the map
\[ \Nat_{\cat{C}}(G,\creff_1F) \to \Nat_{\cat{C}}(G,F) \]
has inverse given by
\[ \Nat_{\cat{C}}(G,F) \to \Nat_{\cat{C}}(\creff_1G,\creff_1F) \homeq \Nat_{\cat{C}}(G,\creff_1F). \]
Since we also have an equivalence $\creff_n(\creff_1F) \weq \creff_nF$, we can replace $F$ with $\creff_1F$ in the lemma, i.e. we may assume that $F$ is reduced. 

The case $n = 1$ now follows from Lemma~\ref{lem:yoneda}. We describe the case $n = 2$. The general case is virtually identical.

The \ord{$n$} cross-effect is defined as the total fibre of an $n$-cube (see \cite{goodwillie:1991}); for $n = 2$, this cube takes the form
\[ \creff_2F(x_1,x_2) \homeq \thofib \left( \begin{diagram} \node{F(x_1 \wdge x_2)} \arrow{e} \arrow{s} \node{F(x_1)} \arrow{s} \\ \node{F(x_2)} \arrow{e} \node{F(*)} \end{diagram} \right). \]

Using Lemma~\ref{lem:yoneda} we can write this square as
\[ \begin{diagram} \dgARROWLENGTH=2em
  \node{\Nat_{\cat{C}}(\Sigma^\infty \Hom_{\cat{C}} \left( x_1 \wdge x_2, -\right), F(-))} \arrow{e} \arrow{s}
    \node{\Nat_{\cat{C}}(\Sigma^\infty \Hom_{\cat{C}} \left(x_1, -\right), F(-))} \arrow{s} \\
  \node{\Nat_{\cat{C}}(\Sigma^\infty \Hom_{\cat{C}} \left( x_2, -\right), F(-))} \arrow{e}
    \node{\Nat_{\cat{C}}(\Sigma^\infty \Hom_{\cat{C}} \left( *, -\right), F(-)).}
\end{diagram} \]
Since $\Nat_{\cat{C}}(\Sigma^\infty -,F)$ takes colimits (of $\based$-valued functors on $\cat{C}$) to limits (of spectra), we therefore have
\begin{equation} \label{eq:creff2} \creff_2F(x_1,x_2) \homeq \Nat_{\cat{C}}\left(\Sigma^\infty A(-), F(-)\right) \end{equation}
where $A(-)$ is the total \emph{cofibre} of the $2$-cube of spaces of the form
\[ \begin{diagram} \dgARROWLENGTH=2em
  \node{\Hom_{\cat{C}}(*,-)} \arrow{e} \arrow{s} \node{\Hom_{\cat{C}}(x_1,-)} \arrow{s} \\
  \node{\Hom_{\cat{C}}(x_2,-)} \arrow{e} \node{\Hom_{\cat{C}}(x_1 \wdge x_2,-).}
\end{diagram} \]
That square can be written in the form
\[ \begin{diagram} \dgARROWLENGTH=2em
  \node{*} \arrow{e} \arrow{s} \node{\Hom_{\cat{C}}(x_1,-)} \arrow{s} \\
  \node{\Hom_{\cat{C}}(x_2,-)} \arrow{e} \node{\Hom_{\cat{C}}(x_1,-) \times \Hom_{\cat{C}}(x_2,-).}
\end{diagram} \]
But then the total cofibre $A(-)$ is equivalent to the smash product
\[ \Hom_{\cat{C}}(x_1,-) \smsh \Hom_{\cat{C}}(x_2,-) \]
which, together with (\ref{eq:creff2}), provides the desired equivalence. For the case of general $n$, the key observation is that the total cofibre of an $n$-cube of pointed spaces of the form
\[ \left\{ \prod_{i \in S} A_i \right\}_{S \subseteq \{1,\dots,n\}} \]
is equivalent to the smash product $A_1 \smsh \dots \smsh A_n$.
\end{proof}

It follows from Lemma~\ref{lem:creff} that the terms appearing in the homotopy colimit of (\ref{eq:goodwillie}) can also be expressed in terms of natural transformation objects:
\begin{equation} \label{eq:creff} \begin{split}
    \Omega^{nL} \creff_n(F)(\Sigma^L x_1,\dots,\Sigma^L x_n)
    &\homeq \Omega^{nL} \Nat_{\cat{C}}\left( \Smsh_{i = 1}^{n} \Sigma^\infty \Hom_{\cat{C}}(\Sigma^L x_i,-),F \right) \\
    &\homeq \Nat_{\cat{C}}\left( \Smsh_{i = 1}^{n} \Sigma^\infty \Sigma^L \Omega^L \Hom_{\cat{C}}(x_i,-),F \right)
\end{split} \end{equation}
where the first equivalence is that of Lemma~\ref{lem:creff}, and the second is built from several instances of the adjunction $(\Sigma^L,\Omega^L)$.

It remains to identify how these equivalences interact with the maps in the colimit in (\ref{eq:goodwillie}).

\begin{lemma} \label{lem:calc}
For reduced $F: \cat{C} \to \spectra$ and objects $x_1,\dots,x_n$, the following diagram of spectra commutes up to equivalence:
\[ \begin{diagram} \dgARROWLENGTH=1em
  \node{\Omega^{nL}\creff_nF(\Sigma^L x_1,\dots,\Sigma^L x_n)} \arrow{e,t}{t_{1,\dots,1}(\creff_nF)} \arrow{s,l}{\sim} \node{\Omega^{n(L+1)}\creff_nF(\Sigma^{L+1}x_1,\dots,\Sigma^{L+1}x_n)}  \arrow{s,r}{\sim} \\
  \node{\Nat_{\cat{C}}\left( \Smsh_{i = 1}^{n} \Sigma^\infty \Sigma^L \Omega^L \Hom_{\cat{C}}(x_i,-),F \right)} \arrow{e,t}{\epsilon^*} \node{\Nat_{\cat{C}}\left( \Smsh_{i = 1}^{n} \Sigma^\infty \Sigma^{L+1} \Omega^{L+1} \Hom_{\cat{C}}(x_i,-),F \right)}
\end{diagram} \]
where the vertical maps are the equivalences of (\ref{eq:creff}), the top horizontal map is (the multivariable version of) the stabilization map $t_1$ appearing in Goodwillie's construction of the linearization of a functor (see \cite[1.10]{goodwillie:1990} or \cite[6.1.1.27]{lurie:2017}), and the bottom horizontal map is that induced by the counit
\[ \epsilon: \Sigma^{L+1} \Omega^{L+1} = \Sigma^L(\Sigma \Omega)\Omega^L \to \Sigma^L \Omega^L. \]
\end{lemma}
\begin{proof}
We illustrate with the case $n = 1$. The general case is similar. Since $\creff_1F \homeq F$ for $F$ reduced, our diagram takes the form
\begin{equation} \label{eq:calc} \begin{diagram} \dgARROWLENGTH=1.5em
  \node{\Omega^{L}F\Sigma^L x_1} \arrow{e,t}{t_1F} \arrow{s,l}{\sim} \node{\Omega^{L+1}F\Sigma^{L+1}x_1} \arrow{s,r}{\sim} \\
  \node{\Omega^{L}\Nat_{\cat{C}}(\Sigma^\infty \Hom_{\cat{C}}(\Sigma^L x_1,-),F)} \arrow{e,t}{\epsilon^*} \arrow{s,l}{\sim} \node{\Omega^L\Nat_{\cat{C}}(\Sigma^\infty \Sigma \Omega \Hom_{\cat{C}}(\Sigma^L x_1,-),F)} \arrow{s,r}{\sim} \\
  \node{\Nat_{\cat{C}}(\Sigma^\infty \Sigma^L \Omega^L \Hom_{\cat{C}}(x_1,-),F(-))} \arrow{e,t}{\epsilon^*} \node{\Nat_{\cat{C}}(\Sigma^\infty \Sigma^{L+1} \Omega^{L+1} \Hom_{\cat{C}}(x_1,-),F(-))}
\end{diagram} \end{equation}
where the bottom square commutes by naturality of $\epsilon$, and the top square is ($\Omega^L$ applied to) the $L = 0$ case of the Lemma, with $x_1$ replaced by $\Sigma^L x_1$. It thus is sufficient to consider the case $L = 0$.

To do this, first recall how the map $t_1F$ is constructed. Let $\mathbf{D}$ denote the diagram
\[ \begin{diagram} \dgARROWLENGTH=1em
  \node[2]{I} \arrow{s} \\
  \node{I} \arrow{e} \node{S^1}
\end{diagram} \]
where $I$ is a closed interval with one of its endpoints as the basepoint, and the two maps are the inclusions into the two halves of the circle $S^1$.

Two copies of the inclusion $S^0 \to I$ form a cone over $\mathbf{D}$ and induce the horizontal (and diagonal) maps in the following commutative diagram:
\begin{equation} \label{eq:D} \begin{diagram} \dgARROWLENGTH=1em
  \node{F(x_1)} \arrow{e} \arrow{s,lr}{\sim}{Y_{x_1}} \node{\holim_{A \in \mathbf{D}}F(A \smsh x_1)} \arrow{s,lr}{\sim}{\holim Y_{A \smsh x_1}} \\
  \node{\Nat_{\cat{C}}(\Sigma^\infty \Hom_{\cat{C}}(x_1,-),F)} \arrow{e} \arrow{se} \node{\holim_{A \in \mathbf{D}}\Nat_{\cat{C}}(\Sigma^\infty \Hom_{\cat{C}}(A \smsh x_1,-),F)} \arrow{s,r}{\sim} \\
  \node[2]{\Nat_{\cat{C}}(\Sigma^\infty \hocolim_{A \in \mathbf{D}} \Hom_{\cat{C}}(x_1,-)^A, F)}
\end{diagram} \end{equation}
where $Y_x$ denotes the stable Yoneda embedding from Lemma~\ref{lem:yoneda} applied with object $x$, and the right-hand bottom vertical map is a canonical equivalence involving the tensoring adjunction for objects in $\cat{C}$.

We now argue that the right-hand column of (\ref{eq:D}) can be identified with the right-hand column in (\ref{eq:calc}) by observing that since $I$ is contractible, each homotopy limit/colimit is a loop-space/suspension. That is, we have the following commutative diagram in which the horizontal maps are induced by the equivalences $* \weq I$:
\[ \begin{diagram} \dgARROWLENGTH=1.5em
  \node{\Omega F \Sigma x_1} \arrow{s,lr}{\sim}{\Omega Y_{\Sigma x_1}} \arrow{e,t}{\sim} \node{\holim_{A \in \mathbf{D}}F(A \smsh x_1)} \arrow{s,lr}{\sim}{\holim Y_{A \smsh x_1}} \\
  \node{\Omega \Nat_{\cat{C}}(\Sigma^\infty \Hom_{\cat{C}}(\Sigma x_1,-),F)} \arrow{s,r}{\sim} \arrow{e,t}{\sim} \node{\holim_{A \in \mathbf{D}}\Nat_{\cat{C}}(\Sigma^\infty \Hom_{\cat{C}}(A \smsh x_1,-),F)} \arrow{s,r}{\sim} \\
  \node{\Nat_{\cat{C}}(\Sigma^\infty \Sigma\Omega \Hom_{\cat{C}}(x_1,-),F)} \arrow{e,t}{\sim} \node{\Nat_{\cat{C}}(\Sigma^\infty \hocolim_{A \in \mathbf{D}} \Hom_{\cat{C}}(x_1,-)^A, F)}
\end{diagram} \]
The top map of (\ref{eq:D}) is precisely $t_1F$, and it is easy to identify the bottom map with $\epsilon^*$ when combined with the bottom map in the diagram above.
\end{proof}

Taking homotopy colimits as $L \to \infty$ over the diagrams in Lemma~\ref{lem:calc}, and applying (\ref{eq:goodwillie}), we get the following result.

\begin{proposition} \label{prop:der_n}
For reduced $F: \cat{C} \to \spectra$ and objects $x_1,\dots,x_n \in \cat{C}$, there is an equivalence
\[ \Delta_nF(\Sigma^\infty_{\cat{C}}x_1,\dots,\Sigma^\infty_{\cat{C}}x_n) \homeq \hocolim_{L \to \infty} \Nat_{\cat{C}}\left(\Smsh_{i = 1}^{n} \Sigma^\infty \Sigma^L \Omega^L \Hom_{\cat{C}}(x_i,-), F \right) \]
where the maps in the homotopy colimit are induced by the counit map $\epsilon: \Sigma\Omega \to I$.
\end{proposition}

We also require the following result about Day convolution of representable functors.

\begin{lemma} \label{lem:Day}
For $F_1,\dots,F_n \in \cat{F}_{\cat{C}}$, we have an equivalence:
\[ \Nat_{\cat{C}}(F_1 \smsh \dots \smsh F_n,-) \homeq \Nat_{\cat{C}}(F_1,-) \otimes \dots \otimes \Nat_{\cat{C}}(F_n,-). \]
\end{lemma}
\begin{proof}
We describe the case $n = 2$. The general case is virtually identical. According to Definition~\ref{def:Day}, we have to produce a natural transformation
\[ \alpha: \Nat_{\cat{C}}(F_1,-) \smsh \Nat_{\cat{C}}(F_2,-) \to \Nat_{\cat{C}}(F_1 \smsh F_2,- \smsh -) \]
which is given by taking the smash product of natural transformations.

We then have to show that $\alpha$ induces equivalences
\[ \begin{diagram}\dgARROWLENGTH=1em \node{\Hom_{\Fun(\cat{F}_{\cat{C}},\spectra)}(\Nat_{\cat{C}}(F_1 \smsh F_2,-), {A})} \arrow{s,l}{\alpha^*} \\ \node{\Hom_{\Fun(\cat{F}_{\cat{C}} \times \cat{F}_{\cat{C}},\spectra)} (\Nat_{\cat{C}}(F_1,-) \smsh \Nat_{\cat{C}}(F_2,-), {A}(- \smsh -))} \end{diagram} \]
for arbitrary ${A}: \cat{F}_{\cat{C}} \to \spectra$.

First note that since $\Nat_{\cat{C}}(F_1 \smsh F_2,-)$ and $\Nat_{\cat{C}}(F_1,-) \smsh \Nat_{\cat{C}}(F_2,-)$ are reduced, it is sufficient to prove $\alpha^*$ is an equivalence when ${A}$ is reduced. (This is because any natural transformation out of a reduced functor between pointed $\infty$-categories factors, up to equivalence, via the universal reduction of its target.)

Now notice that a functor of the form $\Nat_{\cat{C}}(G,-)$ is linear and hence is equivalent to the linearization of
\[ \Sigma^\infty \Omega^\infty \Nat_{\cat{C}}(G,-) \homeq \Sigma^\infty \Hom_{\cat{F}_{\cat{C}}}(G,-). \]
We therefore have an equivalence
\[ \begin{split}
  \Nat_{\cat{C}}(G,-)
  &\homeq \hocolim_{k \to \infty} \Sigma^{-k} \Sigma^\infty \Hom_{\cat{F}_{\cat{C}}}(G,\Sigma^k(-)) \\
  &\homeq \hocolim_{k \to \infty} \Sigma^{-k} \Sigma^\infty \Hom_{\cat{F}_{\cat{C}}}(\Sigma^{-k}G,-).
\end{split} \]
Similarly, the natural transformation $\alpha$ can be identified with the map
\[ \begin{diagram} \dgARROWLENGTH=1em
  \node{\hocolim_{k \to \infty} \Sigma^{-k} \Sigma^\infty \Hom_{\cat{F}_{\cat{C}}}(\Sigma^{-k}F_1 \smsh F_2,- \smsh -)} \\
  \node{\hocolim_{k_1,k_2 \to \infty} \Sigma^{-k_1-k_2} \Sigma^\infty \Hom_{\cat{F}_{\cat{C}}}(\Sigma^{-k_1}F_1,-) \smsh \Hom_{\cat{F}_{\cat{C}}}(\Sigma^{-k_2}F_2,-)} \arrow{n}
\end{diagram} \]
given by inclusion into the term with $k = k_1+k_2$, and therefore, by the Yoneda Lemma (\ref{lem:yoneda}), the map $\alpha^*$ is equivalent to:
\[ \begin{diagram} \dgARROWLENGTH=1em
  \node{\holim_{k \to \infty} \Omega^\infty \Sigma^k {A}(\Sigma^{-k}F_1 \smsh F_2)} \arrow{s} \\
  \node{\holim_{k_1,k_2 \to \infty} \Omega^\infty \Sigma^{k_1+k_2}{A}(\Sigma^{-k_1}F_1 \smsh \Sigma^{-k_2}F_2)}
\end{diagram} \]
induced by projecting onto the term $k = k_1+k_2$. This map is an equivalence since the diagonal map $\mathbb{N} \to \mathbb{N} \times \mathbb{N}$ is final.
\end{proof}

\begin{remark}
One way to interpret the proof of Lemma~\ref{lem:Day} is that the ordinary Day convolution, as in Definition~\ref{def:Day}, is equivalent to a spectrally-enriched version of Day convolution, where the universal property is satisfied only with respect to spectrally-enriched (or, by~\cite[4.22]{blumberg/gepner/tabuada:2013}, exact) functors $\cat{F}_{\cat{C}} \to \spectra$.
\end{remark}

Theorem~\ref{thm:dn} is now a fairly simple consequence of Proposition~\ref{prop:der_n} and Lemma~\ref{lem:Day}, though we have to be careful about how we extend to $X_1,\dots,X_n \in \spectra(\cat{C})$.

\begin{proof}[Proof of Theorem~\ref{thm:dn}]
We have to prove that for $X_1,\dots,X_n \in \spectra(\cat{C})$, the functor $\der_n(-)(X_1,\dots,X_n)$ is a Day convolution of the form
\[ \der_1(-)(X_1) \otimes \dots \otimes \der_1(-)(X_n). \]
First suppose that $X_i = \Sigma^\infty_{\cat{C}} x_i$ where $x_1,\dots,x_n$ are compact objects in $\cat{C}$. In this case, we can apply Lemma~\ref{lem:Day} with the functors $F_i = \Sigma^\infty \Sigma^L\Omega^L\Hom_{\cat{C}}(x_i,-)$ (which are in $\cat{F}_{\cat{C}}$ since each $x_i$ is compact).

It follows easily from Definition~\ref{def:Day} that the Day convolution commutes with colimits in each variable, so we can take the homotopy colimit as $L \to \infty$ of the result of Lemma~\ref{lem:Day} and, using Proposition~\ref{prop:der_n}, we get
\begin{equation} \label{eq:Delta-conv} \Delta_n(-)(X_1,\dots,X_n) \homeq \Delta_1(-)(X_1) \otimes \dots \otimes \Delta_1(-)(X_n). \end{equation}
Recalling that $\der_nF(X_1,\dots,X_n) \homeq \Delta_nF(X_1,\dots,X_n)$ for spectrum-valued functors, we see that (\ref{eq:Delta-conv}) is precisely the desired equivalence in this case.

Next, note that arbitrary objects $x_1,\dots,x_n$ in the compactly-generated $\infty$-category $\cat{C}$ can be written as filtered colimits of compact objects. We can therefore recover the case of general $x_1,\dots,x_n$ from (\ref{eq:Delta-conv})) since again the Day convolution commutes with filtered colimits.

Finally, consider the case of general $X_1,\dots,X_n \in \spectra(\cat{C})$. The linearization of the functor $\Sigma^\infty_{\cat{C}}\Omega^\infty_{\cat{C}}$ is equivalent to the identity functor on $\spectra(\cat{C})$, that is:
\[ X_i \homeq P_1(\Sigma^\infty_{\cat{C}}\Omega^\infty_{\cat{C}})(X_i) \homeq \hocolim_{k \to \infty} \Omega^k \Sigma^\infty_{\cat{C}}\Omega^\infty_{\cat{C}}\Sigma^k X_i. \]
In other words, an arbitrary object of $\spectra(\cat{C})$ can be built from suspension spectrum objects by filtered colimit and desuspension, both of which commute with the Day convolution. We thus obtain the case of general $X_1,\dots,X_n$ from that of suspension spectrum objects.
\end{proof}

\section{Derivatives of arbitrary functors} \label{sec:generalF}

In this section we describe models for the derivatives of an arbitrary (reduced) functor $F: \cat{C} \to \cat{D}$ between pointed compactly-generated $\infty$-categories. In particular, we deduce that the terms in the coloured operad $\mathbb{I}_{\cat{C}}$ described in Remark~\ref{rem:I_C} are given by the derivatives of the identity functor on $\cat{C}$. This claim is a consequence of the following theorem.

\begin{theorem} \label{thm:dF}
For reduced $F: \cat{C} \to \cat{D}$, $X_1,\dots,X_n \in \spectra(\cat{C})$ and $Y \in \spectra(\cat{D})$, we have
\[ \der_nF(X_1,\dots,X_n;Y) \homeq \Nat_{\cat{F}_{\cat{D}}}(\der_1(-)(Y), \der_n(- \circ F)(X_1,\dots,X_n)). \]
\end{theorem}

\begin{corollary} \label{cor:dernI}
	For any pointed compactly-generated $\infty$-category $\cat{C}$, we have
	\[ \der_nI_{\cat{C}}(X_1,\dots,X_n;Y) \homeq \Nat_{\cat{F}_{\cat{C}}}(\der_1(-)(Y), \der_n(-)(X_1,\dots,X_n)). \]
\end{corollary}

\begin{example}
	When $\cat{C} = \cat{D} = \based$, we have
	\[ \der_nF \homeq \Nat(\der_1,\der_n(- \circ F)) \]
	and, in particular,
	\[ \der_nI_{\cat{C}} \homeq \Nat(\der_1,\der_n) \homeq \Nat(\der_1,\der_1^{\otimes n}). \]
	In other words, the derivatives of the identity functor on $\based$ form the coendomorphism operad of the functor $\der_1: \Fun(\based,\spectra) \to \spectra$ with respect to Day convolution. In \cite{ching:2005} an operad structure on these derivatives was constructed by taking the Koszul dual of the commutative operad in spectra. It is not obvious that these two operad structures on $\der_*I_{\based}$ are equivalent, though, as we will see in the proof of Theorem~\ref{thm:dF} below, both depend on the cosimplicial resolution of the identity functor via the adjunction $(\Sigma^\infty,\Omega^\infty)$, making a connection plausible.
\end{example}

Before giving the proof of \ref{thm:dF}, let us construct the map that realizes the desired equivalence. That map is based on a natural transformation
\begin{equation} \label{eq:delta} c: \Delta_1G(\Delta_nF(X_1,\dots,X_n)) \to \Delta_n(GF)(X_1,\dots,X_n) \end{equation}
which we now define for $F: \cat{C} \to \cat{D}$ reduced and $G: \cat{D} \to \spectra$ reduced and preserving filtered colimits.

\begin{definition} \label{def:delta}
With $F$ and $G$ as above, there is a natural transformation of symmetric functors $\cat{C}^n \to \spectra$
\[ c'': G\creff_nF \to \creff_n(GF) \]
coming from the definition of the cross-effect as a total homotopy fibre. Taking multilinearization of $c''$, we get a map of symmetric multilinear functors $\cat{C}^n \to \spectra$:
\[ c': P_{1,\dots,1}(G\creff_nF) \to P_{1,\dots,1}(\creff_n(GF)) \]
which corresponds to a map of symmetric multilinear functors $\spectra(\cat{C})^n \to \spectra$ which we can write
\[ c: \Delta_{1,\dots,1}(G\creff_nF) \to \Delta_{1,\dots,1}(\creff_n(GF)). \]
The target of $c$ is equivalent to $\Delta_n(GF)$ by \cite[6.1]{goodwillie:2003} (or~\cite[6.1.3.23]{lurie:2017} in the $\infty$-categorical setting). To identify the source of $c$ with $\Delta_1G(\Delta_nF)$ we apply~\cite[6.2.1.22]{lurie:2017}, i.e. the multivariable version of the Klein-Rognes~\cite{klein/rognes:2002} chain rule for first derivatives.
\end{definition}

\begin{definition} \label{def:dF-equiv}
For any reduced $G: \cat{D} \to \spectra$, the functor $\Delta_1G = \der_1G: \spectra(\cat{D}) \to \spectra$ is, by definition, linear, and hence enriched over $\spectra$, at least up to homotopy. In other words we have natural maps
\[ \Map_{\spectra(\cat{D})}(Y,\Delta_nF(X_1,\dots,X_n)) \to \Nat_{\cat{F}_{\cat{D}}}(\Delta_1(-)(Y),\Delta_1(-) (\Delta_nF(X_1,\dots,X_n))). \]
Composing with our map $c$ from (\ref{eq:delta}) we get
\begin{equation} \label{eq:map} \Map_{\spectra(\cat{D})}(Y,\Delta_nF(X_1,\dots,X_n)) \to \Nat_{\cat{F}_{\cat{D}}}(\Delta_1(-)(Y),\Delta_n(- \circ F)(X_1,\dots,X_n)) \end{equation}
which is the form of the equivalence required by Theorem~\ref{thm:dF}.
\end{definition}

\begin{proof}[Proof of Theorem~\ref{thm:dF}]
First note that each side of the desired equivalence commutes with desuspension and filtered colimits in the variable $Y$. The argument in the proof of \ref{thm:dn} implies that it is sufficient to consider the case $Y = \Sigma^\infty_{\cat{D}}y$ for some compact object $y$ in the compactly-generated $\infty$-category $\cat{D}$.

Using Proposition~\ref{prop:der_n}, the right-hand side of the desired equivalence can be written in the form
\[ \holim_{L \to \infty} \Nat_{\cat{F}_{\cat{D}}}(\Nat_{\cat{D}}(\Sigma^\infty \Sigma^L\Omega^L \Hom_{\cat{D}}(y,\bullet),-),\der_n(- \circ F)(X_1,\dots,X_n)) \]
which, by a stable version of the Yoneda Lemma \cite[6.4]{nikolaus:2016}, is equivalent to
\[ \holim_{L \to \infty} \der_n(\Sigma^\infty \Sigma^L\Omega^L \Hom_{\cat{D}}(y,F))(X_1,\dots,X_n). \]
On the other hand, for the left-hand side of the desired result, we have an equivalence
\[ \der_nF(X_1,\dots,X_n;\Sigma^\infty_\cat{D}y) \homeq \der_n(\Hom_{\cat{D}}(y,F))(X_1,\dots,X_n) \]
which follows from the fact that
\[ \Hom_{\cat{D}}(y,\Delta_nF) \homeq \Delta_n(\Hom_{\cat{D}}(y,F)) \]
for a compact object $y \in \cat{D}$.

It is now sufficient to show that, for reduced $G: \cat{C} \to \based$, there is a natural equivalence
\begin{equation} \label{eq:holim} \alpha: \der_nG \weq \holim_{L \to \infty} \der_n(\Sigma^\infty \Sigma^L \Omega^L G) \end{equation}
where the map $\alpha$ has components given by
\[ \Delta_n(G) \homeq \Delta_1(\Sigma^\infty \Sigma^L \Omega^L)\Delta_n(G) \arrow{e,t}{c} \Delta_n(\Sigma^\infty \Sigma^L \Omega^L G) \]
with $c$ as in (\ref{eq:delta}). This claim contains the real substance of the result we are trying to prove, and it occupies the majority of our effort here.

We first show that $\alpha$ is an equivalence when $G = \Omega^\infty \mathbb{G}$ for some $\mathbb{G}: \cat{C} \to \spectra$ (in which case note that $\der_nG \homeq \der_n\mathbb{G}$). Then there is a map
\[ \beta: \holim_{L \to \infty} \der_n(\Sigma^\infty \Sigma^L \Omega^L \Omega^\infty \mathbb{G}) \to \der_n\mathbb{G} \]
given by projection onto the $L = 0$ term followed by the counit of the adjunction $(\Sigma^\infty,\Omega^\infty)$. It is easy to check from the definitions that $\beta\alpha$ is equivalent to the identity. It is therefore sufficient to show that $\beta$ is an equivalence.

For this task, we need an instance of the chain rule for spectrum-valued functors which tells us that there is an equivalence
\[ \der_n(\Sigma^\infty \Sigma^L \Omega^L \Omega^\infty \mathbb{G}) \homeq \prod_{\mathsf{P}(n)} \der_k(\Sigma^\infty \Sigma^L \Omega^L \Omega^\infty) \smsh \der_{n_1}\mathbb{G} \smsh \dots \smsh \der_{n_k}\mathbb{G} \]
where $\mathsf{P}(n)$ is the set of unordered partitions of the set $\{1,\dots,n\}$, where $n_1,\dots,n_k$ denote the sizes of the pieces of a partition, and where we have suppressed the dependence on variables $X_1,\dots,X_n \in \spectra(\cat{C})$ for the sake of readability. This result is a generalization of the main theorem of \cite{ching:2010} with a similar proof. Details are provided in Appendix~\ref{app:chain}.

The source of the map $\beta$ thus splits as
\begin{equation} \label{eq:beta} \prod_{\mathsf{P}(n)} \; \holim_{L \to \infty} \; [\der_k(\Sigma^\infty \Sigma^L \Omega^L \Omega^\infty) \smsh \der_{n_1}\mathbb{G} \smsh \dots \smsh \der_{n_k}\mathbb{G}] \end{equation}
and $\beta$ is given by projection onto the term corresponding to the indiscrete partition, i.e. with $k = 1$. (Notice that in this term all the maps in the inverse system are equivalences and the homotopy limit is just $\der_n\mathbb{G}$.)

A standard calculation shows that $\der_k(\Sigma^\infty \Sigma^L \Omega^L \Omega^\infty) \homeq S^{-L(k-1)}$, that is, a negative-dimensional sphere spectrum. The maps in the inverse systems in (\ref{eq:beta}) are induced by the counit $\Sigma \Omega \to I$ via maps
\[ \der_k(\Sigma^\infty \Sigma^{L+1} \Omega^{L+1} \Omega^\infty) \to \der_k(\Sigma^\infty \Sigma^L \Omega^L \Omega^\infty) \]
and hence are trivial when $k > 1$ for dimension reasons. It follows that the homotopy limits appearing in (\ref{eq:beta}) are trivial when $k > 1$, and hence that the projection map $\beta$ is an equivalence. This completes the proof that the map $\alpha$ is an equivalence when $G = \Omega^\infty\mathbb{G}$.

Now consider arbitrary reduced $G: \cat{C} \to \based$. There is a commutative diagram
\[ \begin{diagram} \dgARROWLENGTH=2em
  \node{\der_nG} \arrow{e} \arrow{s,l}{\alpha}
    \node{\Tot \der_n(\Omega^\infty(\Sigma^\infty\Omega^\infty)^\bullet \Sigma^\infty G)} \arrow{s,r}{\Tot \alpha} \\
  \node{\holim_{L \to \infty} \der_n(\Sigma^\infty \Sigma^L \Omega^L G)} \arrow{e}
    \node{\Tot \holim_{L \to \infty} \der_n(\Sigma^\infty \Sigma^L \Omega^L \Omega^\infty (\Sigma^\infty\Omega^\infty)^\bullet \Sigma^\infty G)}
\end{diagram} \]
where $\Tot$ denotes the totalization of cosimplicial spectra which are built from the $(\Sigma^\infty,\Omega^\infty)$ adjunction. The horizontal maps are equivalences by induction on the Taylor tower of $G$ (by the argument of \cite[4.1.1]{arone/ching:2011} and using the fact that $\Tot$ commutes with $\holim$), and the right-hand vertical map is an equivalence by the case already considered. Therefore the map $\alpha$ is an equivalence for arbitrary $G$. This completes the proof of Theorem~\ref{thm:dF}.
\end{proof}

\begin{remark} \label{rem:pro-operad}
A key part of the proof of Theorem~\ref{thm:dF} was the construction of the equivalence
\[ \alpha: \der_nG \arrow{e,t}{\sim} \holim_{L \to \infty} \der_n(\Sigma^\infty \Sigma^L \Omega^L G) \]
for a functor $G: \cat{C} \to \based$. In particular, when $G = I_{\based}$, we get
\[ \der_*I_{\based} \homeq \holim_{L \to \infty} \der_*(\Sigma^\infty \Sigma^L \Omega^L). \]
The terms in the homotopy limit on the right-hand side turn out to be equivalent to the operads $\mathsf{K}(\mathsf{E}_L)$ given by the Koszul duals of the stable little $L$-discs operad, themselves equivalent to desuspensions of those little disc operads by~\cite{ching/salvatore:2020}, and this formula expresses $\der_*I_{\based}$ as the inverse limit of a `pro-operad'. Similarly, we have an equivalence
\[ \der_*I_{\spectra} \homeq \der_*\Omega^\infty \homeq \holim_{L \to \infty} \der_*(\Sigma^\infty \Sigma^L \Omega^L \Omega^\infty) \homeq \holim_{L \to \infty} \mathbf{S}^{-L} \]
which expresses $\der_*I_{\spectra}$ as the inverse limit of a pro-operad whose components are certain operads $\mathbf{S}^{-L}$ formed by desuspensions of the sphere operad.

In \cite{arone/ching:2016}, Arone and the author showed that these two pro-operads classify the Taylor towers of functors $\based \to \spectra$ and $\spectra \to \spectra$ respectively. We believe that an analogous pro-operad can be constructed for any pointed, compactly-generated $\infty$-category $\cat{C}$. The inverse limit of this pro-operad is equivalent to the operad $\der_*I_{\cat{C}}$ and modules over the pro-operad should classify the Taylor towers of functors $\cat{C} \to \spectra$.
\end{remark}

\begin{remark} \label{rem:bimod}
Theorem~\ref{thm:dF} provides models of the derivatives of a functor $F: \cat{C} \to \cat{D}$ that admit natural composition maps in the following sense. Define a collection $\mathbb{D}_F$ of spectra by
\[ \mathbb{D}_F(X_1,\dots,X_n;Y) := \Nat_{\cat{F}_{\cat{D}}}(\der_1(-)(Y),(\der_1(-)(X_1) \otimes \dots \otimes \der_1(-)(X_n))(- \circ F)) \]
for $X_1,\dots,X_n \in \spectra(\cat{C})$ and $Y \in \spectra(\cat{D})$. By Theorem~\ref{thm:dF}, these spectra are the derivatives of $F$. Notice that $\mathbb{D}_{I_\cat{C}}$ is the same collection of spectra as in the coloured operad $\mathbb{I}_{\cat{C}}$ in Remark~\ref{rem:I_C}.

Now suppose we have reduced functors $F: \cat{C} \to \cat{D}$ and $G: \cat{D} \to \cat{E}$ that preserve filtered colimits. Then we can build maps of the form
\[ \begin{diagram} \node{\mathbb{D}_G(Y_1,\dots,Y_k;Z) \smsh \mathbb{D}_F(\un{X}_1;Y_1) \smsh \dots \smsh \mathbb{D}_F(\un{X}_k;Y_k)} \arrow{s} \\ \node{\mathbb{D}_{GF}(\un{X}_1,\dots,\un{X}_k;Z)} \end{diagram} \]
where $Z \in \spectra(\cat{E})$, $Y_1,\dots,Y_k \in \spectra(\cat{D})$ and each $\un{X}_i$ is a sequence of objects in $\spectra(\cat{C})$. In particular, the derivatives $\mathbb{D}_F$ form a bimodule over the operads $\mathbb{I}_{\cat{C}}$ and $\mathbb{I}_{\cat{D}}$ described in Remark~\ref{rem:I_C}, at least up to homotopy. This structure is made precise, in the context of $\infty$-operads, in Section~\ref{sec:bimodule}.
\end{remark}

\section{Stable \texorpdfstring{$\infty$}{infinity}-operads and derivatives of the identity} \label{sec:operad}

In this section we provide a formal definition of the operad $\mathbb{I}_{\cat{C}}$ of Remark~\ref{rem:I_C} in the context of Lurie's theory of $\infty$-operads. Here is an outline of our main construction.

We start by describing a symmetric monoidal $\infty$-category that represents the objectwise smash product of functors $\cat{C} \to \spectra$, and hence the desired monoidal product on $\cat{F}_{\cat{C}}$. Then we turn to the Day convolution, using a construction of Lurie to describe a symmetric monoidal $\infty$-category that represents the convolution of functors $\cat{F}_{\cat{C}} \to \spectra$.

Some care is needed here because the $\infty$-category $\cat{F}_{\cat{C}}$ is not small. However, it is generated under filtered colimits by a small symmetric monoidal subcategory $\cat{F}_{\cat{C}}^{\omega}$. We construct a symmetric monoidal $\infty$-category $\Fun(\cat{F}_{\cat{C}}^{\omega},\spectra)^{\otimes}$ that represents the Day convolution of functors $\cat{F}_{\cat{C}}^{\omega} \to \spectra$, and note that the proof of Theorem~\ref{thm:dn} carries over to this context.

As Remark~\ref{rem:I_C} shows, we are interested in morphisms \emph{into} the Day convolution rather than out of it, so we next apply work of Barwick, Glasman and Nardin \cite{barwick/glasman/nardin:2018} to construct a symmetric monoidal $\infty$-category $\Fun(\cat{F}_{\cat{C}}^{\omega},\spectra)^{op,\otimes}$ with the same monoidal product, i.e. Day convolution, but with the \emph{opposite} underlying $\infty$-category.

Finally, in Definition~\ref{def:I}, we restrict to the full subcategory of $\Fun(\cat{F}_{\cat{C}}^{\omega},\spectra)^{op,\otimes}$ generated by those objects of the form $\der_1(-)(X)$ for $X \in \spectra(\cat{C})$. The resulting $\infty$-operad $\mathbb{I}_{\cat{C}}^{\otimes}$ is a precise version of the operad described informally in Remark~\ref{rem:I_C}.

We should also note that all the $\infty$-operads appearing in our work, including the symmetric monoidal structures, are \emph{non-unital} in the sense that they do not encode unit objects. We start our description of these constructions by recalling the basic theory of $\infty$-operads from \cite{lurie:2017} with some slight adjustment to take into account our focus on the non-unital case. The language of (co)cartesian edges and fibrations from~\cite[2.4]{lurie:2009} features heavily in the remainder of this paper.

\begin{definition} \label{def:op}
Let $\Surj$ be the category whose objects are pointed finite sets and whose morphisms are surjections which preserve the basepoint. We write $\langle n \rangle := \{*,1,\dots,n\}$. A morphism in $\Surj$ is \emph{inert} if the inverse image of every non-basepoint contains exactly one element. For example, let $\rho_i: \langle n \rangle \to \langle 1 \rangle$ denote the inert morphism with $\rho_i(i) = 1$ and $\rho_i(j) = *$ for $j \neq i$. A morphism is \emph{active} if the inverse image of the basepoint consists only of the basepoint.

A \emph{non-unital $\infty$-operad} is a map of $\infty$-categories of the form
\[ p: \cat{O}^{\otimes} \to \Surj \]
that satisfies the following conditions:
\begin{enumerate}
  \item for every object $\un{X} \in \cat{O}^{\otimes}$, every inert morphism $\alpha$ in $\Surj$ with source $p(\un{X})$ has a $p$-cocartesian lift $\bar{\alpha}$ in $\cat{O}^{\otimes}$ with source $\un{X}$;
  \item for every $n \geq 0$, the $p$-cocartesian lifts $\bar{\rho}_i$ determine an equivalence of $\infty$-categories
      \[ \bar{\rho}: \cat{O}^{\otimes}_{\langle n \rangle} \homeq (\cat{O}^{\otimes}_{\langle 1 \rangle})^{n} \]
      where $\cat{O}^{\otimes}_{\langle n \rangle}$ denotes the fibre $p^{-1}(\langle n \rangle)$;
  \item for every morphism $\alpha: \langle m \rangle \to \langle n \rangle$ in $\Surj$ and every pair of objects $\un{X},\un{Y} \in \cat{O}^{\otimes}$ with $p(\un{Y}) = \langle n \rangle$, the $p$-cocartesian lifts $\bar{\rho}_i : \un{Y} \to Y_i$ determine an equivalence
      \[ \Hom_{\cat{O}^{\otimes}}(\un{X},\un{Y})_{\alpha} \weq \prod_{i = 1}^{n} \Hom_{\cat{O}^{\otimes}}(\un{X},Y_i)_{\rho_i\alpha} \]
      between spaces of morphisms for the $\infty$-category $\cat{O}^{\otimes}$; we use a subscript to denote the subspace consisting of those morphisms that map by $p$ to the given map in $\Surj$.
\end{enumerate}
Since all $\infty$-operads appearing in this paper are non-unital, we drop that adjective. We also commonly leave the map $p$ implied and refer to \emph{the $\infty$-operad $\cat{O}^{\otimes}$}. We write $\cat{O} = \cat{O}^{\otimes}_{\langle 1 \rangle}$ and refer to this as the \emph{underlying $\infty$-category} for the $\infty$-operad $\cat{O}^{\otimes}$.
\end{definition}

\begin{remark} \label{rem:operad}
An object $\un{X} \in \cat{O}^{\otimes}$ with $p(\un{X}) = S$ can be identified with a collection of objects of $\cat{O}$ indexed by $S$: a bijection $\alpha: S \isom \langle n \rangle$ induces a sequence of equivalences
\[ \cat{O}^{\otimes}_S \arrow{e,tb}{\bar{\alpha}}{\homeq} \cat{O}^{\otimes}_{\langle n \rangle} \homeq \cat{O}^n \homeq \prod_{S} \cat{O}. \]
Based on this observation, we typically use a finite sequence of objects in $\cat{O}$ as a representative for an arbitrary object of $\cat{O}^{\otimes}$. For example, in part (iii) of Definition~\ref{def:op} we identify the object $\un{Y}$ with the sequence $(Y_1,\dots,Y_n)$.
\end{remark}

\begin{remark} \label{rem:multi}
An $\infty$-operad $\cat{O}^{\otimes}$ is an $\infty$-categorical version of a simplicial coloured operad whose colours are the objects of the underlying $\infty$-category $\cat{O}$. Given objects  $X_1,\dots,X_n,Y \in \cat{O}$, for $n \geq 1$, we write
\[ \Hom_{\cat{O}^{\otimes}}(X_1,\dots,X_n;Y) := \Hom_{\cat{O}^{\otimes}}((X_1,\dots,X_n),Y)_{\langle n \rangle \to \langle 1 \rangle} \]
for the fibre of the morphism space in $\cat{O}^{\otimes}$ over the unique active morphism $\langle n \rangle \to \langle 1 \rangle$ in $\Surj$. We call these spaces the \emph{multi-morphism spaces} of the $\infty$-operad $\cat{O}^{\otimes}$. They admit composition maps that are associative up to homotopy and through which we can view $\cat{O}^{\otimes}$ as the analogue of a coloured operad (or symmetric multicategory) enriched in simplicial sets. The definition of $\infty$-operad ensures that all mapping spaces for $\cat{O}^{\otimes}$ are determined by the multi-morphism spaces described here.
\end{remark}

\begin{definition}
Given $\infty$-operads $p_1: \cat{O}_1^{\otimes} \to \Surj$ and $p_2: \cat{O}_2^{\otimes} \to \Surj$, a \emph{map of $\infty$-operads} $g: \cat{O}_1^{\otimes} \to \cat{O}_2^{\otimes}$ is a functor $g$ such that $p_2 \circ g = p_1$, and that sends $p_1$-cocartesian lifts in $\cat{O}_1^{\otimes}$ of inert maps in $\Surj$ to $p_2$-cocartesian lifts in $\cat{O}_2^{\otimes}$. An \emph{equivalence} of $\infty$-operads is a map of $\infty$-operads that is an equivalence on the underlying $\infty$-categories.
\end{definition}

\begin{definition} \label{def:suboperad}
Let $p: \cat{O}^{\otimes} \to \Surj$ be an $\infty$-operad, and let $\cat{O}'$ be a full subcategory of the underlying $\infty$-category $\cat{O}$. Then we let $\cat{O}'^{\otimes}$ be the full subcategory of $\cat{O}^{\otimes}$ whose objects are those equivalent (via the identifications of Remark~\ref{rem:operad}) to sequences $(X_1,\dots,X_n)$ where $X_1,\dots,X_n \in \cat{O}'$. The restriction of $p$ to $\cat{O}'^{\otimes}$ is also an $\infty$-operad, and the inclusion $\cat{O}'^{\otimes} \to \cat{O}^{\otimes}$ is a map of $\infty$-operads. We refer to $\cat{O}'^{\otimes}$ as \emph{the suboperad of $\cat{O}^{\otimes}$ generated by $\cat{O}'$}.
\end{definition}

\begin{definition} \label{def:symmon}
A \emph{non-unital symmetric monoidal $\infty$-category} is a non-unital $\infty$-operad $p: \cat{C}^{\otimes} \to \Surj$ such that $p$ is a cocartesian fibration. This condition implies that for $X_1,\dots,X_n,Y \in \cat{C}$, we have
\[ \Hom_{\cat{C}^{\otimes}}(X_1,\dots,X_n;Y) \homeq \Hom_{\cat{C}}(X_1 \otimes \dots \otimes X_n, Y) \]
for some object $X_1 \otimes \dots \otimes X_n$ that depends functorially on $X_1,\dots,X_n$, and such that the operation $\otimes$ is associative and commutative up to higher coherent homotopies. This definition mimics the way in which a symmetric monoidal category can be expressed as a coloured operad.

A map of $\infty$-operads $g: \cat{C}_1^{\otimes} \to \cat{C}_2^{\otimes}$ between non-unital symmetric monoidal $\infty$-categories is \emph{symmetric monoidal} if it takes all cocartesian morphisms in $\cat{C}_1^{\otimes}$ to cocartesian morphisms in $\cat{C}_2^{\otimes}$.
\end{definition}

The $\infty$-operads we study in this paper are stable in the following sense.

\begin{definition} \label{def:stable-operad}
An $\infty$-operad $\cat{O}^{\otimes}$ is \emph{stable} if the underlying $\infty$-category $\cat{O}$ is stable (in the sense of~\cite[1.1.1.9]{lurie:2017}) and, for each $n \geq 1$, the functor
\[ (\cat{O}^{op})^n \times \cat{O} \to \spaces; \quad (X_1,\dots,X_n,Y) \mapsto \Hom_{\cat{O}^{\otimes}}(X_1,\dots,X_n;Y) \]
preserves finite limits in each variable. In that case, those functors are linear in each variable and so factor via corresponding spectrum-valued functors which we denote
\[ \Map_{\cat{O}^{\otimes}}(X_1,\dots,X_n;Y). \]
We refer to these objects as the \emph{multi-morphism spectra} of the stable $\infty$-operad $\cat{O}^{\otimes}$.
\end{definition}

\begin{example}
A symmetric monoidal $\infty$-category $\cat{C}^{\otimes}$ is stable if and only if $\cat{C}$ is stable and the monoidal product $\otimes$ is exact (i.e. preserves finite limits and finite colimits) in each variable. In that case we have
\[ \Map_{\cat{C}^{\otimes}}(X_1,\dots,X_n,Y) \homeq \Map_{\cat{C}}(X_1 \otimes \dots \otimes X_n,Y). \]
\end{example}

\begin{example}
There is a non-unital symmetric monoidal $\infty$-category $\spectra^{\smsh} \to \Surj$ whose underlying $\infty$-category is $\spectra$ and whose monoidal structure represents the ordinary smash product of spectra. See \cite[4.8.2]{lurie:2017} for the unital version which is a cocartesian fibration $\spectra^{\smsh}_{\mathsf{u}} \to \Fin$. The corresponding non-unital $\infty$-operad is given by pulling back this fibration along the inclusion $\Surj \to \Fin$.
\end{example}

\begin{example} \label{ex:finspec}
Let $\cat{O}^{\otimes}$ be a stable $\infty$-operad with $\cat{O}$ equivalent to the $\infty$-category of finite spectra. Then the multi-morphism spectra for $\cat{O}^{\otimes}$ are determined by their values on the sphere spectrum. In particular, the data of $\cat{O}^{\otimes}$ are determined by the symmetric sequence of spectra
\[ \mathbf{O}(n) := \Map_{\cat{O}^{\otimes}}(\underbrace{S^0,\dots,S^0}_n;S^0) \]
together with appropriate composition maps (that are associative up to higher coherent homotopies). In this way, $\cat{O}^{\otimes}$ can be viewed as the $\infty$-categorical version of an ordinary monochromatic operad of spectra.
\end{example}

We now turn to the main subject of this section, and we start with the construction of a symmetric monoidal $\infty$-category that represents the objectwise smash product on $\cat{F}_{\cat{C}}$.

\begin{construction} \label{cons:objectwise}
Consider the pullback of simplicial sets of the form
\[ \begin{diagram} \dgARROWLENGTH=1em
  \node{\Fun(\cat{C},\spectra)^{\smsh}} \arrow{e} \arrow{s,l}{p_{\cat{C}}} \node{\Fun(\cat{C},\spectra^{\smsh})} \arrow{s} \\
  \node{\Surj} \arrow{e} \node{\Fun(\cat{C},\Surj)}
\end{diagram} \]
where the right-hand map is induced by the cocartesian fibration $\spectra^{\smsh} \to \Surj$ and the bottom map sends a finite pointed set to the constant functor with that value. The induced map $p_{\cat{C}}$ is then also a cocartesian fibration of $\infty$-operads, with fibres
\[ \Fun(\cat{C},\spectra)^{\smsh}_{\langle n \rangle} \homeq \Fun(\cat{C},\spectra^{\smsh}_{\langle n \rangle}). \]
Thus $p_{\cat{C}}$ is a (non-unital) symmetric monoidal $\infty$-category with underlying $\infty$-category $\Fun(\cat{C},\spectra)$ and monoidal product given by the objectwise smash product of functors; see \cite[2.1.3.4]{lurie:2017}.
\end{construction}

\begin{definition} \label{def:obinfty}
Let $\cat{F}_{\cat{C}}^{\smsh} \to \Surj$ denote the restriction of the symmetric monoidal $\infty$-category $p_{\cat{C}}$ of Construction~\ref{cons:objectwise} to the full subcategory generated by those functors $\cat{C} \to \spectra$ that are reduced and preserve filtered colimits. Since this collection of functors is closed under the objectwise smash product, $\cat{F}_{\cat{C}}^{\smsh}$ is also a non-unital symmetric monoidal $\infty$-category. (Note that we are forced to deal with a non-unital symmetric monoidal structure because the unit object for the objectwise smash product is not a reduced functor.)
\end{definition}

Our next goal is to describe a non-unital symmetric monoidal $\infty$-category that represents the Day convolution of functors introduced in Definition~\ref{def:Day}. We use a construction of Lurie from~\cite[2.2.6]{lurie:2017}. Glasman~\cite{glasman:2016} describes a similar construction for ordinary (i.e. not non-unital) symmetric monoidal $\infty$-categories.

\begin{construction} \label{cons:day-infty}
Let $\cat{C}^{\otimes}$ and $\cat{D}^{\otimes}$ be non-unital symmetric monoidal $\infty$-categories such that $\cat{C}$ is small, $\cat{D}$ admits all small colimits, and the monoidal structure on $\cat{D}$ preserves colimits in each variable. Applying~\cite[2.2.6.7]{lurie:2017}  with $\cat{O}^{\otimes} = \Surj$, we get a non-unital symmetric monoidal $\infty$-category $\Fun(\cat{C},\cat{D})^{\otimes}$ with the following universal property: for any $\infty$-operad $\cat{A}^{\otimes}$, there is an equivalence between the $\infty$-categories of $\infty$-operad maps $\cat{A}^{\otimes} \to \Fun(\cat{C},\cat{D})^{\otimes}$ and $\infty$-operad maps $\cat{A}^{\otimes} \times_{\Surj} \cat{C}^{\otimes} \to \cat{D}^{\otimes}$.
\end{construction}

We would like to apply \ref{cons:day-infty} to functors $\cat{F}_{\cat{C}} \to \spectra$, but since $\cat{F}_{\cat{C}}$ is not small, we cannot do this directly. However, $\cat{F}_{\cat{C}}$ is a compactly-generated $\infty$-category, i.e. is generated under filtered colimits by the small subcategory $\cat{F}^{\omega}_{\cat{C}}$ of compact objects.

\begin{lemma} \label{lem:compact}
Let $\cat{C}$ be a pointed compactly-generated $\infty$-category. Then the $\infty$-category $\cat{F}_{\cat{C}}$ is compactly-generated and the subcategory $\cat{F}^{\omega}_{\cat{C}}$ of compact objects is closed under the objectwise smash product of functors.
\end{lemma}
\begin{proof}
Let $\cat{R} \subseteq \cat{F}_{\cat{C}}$ be the full subcategory generated by the \emph{representable} functors, i.e. those of the form $R_x := \Sigma^\infty \Hom_{\cat{C}}(x,-)$ for some compact object $x \in \cat{C}$. It is a standard consequence of the Yoneda Lemma~\ref{lem:yoneda}, and the fact that equivalences in $\cat{F}_{\cat{C}}$ are detected objectwise on compact objects in $\cat{C}$, that an arbitrary $F \in \cat{F}_{\cat{C}}$ is the colimit of the canonical diagram
\[ F \homeq \colim_{R_x \to F} R_x, \]
indexed by the overcategory $\cat{R}_{/F}$. It follows, by the argument of \cite[5.3.4.17]{lurie:2009}, that an arbitrary $F$ is a filtered colimit of finite colimits of diagrams in $\cat{R}$, and therefore that the compact objects in $\cat{F}_{\cat{C}}$ are the retracts of those finite colimits. In particular, $\cat{F}_{\cat{C}}$ is compactly-generated.

Finally, from Lemma~\ref{lem:creff} it follows that the objectwise smash product of two representable functors is compact, since the cross-effect construction commutes with filtered colimits. Thus, the objectwise smash product of any two compact functors is compact.
\end{proof}

\begin{definition} \label{def:compact}
Let $(\cat{F}^{\omega}_{\cat{C}})^{\smsh} \to \Surj$ be the suboperad of the symmetric monoidal $\infty$-category $p_{\cat{C}}$ of Definition~\ref{def:obinfty} generated by the compact objects in $\cat{F}_{\cat{C}}$. By \cite[2.2.1.1]{lurie:2017}, this suboperad is an essentially small stable symmetric monoidal $\infty$-category.
\end{definition}

\begin{definition} \label{def:conv}
Applying Construction~\ref{cons:day-infty} to the symmetric monoidal $\infty$-category of the previous paragraph, we get a new stable non-unital symmetric monoidal $\infty$-category
\[ q_{\cat{C}}: \Fun(\cat{F}_{\cat{C}}^{\omega},\spectra)^{\otimes} \to \Surj. \]
\end{definition}

To proceed to the definition of the $\infty$-operad $\mathbb{I}_{\cat{C}}^{\otimes}$, we need one more general construction.

\begin{construction} \label{cons:symmon-op}
Let $q: \cat{E}^{\otimes} \to \Surj$ be a non-unital symmetric monoidal $\infty$-category. Then Barwick, Glasman and Nardin \cite[3.6]{barwick/glasman/nardin:2018} define another non-unital symmetric monoidal $\infty$-category $q^{\vee,op}: \cat{E}^{op,\otimes} \to \Surj$ that represents the induced symmetric monoidal structure on the opposite $\infty$-category of $\cat{E}$. Note that when $\cat{E}^{\otimes}$ is stable, so is $\cat{E}^{op,\otimes}$.
\end{construction}

\begin{definition}
Applying \ref{cons:symmon-op} to the non-unital symmetric monoidal $\infty$-category $q_{\cat{C}}$ of Definition~\ref{def:conv}, there is a stable non-unital symmetric monoidal $\infty$-category
\[ q_{\cat{C}}^{\vee,op}: \Fun(\cat{F}_{\cat{C}}^{\omega},\spectra)^{op,\otimes} \to \Surj \]
that represents the monoidal structure corresponding to Day convolution on the opposite $\infty$-category of the category of functors $\cat{F}_{\cat{C}}^{\omega} \to \spectra$. Note that the multi-morphism spectra in $\Fun(\cat{F}_{\cat{C}}^{\omega},\spectra)^{op,\otimes}$ are given by the mapping spectra
\[ \Nat_{\cat{F}_{\cat{C}}^{\omega}}({A},{B}_1 \otimes \dots \otimes {B}_n) \]
where $\otimes$ denotes the Day convolution of functors $\cat{F}^{\omega}_{\cat{C}} \to \spectra$. Comparing with Remark~\ref{rem:I_C}, this observation motivates the following definition, which is the central construction of this paper.
\end{definition}

\begin{definition} \label{def:I}
Let $\cat{C}$ be a pointed compactly-generated $\infty$-category. Then let $\mathbb{I}^{\otimes}_{\cat{C}}$ be the suboperad of the symmetric monoidal $\infty$-category $\Fun(\cat{F}_{\cat{C}}^{\omega},\spectra)^{op,\otimes}$ generated, in the sense of~\ref{def:suboperad}, by those objects essentially of the form
\[ \der_1(-)(X): \cat{F}_{\cat{C}}^{\omega} \to \spectra \]
for $X \in \spectra(\cat{C})$. We usually denote the object $\der_1(-)(X)$ of the underlying $\infty$-category $\mathbb{I}_{\cat{C}}$ simply by $X$.
\end{definition}

\begin{proposition} \label{prop:sp-op}
The $\infty$-operad $\mathbb{I}_{\cat{C}}^{\otimes}$ is stable and has multi-morphism spectra
\[ \Map_{\mathbb{I}_{\cat{C}}^{\otimes}}(X_1,\dots,X_n;Y) \homeq \der_nI_{\cat{C}}(X_1,\dots,X_n;Y). \]
In particular, the underlying $\infty$-category of $\mathbb{I}_{\cat{C}}^{\otimes}$ is equivalent to $\spectra(\cat{C})^{op}$.
\end{proposition}
\begin{proof}
Since $\mathbb{I}_{\cat{C}}^{\otimes}$ is a full subcategory of a stable symmetric monoidal $\infty$-category, it has multi-morphism spectra given by
\[ \Map_{\mathbb{I}_{\cat{C}}^{\otimes}}(X_1,\dots,X_n;Y) \homeq \Nat_{\cat{F}^{\omega}_{\cat{C}}}(\der_1(-)(Y), \der_1(-)(X_1) \otimes \dots \otimes \der_1(-)(X_n)). \]
Note that the tensor symbol on the right-hand side here denotes Day convolution for functors $\cat{F}^{\omega}_{\cat{C}} \to \spectra$, rather than $\cat{F}_{\cat{C}} \to \spectra$, so we cannot directly apply Theorem~\ref{thm:dn}. However, the functors $\der_1(-)(X_i)$ and $\der_n(-)(X_1,\dots,X_n)$ all preserve filtered colimits, so are equivalent to the left Kan extensions of their restrictions to $\cat{F}^{\omega}_{\cat{C}} \subseteq \cat{F}_{\cat{C}}$ (by \cite[5.3.5.8(2)]{lurie:2009}). It follows that the Day convolution calculated in the subcategory $\cat{F}^{\omega}_{\cat{C}}$ is equivalent to that calculated over $\cat{F}_{\cat{C}}$. We thus have, by Theorem~\ref{thm:dn},
\[ \Map_{\mathbb{I}_{\cat{C}}^{\otimes}}(X_1,\dots,X_n;Y) \homeq \Nat_{\cat{F}^{\omega}_{\cat{C}}}(\der_1(-)(Y), \der_n(-)(X_1,\dots,X_n)). \]
Since $\der_1(-)(Y)$ also preserves filtered colimits, a similar argument implies that in fact
\[ \Map_{\mathbb{I}_{\cat{C}}^{\otimes}}(X_1,\dots,X_n;Y) \homeq \Nat_{\cat{F}_{\cat{C}}}(\der_1(-)(Y), \der_n(-)(X_1,\dots,X_n)) \]
which yields the desired formula by Corollary~\ref{cor:dernI}.

In particular, the underlying $\infty$-category of $\mathbb{I}^{\otimes}_{\cat{C}}$ has mapping spectra
\[ \Map_{\mathbb{I}_{\cat{C}}}(X,Y) \homeq \der_1(I_{\cat{C}})(X;Y) \homeq \Map_{\spectra(\cat{C})}(Y,X) \]
so is equivalent to $\spectra(\cat{C})^{op}$. It also now follows that the $\infty$-operad $\mathbb{I}^{\otimes}_{\cat{C}}$ is stable.
\end{proof}

\begin{remark} \label{rem:compact}
It is sometimes convenient to restrict $\mathbb{I}^{\otimes}_{\cat{C}}$ to the \emph{small} $\infty$-operad $\check{\mathbb{I}}^{\otimes}_{\cat{C}} \subseteq \mathbb{I}^{\otimes}_{\cat{C}}$ whose underlying objects are the functors $\der_1(-)(X)$ for \emph{compact} objects $X$ in $\spectra(\cat{C})$. Since $\spectra(\cat{C})$ is compactly-generated by \cite[1.4.3.7]{lurie:2017}, those compact objects generate $\spectra(\cat{C})$ under filtered colimits. Moreover, the functor $\der_nI_{\cat{C}}$ preserves filtered colimits in each of its variables, and so the previous proposition shows that the $\infty$-operad $\mathbb{I}^{\otimes}_{\cat{C}}$ is determined by its restriction to these compact objects in a canonical way.
\end{remark}

\begin{remark} \label{rem:corep}
The $\infty$-operad $\mathbb{I}_{\cat{C}}^{\otimes}$ satisfies the additional property of being \emph{corepresentable} in the sense of \cite[6.2.4.3]{lurie:2017}, that is, the structure map $\mathbb{I}^{\otimes}_{\cat{C}} \to \Surj$ is a locally cocartesian fibration. In other language, we can think of $\mathbb{I}_{\cat{C}}^{\otimes}$ as encoding an \emph{oplax normal} symmetric monoidal structure on the underlying $\infty$-category $\spectra(\cat{C})^{op}$ or, equivalently, a \emph{lax normal} symmetric monoidal structure on $\spectra(\cat{C})$ in the sense of \cite{day/street:2003}.

More explicitly, this lax monoidal structure consists of the functors
\[ \Delta_n(I_{\cat{C}}): \spectra(\cat{C})^n \to \spectra(\cat{C}) \]
associated to the layers in the Taylor tower of $I_{\cat{C}}$, together with suitably compatible natural transformations
\[ \Delta_n(I_{\cat{C}})(X_1,\dots, \Delta_r(I_{\cat{C}})(Y_1,\dots,Y_r),\dots,X_n) \to \Delta_{n+r-1}(I_{\cat{C}})(X_1,\dots,Y_1,\dots,Y_r,\dots,X_n). \]
Such a structure is also referred to sometimes as a \emph{functor-operad} \cite[2.7]{yeakel:2020}.
\end{remark}

\begin{lemma} \label{lem:corep}
The $\infty$-operad map $\mathbb{I}_{\cat{C}}^{\otimes} \to \Surj$ is a locally cocartesian fibration.
\end{lemma}
\begin{proof}
This follows from \cite[Remark~6.2.4.5]{lurie:2017} and the natural equivalences
\[ \Map_{\mathbb{I}_{\cat{C}}^{\otimes}}(X_1,\dots,X_n;Y) \homeq \der_n(I_{\cat{C}})(X_1,\dots,X_n;Y) \homeq \Map_{\spectra(\cat{C})}(Y,\Delta_n(I_{\cat{C}})(X_1,\dots,X_n)) \]
of \ref{prop:sp-op} and \ref{def:deriv}.
\end{proof}

We have therefore proved the following result, which verifies Conjecture~6.3.0.13 of \cite{lurie:2017}.

\begin{proposition} \label{prop:corep}
Let $\cat{C}$ be a pointed compactly-generated $\infty$-category. Then there is a stable corepresentable $\infty$-operad $\mathbb{I}_{\cat{C}}^{\otimes}$ with underlying $\infty$-category $\spectra(\cat{C})^{op}$ whose corresponding lax monoidal structure consists of the symmetric multilinear functors
\[ \Delta_n(I_{\cat{C}}): \spectra(\cat{C})^n \to \spectra(\cat{C}) \]
associated to the Taylor tower of $I_{\cat{C}}$.
\end{proposition}

\begin{remark}
One would expect there to be a close relationship between the $\infty$-operad $\mathbb{I}_{\cat{C}}^{\otimes}$ constructed here and Lurie's $\infty$-operad $\spectra(\cat{C})^{\otimes}$ of \cite[6.3.0.14]{lurie:2017}. As described in \cite[6.3.0.17]{lurie:2017}, we expect these two $\infty$-operads to be Koszul dual, but as far as we know a theory of Koszul duality for (stable) $\infty$-operads has not yet been sufficiently developed to allow this conjecture to be verified.
\end{remark}

\section{Stable algebras over \texorpdfstring{$\infty$}{infinity}-operads} \label{sec:ex-AlgO}

We now turn to our main example: the case where $\cat{C}$ is an $\infty$-category of algebras over a (non-unital) stable $\infty$-operad. In particular, this covers the `classical' case of algebras over a (reduced) operad of spectra: for example, the reader may have in mind (non-unital) $A_\infty$- or $E_\infty$-ring spectra.

Let $\cat{O}^{\otimes}$ be a stable non-unital $\infty$-operad, and let $\Alg_{\cat{O}^{\otimes}}$ be the category of stable $\cat{O}^{\otimes}$-algebras defined in \ref{def:stable-algebra} below. It is a well-known slogan in Goodwillie calculus that the derivatives of the identity functor on $\Alg_{\cat{O}^{\otimes}}$ should be equivalent to $\cat{O}^{\otimes}$ itself. For example, in the monochromatic case a model for the Taylor tower for the identity on $\Alg_{\cat{O}^{\otimes}}$ is constructed by Pereira in~\cite{pereira:2013}, where the derivatives can be read off directly as the terms of the operad $\cat{O}^{\otimes}$. This tower was also studied by Harper and Hess in~\cite{harper/hess:2013}. The goal of this section is to improve that slogan to a version that takes the operad structures into account.

We are therefore interested in comparing the $\infty$-operad $\mathbb{I}_{\Alg_{\cat{O}^{\otimes}}}^{\otimes}$, given by applying Definition~\ref{def:I} to $\Alg_{\cat{O}^{\otimes}}$, with the $\infty$-operad $\cat{O}^{\otimes}$ itself. We will see, however, that the underlying $\infty$-categories of these two $\infty$-operads are not equivalent, thus precluding an actual equivalence of $\infty$-operads. Instead, we can identify $\mathbb{I}_{\Alg_{\cat{O}^{\otimes}}}^{\otimes}$ with an $\infty$-operad $\Pro(\cat{O})^{\otimes^{\ex}}$ whose underlying $\infty$-category is that of pro-objects in $\cat{O}$. There is a fully-faithful embedding of $\cat{O}^{\otimes}$ into $\Pro(\cat{O})^{\otimes^{\ex}}$, so we can identify $\cat{O}^{\otimes}$ with a full suboperad of $\mathbb{I}_{\Alg_{\cat{O}^{\otimes}}}^{\otimes}$.

Moreover, we show that the inclusion of this suboperad induces an equivalence
\[ \Alg_{\mathbb{I}_{\Alg_{\cat{O}^{\otimes}}}^{\otimes}} \weq \Alg_{\cat{O}^{\otimes}} \]
which we take to mean that $\mathbb{I}_{\Alg_{\cat{O}^{\otimes}}}^{\otimes}$ is \emph{Morita}-equivalent to $\cat{O}^{\otimes}$. This is our most precise version of the slogan mentioned at the beginning of this section.

\begin{definition} \label{def:stable-algebra}
Let $\cat{O}^{\otimes}$ be a small stable non-unital $\infty$-operad. A \emph{stable $\cat{O}^{\otimes}$-algebra} is a map of (non-unital) $\infty$-operads
\[ X: \cat{O}^{\otimes} \to \spectra^{\smsh} \]
such that the underlying functor $X: \cat{O} \to \spectra$ is exact. Let $\Alg_{\cat{O}^{\otimes}}$ denote the $\infty$-category of stable $\cat{O}^{\otimes}$-algebras, a full subcategory of $\Fun(\cat{O}^{\otimes},\spectra^{\smsh})$.
\end{definition}

\begin{example}
Let $\mathbf{O}$ be an ordinary reduced monochromatic operad in the symmetric monoidal model category of spectra. As ddescribed in Example~\ref{ex:finspec}, $\mathbf{O}$ corresponds to a stable $\infty$-operad $\cat{O}^{\otimes}$ whose underlying $\infty$-category is equivalent to $\spectra^{\omega}$, the $\infty$-category of finite spectra. The $\infty$-category $\Alg_{\cat{O}^{\otimes}}$ is then equivalent to that of fibrant-cofibrant objects in the projective model structure on the category of (non-unital) $\mathbf{O}$-algebras. For example, if $\mathbf{O}$ is the ordinary commutative operad, then $\Alg_{\cat{O}^{\otimes}}$ is equivalent to the $\infty$-category of (non-unital) $E_\infty$-ring spectra.
\end{example}

\begin{lemma} \label{lem:compgen}
For a small stable non-unital $\infty$-operad $\cat{O}^{\otimes}$, the $\infty$-category $\Alg_{\cat{O}^{\otimes}}$ is pointed and compactly-generated.
\end{lemma}
\begin{proof}
Let $\Alg^{\mathrm{un}}_{\cat{O}^{\otimes}}$ be the $\infty$-category of \emph{all} $\cat{O}^{\otimes}$-algebras in $\spectra$ (that is, the $\infty$-operad maps $\cat{O}^{\otimes} \to \spectra^{\smsh}$ with no restriction on exactness of the underlying functor). Then we have a pullback of $\infty$-categories:
\[ \begin{diagram} \dgARROWLENGTH=2em
	\node{\Alg_{\cat{O}^{\otimes}}} \arrow{s} \arrow{e} \node{\Alg^{\mathrm{un}}_{\cat{O}^{\otimes}}} \arrow{s,r}{U} \\
	\node{\Fun^{\ex}(\cat{O},\spectra)} \arrow{e,V} \node{\Fun(\cat{O},\spectra)} 
\end{diagram} \]
where $\Fun^{\ex}(\cat{O},\spectra) \subseteq \Fun(\cat{O},\spectra)$ is the subcategory of exact functors.

The $\infty$-category $\Alg^{\mathrm{un}}_{\cat{O}^{\otimes}}$ is compactly-generated by \cite[5.3.1.17]{lurie:2017}, and the forgetful functor $U$ has a left adjoint by the Adjoint Functor Theorem~\cite[5.5.2.9]{lurie:2009}, since limits and filtered colimits of $\cat{O}$-algebras are calculated objectwise by~\cite[3.2]{lurie:2017}. The $\infty$-category $\Fun^{\ex}(\cat{O},\spectra)$ is closed under all limits and colimits in $\Fun(\cat{O},\spectra)$ and the inclusion admits a left adjoint $D_{\ex}$ given by the Goodwillie excisive approximation (applied to the reduction of a functor $\cat{O} \to \spectra$). It follows from~\cite[5.5.7.3]{lurie:2009} that $\Fun^{\ex}(\cat{O},\spectra)$ is compactly-generated.

By~\cite[5.5.3.18]{lurie:2009}, the pullback diagram above is also a pullback in the $\infty$-category of presentable $\infty$-categories and right adjoint functors. It follows that the inclusion $\Alg_{\cat{O}^{\otimes}} \to \Alg^{\mathrm{un}}_{\cat{O}^{\otimes}}$ has a left adjoint, and so, by~\cite[5.5.7.3]{lurie:2009} again, $\Alg_{\cat{O}^{\otimes}}$ is compactly-generated. Finally, $\Alg_{\cat{O}^{\otimes}}$ is pointed with the constant trivial functor as a null object.
\end{proof}

\begin{remark} \label{rem:freealg}
The proof of Lemma~\ref{lem:compgen} also shows that the forgetful functor
\[ U: \Alg_{\cat{O}^{\otimes}} \to \Fun^{\ex}(\cat{O},\spectra) \]
has a left adjoint $F$, the \emph{free (stable) $\cat{O}^{\otimes}$-algebra} functor.
\end{remark}

We now wish to calculate the derivatives of the identity functor on $\Alg_{\cat{O}^{\otimes}}$, and the first step is to identify the stabilization of this $\infty$-category of algebras. The following result is a generalization of a theorem of Basterra and Mandell~\cite{basterra/mandell:2005}.

\begin{proposition} \label{prop:stable-alg}
Let $\cat{O}^{\otimes}$ be a small stable non-unital $\infty$-operad. Then there is an equivalence of $\infty$-categories
\[ \spectra(\Alg_{\cat{O}^{\otimes}}) \homeq \Fun^{\ex}(\cat{O},\spectra) \]
where $\cat{O}$ is the underlying $\infty$-category of $\cat{O}^{\otimes}$, and the right-hand side is the $\infty$-category of exact functors $\cat{O} \to \spectra$.
\end{proposition}
\begin{proof}
We apply a version of~\cite[7.3.4.7]{lurie:2017} to the stable $\cat{O}$-monoidal $\infty$-category
\[ \cat{O}^{\otimes} \times_{\Surj} \spectra^{\smsh} \to \cat{O}^{\otimes}. \]
That result is stated by Lurie for \emph{unital} $\infty$-operads, but a similar proof applies in the non-unital case. In particular, we apply~\cite[6.2.2.17]{lurie:2017} to the free-forgetful adjunction between $\Alg_{\cat{O}^{\otimes}}$ and $\Fun^{\ex}(\cat{O},\spectra)$.
\end{proof}

Proposition~\ref{prop:stable-alg} implies that $\mathbb{I}_{\Alg_{\cat{O}^{\otimes}}}^{\otimes}$ is an $\infty$-operad whose underlying $\infty$-category is equivalent to $\Fun^{\ex}(\cat{O},\spectra)^{op}$. That $\infty$-category, however, has another interpretation.

\begin{definition}
Let $\cat{C}$ be a small $\infty$-category that admits finite limits. Then the $\infty$-category of \emph{pro-objects} on $\cat{C}$ is the full subcategory
\[ \Pro(\cat{C}) \subseteq \Fun(\cat{C},\spaces)^{op} \]
consisting of those functors $\cat{C} \to \spaces$ that preserve finite limits. (See~\cite[5.3.5.4]{lurie:2009} for a dual description of ind-objects.) 
\end{definition}

\begin{lemma}
If $\cat{C}$ is a small stable $\infty$-category, then composition with $\Omega^\infty$ determines an equivalence
\[ \Fun^{\ex}(\cat{C},\spectra)^{op} \weq \Pro(\cat{C}). \]
In particular, for a small stable $\infty$-operad $\cat{O}^{\otimes}$, there is an equivalence
\[ \spectra(\Alg_{\cat{O}^{\otimes}})^{op} \homeq \Pro(\cat{O}). \]
\end{lemma}
\begin{proof}
A functor $\cat{C} \to \spaces$ that preserves finite limits is excisive and linear (since $\cat{C}$ is stable) and hence by~\cite[6.1.2.9]{lurie:2017} factors via $\Omega^\infty$ and an exact functor $\cat{C} \to \spectra$. Conversely, any functor that factors in this way preserves finite limits.
\end{proof}

\begin{example}
Let $\cat{O}^{\otimes}$ be a stable non-unital $\infty$-operad whose underlying $\infty$-category is that of finite spectra, corresponding to an ordinary operad $\mathbf{O}$ of spectra as in Example~\ref{ex:finspec}. Then $\spectra(\Alg_{\cat{O}^{\otimes}}) \homeq \spectra$, and $\mathbb{I}_{\Alg_{\cat{O}^{\otimes}}}^{\otimes}$ is a corepresentable stable non-unital $\infty$-operad whose underlying $\infty$-category is equivalent to $\Pro(\spectra^{\omega}) \homeq \spectra^{op}$. Such an object can also be identified with an ordinary operad $\mathbf{O}$.
\end{example}

The goal of the rest of this section is to identify the $\infty$-operad structure on $\Pro(\cat{O})$ that corresponds to the derivatives of the identity on $\Alg_{\cat{O}^{\otimes}}$. To describe that structure we use the monoidal envelope construction.

\begin{definition}[{\cite[2.2.4.1]{lurie:2017}}]
Let $\cat{O}^{\otimes} \to \Surj$ be a small non-unital $\infty$-operad. The \emph{symmetric monoidal envelope} of $\cat{O}^{\otimes}$ is the non-unital symmetric monoidal $\infty$-category
\[ \cat{O}^{\otimes}_{act} \to \Surj \]
where $\cat{O}^{\otimes}_{act} := \cat{O}^{\otimes} \times_{\Surj} \mathrm{Act}(\Surj)$ and $\mathrm{Act}(\Surj)$ is the category of active morphisms in $\Surj$ (where objects are active morphisms and morphisms are commutative squares).

The objects of the underlying $\infty$-category $\cat{O}_{act}$ can be identified with the objects of $\cat{O}^{\otimes}$, i.e. finite sequences of objects in $\cat{O}$. For each $n \geq 0$, we can identify $\cat{O}^n$ with a subcategory of $\cat{O}_{act}$, whose objects are the sequences of length $n$ and whose morphisms are those that cover the identity map on $\langle n \rangle$. Note that the empty sequence (the case $n = 0$) is disjoint from the rest of the $\infty$-category $\cat{O}_{act}$ because there are no active surjections between $\langle 0 \rangle$ and any other object of $\Surj$. 
\end{definition}

\begin{remark} \label{rem:free}
The symmetric monoidal $\infty$-category $\cat{O}_{act}^{\otimes}$ has the following universal property: for any symmetric monoidal $\infty$-category $\cat{D}^{\otimes}$, there is an equivalence between the $\infty$-category of $\infty$-operad maps $\cat{O}^{\otimes} \to \cat{D}^{\otimes}$ and that of symmetric monoidal functors $\cat{O}_{act}^{\otimes} \to \cat{D}^{\otimes}$.
\end{remark}

The $\infty$-operad structure on $\Pro(\cat{O})$ is based on a version of Day convolution for functors $\cat{O}_{act} \to \spectra$ that are \emph{exact} in the following sense.

\begin{definition} \label{def:gen-exact}
We say that $F: \cat{O}_{act} \to \spectra$ is \emph{exact} if, for each $n \geq 1$, the restriction $F|_{\cat{O}^n}: \cat{O}^n \to \spectra$ is exact in each variable, and $F() \homeq *$. (That is, on the empty sequence in $\cat{O}_{act}$, $F$ is null.) Let $\Fun^{\ex}(\cat{O}_{act},\spectra)$ denote the full subcategory of $\Fun(\cat{O}_{act},\spectra)$ consisting of the exact functors.
\end{definition}

\begin{lemma} \label{lem:exact-adjoint}
The inclusion $\Fun^{\ex}(\cat{O}_{act},\spectra) \to \Fun(\cat{O}_{act},\spectra)$ admits a left adjoint $D_{\ex}$.
\end{lemma}
\begin{proof}
We show that $\Fun^{\ex}(\cat{O}_{act},\spectra)$ is a strongly reflective subcategory of $\Fun(\cat{O}_{act},\spectra)$ in the sense of~\cite[5.5.4.17]{lurie:2009}. Now $\Fun^{\ex}(\cat{O}_{act},\spectra)$ is the intersection of a (small) collection of subcategories of $\Fun(\cat{O}_{act},\spectra)$, each determined by the condition that a particular restriction to $\cat{O}$ is exact. To be precise, there is such a restriction for each $n \geq 1$, each $1 \leq i \leq n$, and each sequence of objects $J_1,\dots,J_{i-1},J_{i+1},\dots,J_n \in \cat{O}$. The corresponding restriction is along the inclusion
\[ \cat{O} \into \cat{O}_{act}; \quad I \mapsto (J_1,\dots,J_{i-1},I,J_{i+1},\dots,J_n). \]
The subcategory of functors that restrict to an exact functor along this inclusion is given by a pullback (of $\infty$-categories) of the form
\[ \begin{diagram} \dgARROWLENGTH=2em
	\node{\Fun^{\ex}_{\alpha}(\cat{O}_{act},\spectra)} \arrow{s} \arrow{e} \node{\Fun(\cat{O}_{act},\spectra)} \arrow{s,r}{r_{\alpha}} \\
	\node{\Fun^{\ex}(\cat{O},\spectra)} \arrow{e} \node{\Fun(\cat{O},\spectra)}
\end{diagram} \]
where $r_{\alpha}$ is restriction along the relevant inclusion $\cat{O} \into \cat{O}_{act}$. By the Adjoint Functor Theorem, the map $r_{\alpha}$ has a left adjoint, as does the bottom inclusion; see the proof of Lemma~\ref{lem:compgen}. Thus, by~\cite[5.5.3.18]{lurie:2009}, the above pullback can be done in the $\infty$-category of presentable $\infty$-categories and right adjoints. In particular, $\Fun^{\ex}_{\alpha}(\cat{O}_{act},\spectra)$ is a strongly reflective subcategory of $\Fun(\cat{O}_{act},\spectra)$. Therefore, by~\cite[5.5.4.18]{lurie:2009},
\[ \Fun^{\ex}(\cat{O}_{act},\spectra) = \bigcap_{\alpha} \Fun^{\ex}_{\alpha}(\cat{O}_{act},\spectra) \]
is also strongly reflective.
\end{proof}

\begin{definition} \label{def:gen-Day}
Applying Construction~\ref{cons:day-infty} to $\cat{O}_{act}$ we obtain a (non-unital) stable symmetric monoidal $\infty$-category
\[ \Fun(\cat{O}_{act},\spectra)^{\otimes} \]
with symmetric monoidal structure given by Day convolution. Let
\[ \Fun^{\ex}(\cat{O}_{act},\spectra)^{\otimes^{\ex}} \subseteq \Fun(\cat{O}_{act},\spectra)^{\otimes} \]
be the suboperad generated by the exact functors.
\end{definition}

\begin{lemma} \label{lem:exact-day}
There is a symmetric monoidal functor
\[ D_{\ex}: \Fun(\cat{O}_{act},\spectra)^{\otimes} \to \Fun^{\ex}(\cat{O}_{act},\spectra)^{\otimes^{\ex}} \]
whose underlying functor is the left adjoint to the inclusion.
\end{lemma}
\begin{proof}
The existence of the left adjoint $D_{\ex}$ was demonstrated in Lemma~\ref{lem:exact-adjoint}. We now apply~\cite[2.2.1.9]{lurie:2017} to show that $D_{\ex}$ extends to a symmetric monoidal functor as shown. It is sufficient to show that the Day convolution of functors $\cat{O}_{act} \to \spectra$ preserves $D_{\ex}$-equivalences. That is, suppose $F_1 \to F_2$ and $G_1 \to G_2$ each become equivalences after applying $D_{\ex}$. We have to show that $F_1 \otimes F_2 \to G_1 \otimes G_2$ has the same property.

It is sufficient to show that for any exact $H: \cat{O}_{act} \to \spectra$, the induced map
\[ \Nat_{\cat{O}_{act}}(D_{\ex}(F_2 \otimes G_2),H) \to \Nat_{\cat{O}_{act}}(D_{\ex}(F_1 \otimes G_1),H) \]
is an equivalence. Using the universal properties of $D_{\ex}$ and the Day convolution, this map can be written in the form
\[ \Nat_{\cat{O}_{act} \times \cat{O}_{act}}(F_2(-) \smsh G_2(-),H(-,-)) \to \Nat_{\cat{O}_{act} \times \cat{O}_{act}}(F_1(-) \smsh G_1(-),H(-,-)), \]
or equivalently
\[ \Nat_{\cat{O}_{act}}(F_2(-), \Nat_{\cat{O}_{act}}(G_2(-),H(-,-))) \to \Nat_{\cat{O}_{act}}(F_1(-), \Nat_{\cat{O}_{act}}(G_1(-),H(-,-))) \]
where $H(-,-)$ denotes $H$ applied to the concatenation of its two variables (i.e. the monoidal product on $\cat{O}_{act}$). Since $H$ is exact, it is exact in each of those variables. The hypotheses on the maps $F_1 \to F_2$ and $G_1 \to G_2$ therefore imply that the above map is an equivalence.
\end{proof}

\begin{remark}
Lemma~\ref{lem:exact-day} tells us that the exact functors $\cat{O}_{act} \to \spectra$ have a symmetric monoidal structure $\otimes^{\ex}$ which we call the \emph{exact Day convolution}. This structure satisfies the same universal property as ordinary Day convolution (Definition~\ref{def:Day}) but with respect only to exact functors. It follows from~\ref{lem:exact-day} that $\otimes^{\ex}$ is given by taking the ordinary Day convolution and then applying $D_{\ex}$. In fact, we believe that the ordinary convolution of exact functors is already exact, so that this last step is unnecessary, and the inclusion of the suboperad in Lemma~\ref{lem:exact-day} is a symmetric monoidal functor. However, we do not need this fact here, so we do not include a proof.
\end{remark}

We now use the exact Day convolution to define an $\infty$-operad structure on $\Pro(\cat{O})$.

\begin{lemma} \label{lem:pro}
There is a fully-faithful embedding of stable $\infty$-categories
\[ \Pro(\cat{O}) \homeq \Fun^{\ex}(\cat{O},\spectra)^{op} \arrow{e,t}{\ell^{op}} \Fun^{\ex}(\cat{O}_{act},\spectra)^{op} \]
where $\ell$ is left Kan extension along the inclusion $i: \cat{O} \into \cat{O}_{act}$. The essential image of $\ell$ consists of those functors $F: \cat{O}_{act} \to \spectra$ such that
\begin{enumerate}
  \item $F(I_1,\dots,I_n) \homeq *$ for any $I_1,\dots,I_n \in \cat{O}$ where $n \neq 1$;
  \item the restriction of $F$ to $\cat{O} \subseteq \cat{O}_{act}$ is an exact functor $\cat{O} \to \spectra$.
\end{enumerate}
\end{lemma}
\begin{proof}
Since the embedding $i: \cat{O} \into \cat{O}_{act}$ is fully-faithful, the Kan extension $\ell$ is too, and since there are no morphisms in $\cat{O}_{act}$ of the form $I \to (J_1,\dots,J_n)$ for $n \neq 1$, the functor $\ell(F)$ satisfies condition (i) and its restriction to $\cat{O}$ is given by $F$ itself.	
\end{proof}

\begin{definition} \label{def:O-Day}
Let $\Fun^{\ex}(\cat{O}_{act},\spectra)^{op,\otimes^{\ex}}$ be the symmetric monoidal $\infty$-category given by applying the opposite symmetric monoidal structure of Construction~\ref{cons:symmon-op} to the symmetric monoidal $\infty$-category $\Fun^{\ex}(\cat{O}_{act},\spectra)^{\otimes^{\ex}}$ of Definition~\ref{def:gen-Day}. 

Let $\Pro(\cat{O})^{\otimes^{\ex}}$ be the full suboperad of $\Fun^{\ex}(\cat{O}_{act},\spectra)^{op,\otimes^{\ex}}$ generated by those functors $F: \cat{O}_{act} \to \spectra$ that satisfy conditions (i) and (ii) of Lemma~\ref{lem:pro}. This suboperad has underlying $\infty$-category equivalent (via the map $\ell^{op}$ of \ref{lem:pro}) to $\Pro(\cat{O})$.
\end{definition}

\begin{remark} \label{rem:pro}
We can describe the $\infty$-operad structure on $\Pro(\cat{O})$ in the following way. Suppose we have exact functors $X_1,\dots,X_n,Y: \cat{O} \to \spectra$. Then the multi-morphism spectrum
\[ \Hom_{\Pro(\cat{O})^{\otimes^{\ex}}}(X_1,\dots,X_n;Y) \]
is equivalent to the spectrum of natural transformations (of functors $\cat{O}_{act} \to \spectra$) of the form
\[ \ell Y \to \ell X_1 \otimes^{\ex} \dots \otimes^{\ex} \ell X_n \]
where $\otimes^{\ex}$ denotes the exact Day convolution. Equivalently, this is the spectrum of natural transformations (of functors $\cat{O} \to \spectra)$ of the form
\[ Y \to r(\ell X_1 \otimes^{\ex} \dots \otimes^{\ex} \ell X_n) \]
where $r$ is restriction along $i: \cat{O} \into \cat{O}_{act}$. In other words, the $\infty$-operad structure on $\Pro(\cat{O})$ is corepresented by the construction
\[ (X_1,\dots,X_n) \mapsto r(\ell X_1 \otimes^{\ex} \dots \otimes^{\ex} \ell X_n) \]
which we can write, using coend notation for the Day convolution, as the functor
\[ J \mapsto \Map_{\cat{O}^{\otimes}}(I_1,\dots,I_n;J) \smsh_{I_1,\dots,I_n \in \cat{O}} X_1(I_1) \smsh \dots \smsh X_n(I_n). \]
\end{remark}

The main result of this section is that the $\infty$-operad $\mathbb{I}_{\Alg_{\cat{O}^{\otimes}}}^{\otimes}$ is equivalent to $\Pro(\cat{O})^{\otimes^{\ex}}$. In order to make the comparison, we require the following construction.

\begin{definition} \label{def:ev}
Let $\cat{O}^{\otimes}$ be a small stable non-unital $\infty$-operad. Using Construction~\ref{cons:objectwise} we can produce a map of $\infty$-operads
\[ \ev: \cat{O}^{\otimes} \to \cat{F}_{\Alg_{\cat{O}^{\otimes}}}^{\smsh} \subseteq \Fun(\Alg_{\cat{O}^{\otimes}},\spectra)^{\smsh} \]
given, on underlying $\infty$-categories, by
\[ I \mapsto \ev_I \]
where $\ev_I: \Alg_{\cat{O}^{\otimes}} \to \spectra$ is the functor that evaluates an $\cat{O}^{\otimes}$-algebra $X$ at the object $I \in \cat{O}$. By the universal property of the monoidal envelope, the map $\ev$ extends canonically to a symmetric monoidal functor
\[ \evb : \cat{O}_{act}^{\otimes} \to \cat{F}_{\Alg_{\cat{O}^{\otimes}}}^{\smsh} \]
given by
\[ (I_1,\dots,I_k) \mapsto \ev_{I_1} \smsh \dots \smsh \ev_{I_k}. \]
Note that the proof of~\cite[2.2.4.9]{lurie:2017} demonstrates that $\evb$ is a (relative over $\Surj$) left Kan extension of $\ev$ along the inclusion of $\infty$-operads $\cat{O}^{\otimes} \to \cat{O}^{\otimes}_{act}$.

Finally, restriction along $\evb$ determines a map of $\infty$-operads (or lax symmetric monoidal functor)
\[ \evb^*: \Fun(\cat{F}_{\Alg_{\cat{O}^{\otimes}}},\spectra)^{\otimes} \to \Fun(\cat{O}_{act},\spectra)^{\otimes} \]
by~\cite[3.8]{nikolaus:2016}.
\end{definition}

\begin{remark}
There is one technical wrinkle in the previous definition that we have to be careful with. Recall that the $\infty$-category $\cat{F}_{\Alg_{\cat{O}^{\otimes}}}$ is not small, and so in the construction of $\mathbb{I}_{\Alg_{\cat{O}^{\otimes}}}$ we replaced it with the small subcategory $\cat{F}^{\omega}_{\Alg_{\cat{O}^{\otimes}}}$ of compact objects. Since it is unclear whether the functors $\ev_I$ are compact, we can simply add them in the following way.

Let $\tilde{\cat{F}}^{\omega}_{\Alg_{\cat{O}^{\otimes}}}$ be the (essentially) small full subcategory of $\cat{F}_{\Alg_{\cat{O}^{\otimes}}}$ obtained from $\cat{F}^{\omega}_{\Alg_{\cat{O}^{\otimes}}}$ by adjoining the objects $\ev_{I}$ for all $I \in \cat{O}$, and then taking closure under the objectwise smash product. The functor $\evb$ takes values in this (essentially small) symmetric monoidal $\infty$-category, and we get an associated restriction functor
\begin{equation} \label{eq:ev} \evb^*: \Fun(\tilde{\cat{F}}^{\omega}_{\Alg_{\cat{O}^{\otimes}}},\spectra)^{\otimes} \to \Fun(\cat{O}_{act},\spectra)^{\otimes}. \end{equation}
Since $\tilde{\cat{F}}^{\omega}_{\Alg_{\cat{O}^{\otimes}}}$ generates $\cat{F}_{\Alg_{\cat{O}^{\otimes}}}$ under filtered colimits, the arguments of Proposition~\ref{prop:sp-op} imply that $\mathbb{I}_{\Alg_{\cat{O}^{\otimes}}}^{\otimes}$ can be identified with the suboperad of $\Fun(\tilde{\cat{F}}^{\omega}_{\Alg_{\cat{O}^{\otimes}}},\spectra)^{op,\otimes}$ generated by objects of the form $\der_1(-)(X)$ for $X \in \spectra(\Alg_{\cat{O}^{\otimes}})$.
\end{remark}

Our goal now is to show that $\evb^*$ determines an equivalence of $\infty$-operads between $\mathbb{I}_{\Alg_{\cat{O}^{\otimes}}}^{\otimes}$ and $\Pro(\cat{O})^{\otimes^{\ex}}$ by taking opposites and restricting to the relevant suboperads. The following result contains the main calculation that establishes this equivalence.

\begin{proposition} \label{prop:ev}
	Let $\mathbb{J}^{\otimes} \subseteq \Fun(\tilde{\cat{F}}^{\omega}_{\Alg_{\cat{O}^{\otimes}}},\spectra)^{\otimes}$ be the monoidal subcategory consisting of those objects of the form
	\[ \der_n(-)(X_1,\dots,X_n) \]
	for $X_1,\dots,X_n \in \spectra(\Alg_{\cat{O}^{\otimes}}) \homeq \Fun^{\ex}(\cat{O},\spectra)$.  (This subcategory is closed under Day convolution by Theorem~\ref{thm:dn}.) Then the map of $\infty$-operads (\ref{eq:ev}) restricts to a symmetric monoidal functor
	\[ \evb^*: \mathbb{J}^{\otimes} \to \Fun^{\ex}(\cat{O}_{act},\spectra)^{\otimes^{\ex}}. \]
\end{proposition}
\begin{proof}
First note that
\[ \evb^*(\der_n(-)(X_1,\dots,X_n)) \homeq \der_n(\evb)(X_1,\dots,X_n) \]
is an exact functor $\cat{O}_{act} \to \spectra$ because the underlying functor of a stable $\cat{O}^{\otimes}$-algebra is exact, the smash product preserves colimits in each variable, and $\der_n(-)(X_1,\dots,X_n)$ preserves colimits. So $\evb^*$ takes values in the $\infty$-operad $\Fun^{\ex}(\cat{O}_{act},\spectra)^{\otimes^{\ex}}$ as claimed.

We then have to show that the lax symmetric monoidal functor $\evb^*$ is, in fact, symmetric monoidal. The lax monoidal structure determines maps
\[ \der_1(\evb)(Y_1) \otimes^{\ex} \dots \otimes^{\ex} \der_1(\evb)(Y_k) \to (\der_1(-)(Y_1) \otimes \dots \otimes \der_1(-)(Y_k)(\evb) \]
where $\otimes^{\ex}$ is the exact Day convolution of Lemma~\ref{lem:exact-day}. By Theorem~\ref{thm:dn} we can write this map as
\begin{equation} \label{eq:der-ev} \der_1(\evb)(Y_1) \otimes^{\ex} \dots \otimes^{\ex} \der_1(\evb)(Y_k) \to \der_k(\evb)(Y_1,\dots,Y_k), \end{equation}
and it is sufficient to show that each such map is an equivalence. 

Changing perspective slightly, consider the functor
\[ \evb^{\sharp}: \Alg_{\cat{O}^{\otimes}} \to \Fun^{\ex}(\cat{O}_{act},\spectra); \quad X \mapsto \evb(-)(X). \]
Since limits in the $\infty$-category $\Fun^{\ex}(\cat{O}_{act},\spectra)$ are calculated objectwise, the Taylor tower of $\evb^{\sharp}$ is the $\cat{O}_{act}$-indexed diagram consisting of the Taylor towers of each of the functors $\ev_{J_1} \smsh \dots \smsh \ev_{J_n}: \Alg_{\cat{O}^{\otimes}} \to \spectra$. In other words we have a natural equivalence, of functors $\cat{O}_{act} \to \spectra$
\[ \Delta_k(\evb^{\sharp})(Y_1,\dots,Y_k) \homeq \Delta_k(\evb)(Y_1,\dots,Y_k) \homeq \der_k(\evb)(Y_1,\dots,Y_k). \]

Recall from Definition~\ref{def:ev} that $\evb$ is given by an operadic left Kan extension along the inclusion of $\infty$-operads $\cat{O}^{\otimes} \to \cat{O}^{\otimes}_{act}$. In other words, we can write $\evb^{\sharp}$ as the composite
\[ \Alg_{\cat{O}^{\otimes}} \arrow{e,t}{\ell^{\otimes}} \Alg_{\cat{O}_{act}^{\otimes}} \arrow{e,t}{U'} \Fun^{\ex}(\cat{O}_{act},\spectra) \]
of that Kan extension $\ell^{\otimes}$ with the forgetful functor $U'$.

We can interpret the middle $\infty$-category in the above composite as follows. By the universal property of Day convolution~\cite[2.2.6.9]{lurie:2017}, there is an equivalence
\begin{equation} \label{eq:CAlg} \Alg_{\cat{O}_{act}^{\otimes}} \homeq \CAlg(\Fun^{\ex}(\cat{O}_{act},\spectra)^{\otimes^{\ex}})  \end{equation}
where the right-hand side is the $\infty$-category of (non-unital) commutative algebras in the symmetric monoidal $\infty$-category $\Fun^{\ex}(\cat{O}_{act},\spectra)^{\otimes^{\ex}}$. (To be precise, that universal property gives an equivalence at the level of \emph{all} $\cat{O}$-algebras, not just those that are stable, and for the ordinary, not exact, Day convolution. However, the monoidal adjunction of Lemma~\ref{lem:exact-day} determines an equivalence between the $\infty$-categories of (i) those commutative algebras whose underlying functor is exact, and (ii) the commutative algebras for the exact Day convolution.) 

Now consider the following diagram of adjunctions
\begin{equation} \label{eq:evb} \begin{tikzcd}
	{\Alg_{\cat{O}^{\otimes}}} && {\CAlg(\Fun^{\ex}(\cat{O}_{act},\spectra)^{\otimes^{\ex}})} \\
	\\
	{\Fun^{\ex}(\cat{O},\spectra)} && {\Fun^{\ex}(\cat{O}_{act},\spectra)}
	\arrow[""{name=0, anchor=center, inner sep=0}, "{\ell^{\otimes}}", shift left=2, from=1-1, to=1-3]
	\arrow[""{name=2, anchor=center, inner sep=0}, "U", shift left=2, from=1-1, to=3-1]
	\arrow[""{name=3, anchor=center, inner sep=0}, "\ell", shift left=2, from=3-1, to=3-3]
	\arrow[""{name=4, anchor=center, inner sep=0}, "r", shift left=2, from=3-3, to=3-1]	
	\arrow[""{name=6, anchor=center, inner sep=0}, "U'", shift left=2, from=1-3, to=3-3]
	\arrow[""{name=7, anchor=center, inner sep=0}, "{r^{\otimes}}", shift left=2, from=1-3, to=1-1]
	\arrow[""{name=9, anchor=center, inner sep=0}, "F", shift left=2, from=3-1, to=1-1]
	\arrow[""{name=10, anchor=center, inner sep=0}, "F'", shift left=2, from=3-3, to=1-3]
	\arrow["\dashv"{anchor=center, rotate=0}, draw=none, from=2, to=9]
	\arrow["\dashv"{anchor=center, rotate=0}, draw=none, from=10, to=6]
	\arrow["\dashv"{anchor=center, rotate=-90}, draw=none, from=0, to=7]
	\arrow["\dashv"{anchor=center, rotate=-90}, draw=none, from=3, to=4]
\end{tikzcd} \end{equation}
where:
\begin{enumerate}
	\item $\ell^{\otimes}$ is the (relative over $\Surj$) left Kan extension along $i: \cat{O}^{\otimes} \into \cat{O}^{\otimes}_{act}$, combined with the equivalence (\ref{eq:CAlg}), and $r^{\otimes}$ is its right adjoint, the corresponding restriction;
	\item $(F,U)$ is the free-forgetful adjunction of Remark~\ref{rem:freealg};
	\item $(F',U')$ is the free-forgetful adjunction for commutative algebras in the symmetric monoidal $\infty$-category $\Fun^{\ex}(\cat{O}_{act},\spectra)^{\otimes^{\ex}}$;
	\item the map $l$ is as in Lemma~\ref{lem:pro} with right adjoint $r$ the restriction along $\cat{O} \into \cat{O}_{act}$.
\end{enumerate}
Noting that the diagram of right adjoints commutes, it follows that the diagram of left adjoints also commutes. Also, recall that we have identified the functor $\evb^{\sharp}$ with the composite $U'\ell^{\otimes}$

Our goal is therefore to calculate $\Delta_k(U'\ell^{\otimes})$, or equivalently $\Delta_k(U'\ell^{\otimes}F)$ since the derivatives of a functor on $\cat{O}$-algebras can be identified on free $\cat{O}$-algebras. We thus have
\[ \der_k(\evb)(Y_1,\dots,Y_k) \homeq \Delta_k(U'F'\ell)(Y_1,\dots,Y_k) \homeq \Delta_k(U'F')(\ell Y_1,\dots, \ell Y_k) \]
where the last equivalence follows from the fact that $\ell$ is a left adjoint and~\cite[6.1.1.30]{lurie:2017}.

We have now related the derivatives of $\evb^{\sharp}$ to the derivatives of the free commutative algebra functor in the symmetric monoidal stable $\infty$-category $\Fun^{\ex}(\cat{O}_{act},\spectra)^{\otimes^{\ex}}$ which can be calculated directly. That free algebra functor is given by
\[ U'F'(X) \homeq \bigvee_{n \geq 1} X^{\otimes^{\ex}n}_{\Sigma_n} \]
so has Taylor tower layers
\[ D_k(U'F')(X) \homeq X^{\otimes^{\ex}k}_{\Sigma_k} \]
and derivatives
\[ \Delta_k(U'F')(X_1,\dots,X_k) \homeq X_1 \otimes^{\ex} \dots \otimes^{\ex} X_k. \]
Taking $X_i = \ell Y_i$, we therefore have
\[ \Delta_k(\evb)(Y_1,\dots,Y_k) \homeq (\ell Y_1) \otimes^{\ex} \dots \otimes^{\ex} (\ell Y_k) \]
from which the desired equivalence (\ref{eq:der-ev}) follows.
\end{proof}

\begin{remark}
An informal, but perhaps more intuitive, proof of Proposition~\ref{prop:ev} can be given by providing coend formulas for the various Kan extensions and Day convolutions involved. In particular, we can calculate the derivatives
\[ \der_k(\ev_J)(Y_1,\dots,Y_k) \homeq \Map_{\cat{O}^{\otimes}}(I_1,\dots,I_k;J) \smsh_{I_1,\dots,I_k \in \cat{O}} Y_1(I_1) \smsh \dots \smsh Y_k(I_k). \]
A standard computation for the derivatives of a smash product of spectrum-valued functors such as $\der_k(\ev_{J_1} \smsh \dots \smsh \ev_{J_n})(Y_1,\dots,Y_k)$ then gives
\[ \left[ \bigvee_{\un{k} \epi \un{n}} \Smsh_{i = 1}^{n} \Map_{\cat{O}^{\otimes}}(I_{i_1},\dots,I_{i_{k_i}};J_i) \right] \smsh_{I_1,\dots,I_k \in \cat{O}} Y_1(I_1) \smsh \dots \smsh Y_k(I_k) \]
from which the equivalence (\ref{eq:der-ev}) can be recovered by directly calculating the relevant (exact) Day convolution.
\end{remark}

We can now prove the main result of this section.

\begin{theorem} \label{thm:O}
There is a fully-faithful embedding of $\infty$-operads
\[ \mathbb{I}_{\Alg_{\cat{O}^{\otimes}}}^{\otimes} \into \Fun^{\ex}(\cat{O}_{act},\spectra)^{op,\otimes^{\ex}} \]
whose essential image is $\Pro(\cat{O})^{\otimes^{\ex}}$.
\end{theorem}
\begin{proof}
Applying the opposite symmetric monoidal construction to the map of Proposition~\ref{prop:ev}, we get a symmetric monoidal functor $(\evb)^*: \mathbb{J}^{op,\otimes} \to \Fun^{\ex}(\cat{O}_{act},\spectra)^{op,\otimes^{\ex}}$. Restricting to the suboperad $\mathbb{I}_{\Alg_{\cat{O}^{\otimes}}}^{\otimes} \subseteq \mathbb{J}^{op,\otimes}$, we obtain the required map of $\infty$-operads. Identification of the essential image of that map with $\Pro(\cat{O})$ follows from Lemma~\ref{lem:pro} and the calculation of $\der_1(\evb)(Y)$ in the proof of Proposition~\ref{prop:ev}.

It remains to show that $\evb^*$ induces equivalences on multi-mapping spectra. By Theorem~\ref{thm:dn}, it is enough to show that restriction along $\evb$ induces equivalences
\begin{equation} \label{eq:equiv} \begin{diagram} \node{\Nat_{\tilde{\cat{F}}^{\omega}_{\Alg_{\cat{O}^{\otimes}}}}(\der_1(-)(Y), \der_n(-)(X_1,\dots,X_n))} \arrow{s,l}{\evb^*} \\ \node{\Nat_{\cat{O}_{act}}(\der_1(\evb)(Y), \der_n(\evb)(X_1,\dots,X_n))}
\end{diagram}
\end{equation}
for all $X_1,\dots,X_n,Y \in \spectra(\Alg_{\cat{O}^{\otimes}})$.

Recall from the proof of Proposition~\ref{prop:ev} that $\der_1(\evb)(Y) \homeq \ell(Y)$, so we can write the target of the map $\evb^*$ as a natural transformation object over $\cat{O}$. Precomposing with the equivalence of Theorem~\ref{thm:dF}, we are reduced to showing that evaluation determines an equivalence
\[ \Map_{\spectra(\Alg_{\cat{O}^{\otimes}})} (Y,\Delta_nI_{\Alg_{\cat{O}^{\otimes}}}(X_1,\dots,X_n)) \to \Nat_{\cat{O}}(Y, \der_n(\ev)(X_1,\dots,X_n)). \]
i.e. that the canonical maps
\[ \ev_J\Delta_nI_{\Alg_{\cat{O}^{\otimes}}} \to \Delta_n(\ev_J) = \der_n(\ev_J) \]
are equivalences, which follows from the fact that $\ev_J$ preserves limits.
\end{proof}

Having identified the $\infty$-operad $\mathbb{I}^{\otimes}_{\Alg_{\cat{O}^{\otimes}}}$ with $\Pro(\cat{O})^{\otimes^{\ex}}$, we now wish to relate this calculation to $\cat{O}^{\otimes}$ itself and justify the slogan mentioned at the beginning of this section. To do this, we show that $\cat{O}^{\otimes}$ embeds fully-faithfully, via a stable Yoneda map, into $\Pro(\cat{O})^{\otimes^{\ex}}$, and hence is equivalent to a full suboperad of $\mathbb{I}_{\Alg_{\cat{O}^{\otimes}}}^{\otimes}$.

\begin{definition}
Glasman shows in \cite[Sec. 3]{glasman:2016} that a symmetric monoidal $\infty$-category such as $\cat{O}_{act}^{op,\otimes}$ admits a multiplicative Yoneda embedding, that is, a fully-faithful symmetric monoidal map
\[ Y_{act}: \cat{O}_{act}^{op,\otimes} \to \Fun(\cat{O}_{act},\spaces)^{\otimes}. \]
By~\cite[3.7]{nikolaus:2016}, the suspension spectrum functor $\Sigma^\infty_+: \spaces \to \spectra$ induces a symmetric monoidal functor
\[ \Sigma^\infty_+: \Fun(\cat{O}_{act},\spaces)^{\otimes} \to \Fun(\cat{O}_{act},\spectra)^{\otimes}. \]
Composing $Y_{act}$ with $\Sigma^\infty_+$, and with the symmetric monoidal functor $D_{\ex}$ of Lemma~\ref{lem:exact-day}, and then taking opposites, we get a symmetric monoidal functor, the \emph{stable Yoneda embedding}
\[ Y: \cat{O}_{act}^{\otimes} \to \Fun^{\ex}(\cat{O}_{act},\spectra)^{op,\otimes^{\ex}}. \]
\end{definition}

\begin{lemma} \label{lem:O}
The restriction of $Y$ to $\cat{O}^{\otimes} \subseteq \cat{O}_{act}^{\otimes}$ determines a fully-faithful embedding of $\infty$-operads
\[ Y: \cat{O}^{\otimes} \to \Pro(\cat{O})^{\otimes^{\ex}}; \quad I \mapsto \Map_{\cat{O}}(I,-). \]
\end{lemma}
\begin{proof}
By construction, the value of $Y$ on an object $I \in \cat{O}$ is given by
\[ D_{\ex}(\Sigma^\infty_+\Hom_{\cat{O}_{act}}(I,-)). \]
Since there are no maps in $\cat{O}_{act}$ from $I$ to $(J_1,\dots,J_n)$ with $n \neq 1$, we have
\[ \Sigma^\infty_+\Hom_{\cat{O}_{act}}(I,(J_1,\dots,J_n)) \homeq \begin{cases} \Sigma^\infty_+\Omega^\infty \Map_{\cat{O}}(I,J_1) & \text{if $n = 1$}; \\ \quad \quad \quad * & \text{otherwise}. \end{cases} \]
Since $D_{\ex}(\Sigma^\infty_+ \Omega^\infty \Map_{\cat{O}}(I,-)) \homeq \Map_{\cat{O}}(I,-)$, we deduce the desired formula for $Y$ and see that it takes values in $\Pro(\cat{O})$.

It remains to show that $Y$ induces equivalences of multi-morphism spaces
\[ \Hom_{\cat{O}^{\otimes}}(I_1,\dots,I_n;J) \to \Hom_{\Pro(\cat{O})^{\otimes^{\ex}}}(\Map_{\cat{O}}(I_1,-),\dots,\Map_{\cat{O}}(I_n,-); \Map_{\cat{O}}(J,-)). \]
This claim follows from the calculation of the operad structure on $\Pro(\cat{O})^{\otimes^{\ex}}$ given in Remark~\ref{rem:pro}.
\end{proof}

In combination with Theorem~\ref{thm:O}, Lemma~\ref{lem:O} allows us to identify $\cat{O}^{\otimes}$ with a suboperad of $\mathbb{I}_{\Alg_{\cat{O}^{\otimes}}}^{\otimes}$.

\begin{theorem} \label{thm:O-sub}
Let $\cat{O}^{\otimes}$ be a small stable non-unital $\infty$-operad. The suboperad of $\mathbb{I}_{\Alg_{\cat{O}^{\otimes}}}^{\otimes}$ generated by the representable pro-objects $\Map_{\cat{O}}(I,-) \in \spectra(\Alg_{\cat{O}^{\otimes}})$ is equivalent to $\cat{O}^{\otimes}$.
\end{theorem}

In Remark~\ref{rem:compact} we noted that sometimes it makes sense to restrict to the small $\infty$-operad $\check{\mathbb{I}}^{\otimes}_{\cat{C}} \subseteq \mathbb{I}^{\otimes}_{\cat{C}}$ generated by the compact objects of $\spectra(\cat{C})$. If we do that with $\cat{C} = \Alg_{\cat{O}^{\otimes}}$ then we get the small stable $\infty$-operad with underlying $\infty$-category
\[ \check{\mathbb{I}}_{\Alg_{\cat{O}^{\otimes}}} \homeq (\Ind(\cat{O}^{op})^{\omega})^{op} \]
which, by \cite[5.4.2.4]{lurie:2009}, is an idempotent-completion of $\cat{O}$. In particular, we have the following.

\begin{corollary}
Let $\cat{O}^{\otimes}$ be a small stable $\infty$-operad with $\cat{O}$ idempotent-complete. Then there is an equivalence of $\infty$-operads
\[ \check{\mathbb{I}}^{\otimes}_{\Alg_{\cat{O}^{\otimes}}} \homeq \cat{O}^{\otimes}. \]
\end{corollary}

More generally, we can see that $\check{\mathbb{I}}^{\otimes}_{\Alg_{\cat{O}^{\otimes}}}$ is \emph{Morita}-equivalent to $\cat{O}^{\otimes}$, i.e. that the inclusion induces an equivalence between the corresponding $\infty$-category of (stable) algebras.

\begin{lemma} \label{lem:alg-res}
Let $i: \cat{A}^{\otimes} \to \cat{B}^{\otimes}$ be a fully-faithful map of stable $\infty$-operads. Then the corresponding restriction map
\[ i^*: \Alg_{\cat{B}^{\otimes}} \to \Alg_{\cat{A}^{\otimes}} \]
has a fully-faithful left adjoint.
\end{lemma}
\begin{proof}
Limits and filtered colimits of (stable) algebras are calculate objectwise (by \cite[3.2.2.4]{lurie:2017} and \cite[3.2.3.1]{lurie:2017}), and so $i^*$ preserves those (co)limits. Thus, by the Adjoint Functor Theorem~\cite[5.5.2.9]{lurie:2009}, the desired left adjoint $i_!$ exists. This left adjoint must be given by the same formula as for not-necessarily-stable algebras, described in \cite[3.1.3.1]{lurie:2017}: for an $\cat{A}^{\otimes}$-algebra $X$, and object $J \in \cat{B}$, we have
\[ i_!(X)(J) = \colim_{(\un{I}_1,\dots,\un{I}_n) \in (\cat{A}^{\otimes}_{act})_{/J}} X(\un{I}_1) \smsh \dots \smsh X(\un{I}_n). \]
Since $\cat{A}^{\otimes}_{act}$ is a full subcategory of $\cat{B}^{\otimes}_{act}$, it follows that the counit is an equivalence
\[ i^*i_!(X) \weq X \]
and so $i_!$ is fully-faithful.
\end{proof}

\begin{theorem} \label{thm:morita}
Precomposition with the inclusion of stable $\infty$-operads
\[ \bar{Y}: \cat{O}^{\otimes} \to \check{\mathbb{I}}^{\otimes}_{\Alg_{\cat{O}^{\otimes}}} \]
induces an equivalence of $\infty$-categories
\[ \bar{Y}^*: \Alg_{\check{\mathbb{I}}^{\otimes}_{\Alg_{\cat{O}^{\otimes}}}} \weq \Alg_{\cat{O}^{\otimes}}. \]
In other words, $\check{\mathbb{I}}^{\otimes}_{\Alg_{\cat{O}^{\otimes}}}$ is Morita-equivalent to $\cat{O}^{\otimes}$.
\end{theorem}
\begin{proof}
Given that $\bar{Y}$ is fully-faithful, Lemma~\ref{lem:alg-res} tells us that $\bar{Y}^*$ has a fully-faithful left adjoint $\bar{Y}_!$. It is then sufficient to show that $\bar{Y}^*$ is conservative. Since equivalences of algebras are detected on the underlying $\infty$-categories, this follows from the fact that $\check{\mathbb{I}}_{\Alg_{\cat{O}^{\otimes}}}$ is an idempotent-completion of $\cat{O}$ and \cite[5.1.4.9]{lurie:2009}.
\end{proof}

We conclude this section by noting that there is a natural comparison map between a pointed compactly-generated $\infty$-category $\cat{C}$ and the $\infty$-category of stable algebras over $\mathbb{I}_{\cat{C}}^{\otimes}$.

\begin{definition}
Let $\ev: \cat{C} \to \Fun((\cat{F}^{\omega}_{\cat{C}})^{\smsh},\spectra^{\smsh})$ be adjoint to the evaluation functor
\[ \cat{C} \times (\cat{F}^{\omega}_{\cat{C}})^{\smsh} \to \spectra^{\smsh}; \quad (X,(F_1,\dots,F_n)) \mapsto (F_1(X),\dots,F_n(X)). \]
It is evident that $\ev$ takes values in symmetric monoidal functors $\cat{F}^{\omega}_{\cat{C}} \to \spectra$. Using an equivalence similar to that of (\ref{eq:CAlg}), $\ev$ can be viewed as a map
\[ \ev: \cat{C} \to \CAlg(\Fun(\cat{F}^{\omega}_{\cat{C}},\spectra)^{\otimes}). \]
The Yoneda embedding for the stable symmetric monoidal $\infty$-category $\Fun(\cat{F}^{\omega}_{\cat{C}},\spectra)^{\otimes}$ of Nikolaus~\cite[Sec. 6]{nikolaus:2016} determines a functor, there denoted $j'_{\mathrm{St}}$, from these commutative algebra objects to the $\infty$-category of stable algebras over $\Fun(\cat{F}_\cat{C}^{\omega},\spectra)^{\otimes,op}$. Composing and restricting to $\mathbb{I}_{\cat{C}}^{\otimes}$ we obtain the desired functor
\[ \eta: \cat{C} \to \Alg_{\mathbb{I}_{\cat{C}}^{\otimes}}. \]
For an object $x \in \cat{C}$, the $\mathbb{I}_{\cat{C}}^{\otimes}$-algebra $\eta(x)$ has underlying functor
\[ \spectra(\cat{C})^{op} \to \spectra; \quad Y \mapsto \Nat_{\cat{F}_{\cat{C}}}(\der_1(-)(Y), \ev(x)). \]
Note that when $\cat{C} = \Alg_{\cat{O}^{\otimes}}$ for a stable $\infty$-operad $\cat{O}^{\otimes}$, we do not necessarily expect the map $\eta$ to be an inverse to the equivalence of Theorem~\ref{thm:morita}. 
\end{definition}

\begin{remark}
We conjecture that $\eta$ is the unit of a `quasi-adjunction'~\cite[Sec. 7]{gray:1974}, or `lax $2$-adjunction' between the functors
\[ \cat{C} \mapsto \check{\mathbb{I}}^{\otimes}_{\cat{C}} \]
and
\[ \cat{O}^{\otimes} \mapsto \Alg_{\cat{O}^{\otimes}} \]
which, we conjecture, can be interpreted as relating certain $(\infty,2)$-categories of pointed compactly-generated $\infty$-categories (and reduced finitary functors) and small stable non-unital $\infty$-operads (and bimodules between them). That quasi-adjunction would consist of adjunctions of $\infty$-categories
\[ \der_* : \Fun(\cat{C},\Alg_{\cat{O}^{\otimes}}) \rightleftarrows \Bimod(\cat{O}^{\otimes},\check{\mathbb{I}}_{\cat{C}}^{\otimes}) : \Phi \]
of the form examined by the author and Greg Arone in~\cite{arone/ching:2015} in the classification of Taylor towers of functors of based spaces and spectra. To try to develop this theory any further is beyond the scope of this paper, though we introduce bimodules between $\infty$-operads in the next section.
\end{remark}

\section{Bimodules over \texorpdfstring{$\infty$}{infinity}-operads} \label{sec:bimodule}

Consider a functor $F: \cat{C} \to \cat{D}$ between pointed compactly-generated $\infty$-categories. We wish to show that the derivatives of $F$ have the structure of a \emph{bimodule} over the $\infty$-operads $\mathbb{I}_{\cat{C}}^{\otimes}$ and $\mathbb{I}_{\cat{D}}^{\otimes}$. Bimodules are studied by Lurie in \cite[3.1.2.1]{lurie:2017} under the guise of \emph{correspondences}, or \emph{$\Delta^1$-families}, of $\infty$-operads. In order to state our conjectured chain rule, we need families of $\infty$-operads indexed by an arbitrary $\infty$-category, not just by $\Delta^1$. The following is a non-unital version of \cite[2.3.2.10]{lurie:2017}.

\begin{definition} \label{def:family}
Let $S$ be an $\infty$-category. An \emph{$S$-family of non-unital $\infty$-operads} consists of a categorical fibration
\[ p: \cat{M}^{\otimes} \to S \times \Surj \]
with the following properties:
\begin{enumerate}
  \item the restriction $p_s: \cat{M}^{\otimes}_s \to \Surj$ of $p$ to each object $s \in S$ is an $\infty$-operad;
  \item for each object $(X_1,\dots,X_m) \in \cat{M}^{\otimes}_s$, each inert morphism $\alpha: \langle m \rangle \to \langle n \rangle$ in $\Surj$ has a lift
      \[ \bar{\alpha}: (X_1,\dots,X_m) \to (X'_1,\dots,X'_n) \]
      in $\cat{M}_s^{\otimes}$ such that $\bar{\alpha}$ is $p$-cocartesian (and not merely $p_s$-cocartesian);
  \item for each morphism $f: s \to s'$ in $S$, morphism $\alpha: \langle m \rangle \to \langle n \rangle$ in $\Surj$, and each pair of objects $(X_1,\dots,X_m) \in \cat{M}_{s}$, and $(Y_1,\dots,Y_n) \in \cat{M}_{s'}$, the inert maps $\bar{\rho}_i: (Y_1,\dots,Y_n) \to (Y_i)$ in $\cat{M}_{s'}^{\otimes}$ induce equivalences
      \[ \Hom_{\cat{M}^{\otimes}}((X_1,\dots,X_m),(Y_1,\dots,Y_n))_{f,\alpha} \homeq \prod_{i = 1}^{n} \Hom_{\cat{M}^{\otimes}}((X_1,\dots,X_m),(Y_i))_{f,\rho_i\alpha} \]
      between the spaces of morphisms in $\cat{M}^{\otimes}$ that project via $p$ to $f$ and the given morphisms in $\Surj$.
\end{enumerate}
\end{definition}

\begin{definition} \label{def:stable-family}
We say that an $S$-family of $\infty$-operads $\cat{M}^{\otimes} \to S \times \Surj$ is \emph{stable} if the $\infty$-operad $\cat{M}^{\otimes}_s$ is stable for each $s \in S$ and, for each $n$ and each morphism $f: s \to s'$ in $S$, the functor
\[ \Hom_{\cat{M}^{\otimes}}(-,\dots,-;-)_f : (\cat{M}_s^{op})^n \times \cat{M}_{s'} \to \spaces \]
preserves finite limits in each variable, where the multi-morphism spaces are defined in the same manner as in Remark~\ref{rem:multi}. In this case, as in Definition~\ref{def:stable-operad}, we have corresponding multi-morphism spectra $\Map_{\cat{M}^{\otimes}}(X_1,\dots,X_n;Y)_f$.
\end{definition}

\begin{definition} \label{def:bimodule}
Let $\cat{L}^{\otimes}$ and $\cat{R}^{\otimes}$ be stable non-unital $\infty$-operads. An \emph{$(\cat{L}^{\otimes},\cat{R}^{\otimes})$-bimodule} is a stable $\Delta^1$-family of non-unital $\infty$-operads $p: \cat{M}^{\otimes} \to \Delta^1 \times \Surj$ together with equivalences of $\infty$-operads
\[ \cat{M}^{\otimes}_0 \homeq \cat{R}^{\otimes}, \quad \cat{M}^{\otimes}_1 \homeq \cat{L}^{\otimes}. \]
\end{definition}

\begin{remark}
Let $\cat{M}^{\otimes}$ be an $(\cat{L}^{\otimes},\cat{R}^{\otimes})$-bimodule. Then, for $X_1,\dots,X_n \in \cat{R}$ and $Y \in \cat{L}$, the multi-morphism spectra
\[ \Map_{\cat{M}^{\otimes}}(X_1,\dots,X_n;Y) \]
have actions, on the left by the multi-morphism spectra of $\cat{L}^{\otimes}$, and on the right by the multi-morphism spectra of $\cat{R}^{\otimes}$, that form the structure of a bimodule over the coloured operads corresponding to $\cat{L}^{\otimes}$ and $\cat{R}^{\otimes}$, up to coherent homotopy.
\end{remark}

\begin{remark} \label{rem:chain}
We can think of a stable $S$-family of $\infty$-operads $\cat{M}^{\otimes} \to S \times \Surj$ as a diagram of operads and bimodules over them indexed by the $\infty$-category $S$:
\begin{enumerate}
  \item for each object $s \in S$ we have an $\infty$-operad $\cat{M}^{\otimes}_s$;
  \item for each morphism $f: s_0 \to s_1$ in $S$ we have an $(\cat{M}^{\otimes}_{s_1},\cat{M}^{\otimes}_{s_0})$-bimodule $\cat{M}^{\otimes}_f$
  \item for each $2$-simplex
  \[ \begin{diagram} \dgARROWLENGTH=1.5em
    \node{s_0} \arrow{se,b}{f} \arrow[2]{e,t}{h} \node[2]{s_2} \\
    \node[2]{s_1} \arrow{ne,b}{g}
  \end{diagram} \]
  in $S$, we have a map of $(\cat{M}^{\otimes}_{s_2},\cat{M}^{\otimes}_{s_0})$-bimodules
  \[ \cat{M}^{\otimes}_g \circ_{\cat{M}^{\otimes}_{s_1}} \cat{M}^{\otimes}_f \to \cat{M}^{\otimes}_h\]
  where the left-hand side denotes a \emph{relative composition product} of two bimodules over the $\infty$-operad $\cat{M}^{\otimes}_{s_1}$.
\end{enumerate}
The last part of this description requires some explanation. Composition in the $\infty$-category $\cat{M}^{\otimes}$ determines maps of spectra
\[ \begin{diagram} \node{\Map_{\cat{M}^{\otimes}}(Y_1,\dots,Y_n;Z)_g \smsh \Map_{\cat{M}^{\otimes}}(X_{11},\dots;Y_1)_f \smsh \dots \smsh \Map_{\cat{M}^{\otimes}}(X_{n1},\dots,X_{nk_n};Y_n)_f} \arrow{s} \\ \node{\Map_{\cat{M}^{\otimes}}(X_{11},\dots,X_{nk_n};Z)_h} \end{diagram} \]
for objects $Z \in \cat{M}_{s_2}$, $Y_i \in \cat{M}_{s_1}$ and $X_{ij} \in \cat{M}_{s_0}$. We can interpret these maps as a map of (coloured) symmetric sequences from the composition product of two bimodules to the third bimodule. Associativity of composition in the $\infty$-category $\cat{M}^{\otimes}$ implies that this map factors through a relative composition product over the `middle' $\infty$-operad $\cat{M}^{\otimes}_{s_1}$.
\end{remark}

The aim of this section is to construct families of $\infty$-operads that capture the Goodwillie derivatives of a diagram of $\infty$-categories, together with all the bimodule structures that those derivatives possess. More precisely:

\begin{goal}
Let $S$ be an $\infty$-category, and let
\[ p: \cat{X} \to S^{op} \]
be a cartesian fibration that encodes a diagram of pointed compactly-generated $\infty$-categories, $\cat{X}_s$ for each object $s \in S$, and reduced finitary functors, $\cat{X}_f: \cat{X}_s \to \cat{X}_{s'}$ for each morphism $f: s \to s'$ in $S$. We wish to construct from $p$ a stable family of $\infty$-operads
\[ p_0: \mathbb{D}^{\otimes}_{\cat{X}} \to S \times \Surj \]
with the following properties:
\begin{enumerate}
  \item for each object $s \in S$, the fibre
  \[ (p_0)_s : (\mathbb{D}_{\cat{X}}^{\otimes})_s \to \Surj \]
  is equivalent to the stable $\infty$-operad $\mathbb{I}_{\cat{X}_s}^{\otimes}$ associated to the pointed compactly-generated $\infty$-category $\cat{X}_s$, i.e. this fibre encodes the derivatives of the identity functor on $\cat{X}_s$;
  \item for each morphism $f: s \to s'$ in $S$, the bimodule
  \[ (p_0)_f: (\mathbb{D}_{\cat{X}}^{\otimes})_f \to \Delta^1 \times \Surj \]
  encodes the derivatives of the reduced finitary functor $\cat{X}_f: \cat{X}_s \to \cat{X}_{s'}$.
\end{enumerate}
In particular, when the diagram consists of a single functor $F: \cat{C} \to \cat{D}$, the stable family is an $(\mathbb{I}_\cat{D}^{\otimes},\mathbb{I}_\cat{C}^{\otimes})$-bimodule $\mathbb{D}_F^{\otimes}$ that consists of the derivatives of $F$.

When the diagram consists of a pair of functors $F: \cat{C} \to \cat{D}$, $G: \cat{D} \to \cat{E}$, and their composite $GF: \cat{C} \to \cat{E}$, we obtain the stable family of $\infty$-operads that encodes the three bimodules given by the derivatives of these three functors, together with a map of $(\mathbb{I}^{\otimes}_{\cat{E}},\mathbb{I}^{\otimes}_{\cat{C}})$-bimodules of the form
\[ \mathbb{D}^{\otimes}_{G} \circ_{\mathbb{I}^{\otimes}_{\cat{D}}} \mathbb{D}^{\otimes}_F \to \mathbb{D}^{\otimes}_{GF}. \]
In~\ref{conj:chain} we conjecture a Chain Rule which claims that this map of bimodules is an equivalence.
\end{goal}

\begin{remark} \label{rem:idea}
Here is the idea behind our construction of $\mathbb{D}^{\otimes}_{\cat{X}}$ in the case where $p:\cat{X} \to (\Delta^1)^{op}$ is a cartesian fibration that represents a single reduced finitary functor $F: \cat{C} \to \cat{D}$.

The functor $F$ determines, by precomposition, a functor
\[ F^*: \cat{F}_{\cat{D}} \to \cat{F}_{\cat{C}}; \quad G \mapsto GF \]
which is symmetric monoidal with respect to the pointwise smash product.

The functor $F^*$ in turn induces a \emph{lax} symmetric monoidal functor (i.e. a map of $\infty$-operads) between the Day convolution monoidal structures:
\[ F_*: \Fun(\cat{F}_{\cat{C}},\spectra)^{\otimes} \to \Fun(\cat{F}_{\cat{D}},\spectra)^{\otimes}; \quad {A}(-) \mapsto {A}(-F). \]
The functor $F_*$ has a left adjoint $F^{!}$ given by left Kan extension along $F^*$, and this left adjoint is a symmetric monoidal functor
\[ F^!: \Fun(\cat{F}_{\cat{D}},\spectra)^{\otimes} \to \Fun(\cat{F}_{\cat{C}},\spectra)^{\otimes}. \]
Taking opposites, we get a symmetric monoidal functor
\[ (F^!)^{op}: \Fun(\cat{F}_{\cat{D}},\spectra)^{op,\otimes} \to \Fun(\cat{F}_{\cat{C}},\spectra)^{op,\otimes}. \]
Pulling back the left action of $\Fun(\cat{F}_{\cat{C}},\spectra)^{op,\otimes}$ on itself along the map $(F^!)^{op}$, we get a $(\Fun(\cat{F}_{\cat{D}},\spectra)^{op,\otimes}, \Fun(\cat{F}_{\cat{C}},\spectra)^{op,\otimes})$-bimodule $\cat{M}^{\otimes}$ whose multi-morphism spectra are given by
\[ \Map_{\cat{M}^{\otimes}}(A_1,\dots,A_n;B) = \Nat_{\cat{F}_{\cat{C}}}(F^!B,A_1 \otimes \dots \otimes A_n). \]
Restricting to the $\infty$-operads $\mathbb{I}_{\cat{D}}^{\otimes}$ and $\mathbb{I}_{\cat{C}}^{\otimes}$, we get a bimodule $\mathbb{D}_{\cat{X}}^{\otimes}$ whose multi-morphism spectra are
\[ \Map_{\mathbb{D}_{\cat{X}}^{\otimes}}(X_1,\dots,X_n;Y) \homeq \Nat_{\cat{F}_{\cat{C}}}(F^!\der_1(-)(Y);\der_1(-)(X_1) \otimes \dots \otimes \der_1(-)(X_n)). \]
Via the universal property of the left Kan extension $F^!$ and Theorems~\ref{thm:dn} and \ref{thm:dF}, these spectra are equivalent to the derivatives
\[ \der_nF(X_1,\dots,X_n;Y). \]
Now suppose we are given an $S$-indexed diagram $\cat{X}$ of pointed compactly-generated $\infty$-categories and reduced finitary functors. Then the above procedure yields an $S^{op}$-indexed diagram of stable symmetric monoidal $\infty$-categories $\Fun(\cat{F}_{\cat{X}_s},\spectra)^{op,\otimes}$ and exact symmetric monoidal functors $(\cat{X}_f^!)^{op}$. Such a diagram can be expressed as an $S$-family of stable $\infty$-operads
\begin{equation} \label{eq:Sfam} \Fun(\cat{F}_{\cat{X}},\spectra)^{op,\otimes} \to S \times \Surj \end{equation}
whose fibre over $s \in S$ is the symmetric monoidal $\infty$-category
\[ \Fun(\cat{F}_{\cat{X}_s},\spectra)^{op,\otimes} \to \Surj. \]
(Note the reversal here: an $S^{op}$-indexed diagram of symmetric monoidal functors is encoded as an $S$-family of $\infty$-operads.)
Finally, the desired $S$-family of $\infty$-operads
\[ \mathbb{D}_{\cat{X}}^{\otimes} \to S \times \Surj \]
is the restriction of (\ref{eq:Sfam}) to the full subcategory whose fibre over $s$ is the $\infty$-operad $\mathbb{I}_{\cat{X}_s}^{\otimes} \subseteq \Fun(\cat{F}_{\cat{X}_s},\spectra)^{op,\otimes}$.
\end{remark}

Let us now carry out the precise construction envisaged in Remark~\ref{rem:idea}. For this purpose, we need functorial versions of the three main tools involved in that construction. We start with the following fibrewise mapping space construction.

\begin{definition} \label{def:fibre-map}
Let $X \to S$ and $Y \to S$ be maps of simplicial sets. We define $\Fun_S(X,Y)$ to be the simplicial set over $S$ for which an $n$-simplex consists of:
\begin{enumerate}
  \item an $n$-simplex $\Delta^n \to S$;
  \item a map $\Delta^n \times_S X \to Y$ over $S$;
\end{enumerate}
with simplicial structure given by precomposition with simplex maps $\Delta^m \to \Delta^n$ in a standard way. The construction $\Fun_S(-,-)$ is the internal mapping object in the cartesian closed category of simplicial sets over $S$.
\end{definition}

Here is our fibrewise version of Construction~\ref{cons:objectwise}.

\begin{construction} \label{cons:objectwise-fam}
Let $p: \cat{X} \to S^{op}$ be a cartesian fibration. Define a pullback of simplicial sets
\begin{equation} \label{eq:objectwise-fam} \begin{diagram}
  \node{\Fun_{S^{op}}(\cat{X},\spectra)^{\smsh}} \arrow{s,l}{p'} \arrow{e} \node{\Fun_{S^{op}}(\cat{X},S^{op} \times \spectra^{\smsh})} \arrow{s} \\
  \node{S^{op} \times \Surj} \arrow{e} \node{\Fun_{S^{op}}(\cat{X},S^{op} \times \Surj)}
\end{diagram} \end{equation}
where the right-hand map is induced by the cocartesian fibration $\spectra^{\smsh} \to \Surj$, and the bottom map sends an $n$-simplex $\Delta^n \to S^{op} \times \Surj$ to the composite
\[ \Delta^n \times_{S^{op}} \cat{X} \to \Delta^n \to S^{op} \times \Surj. \]
\end{construction}

\begin{proposition} \label{prop:objectwise-fam}
Let $p: \cat{X} \to S^{op}$ be a cartesian fibration. Then the map
\[ p': \Fun_{S^{op}}(\cat{X},\spectra)^{\smsh} \to S^{op} \times \Surj \]
of Construction~\ref{cons:objectwise-fam} is a stable $S^{op}$-family of $\infty$-operads with the following properties:
\begin{enumerate} \renewcommand{\theenumi}{(\alph{enumi})}
  \item the map $p'$ is a cocartesian fibration;
  \item the fibre of $p'$ over $s \in S$ is the symmetric monoidal $\infty$-category $\Fun(\cat{X}_s,\spectra)^{\smsh}$ of Construction~\ref{cons:objectwise};
  \item the multi-morphism spectra for $p'$, over $f: s' \to s$ in $S^{op}$, are given by
  \[ \Map_{\Fun_{S^{op}}(\cat{X},\spectra)^{\smsh}}(F_1,\dots,F_n;G)_f \homeq \Nat_{\cat{X}_s}(\cat{X}_f^*(F_1) \smsh \dots \smsh \cat{X}_f^*(F_n), G) \]
  where $\cat{X}_f^*$ denotes precomposition with the functor $\cat{X}_f: \cat{X}_s \to \cat{X}_{s'}$ in the diagram of $\infty$-categories classified by the cartesian fibration $p$.
\end{enumerate}
\end{proposition}
\begin{proof}
Firstly, it follows from \cite[3.2.2.12]{lurie:2009} that the right-hand vertical map in (\ref{eq:objectwise-fam}) is a cocartesian fibration, hence the pullback $p'$ is too, giving condition (a). Comparing with Construction~\ref{cons:objectwise}, the fibres of $p'$ are indeed the symmetric monoidal $\infty$-categories $\Fun(\cat{X}_s,\spectra)$, satisfying (b).

To see that $p$ is a stable $S^{op}$-family of $\infty$-operads, we check conditions (i)-(iii) of Definition~\ref{def:family}, and verify the conditions on Definition~\ref{def:stable-family}.

We have already seem that the fibres of $p'$ are stable $\infty$-operads, which gives (i). Since $p'$ is a cocartesian fibration, the $p'_s$-cocartesian lifts of inert morphisms in $\Surj$ are also $p'$-cocartesian, which gives (ii). We also know that $p'$ classifies an $S^{op} \times \Surj$-indexed diagram of $\infty$-categories of the form
\[ (s,\langle n \rangle) \mapsto \Fun(\cat{X}_s,\spectra^{\smsh}_{\langle n \rangle}) \]
where the arrows in this diagram are those induced by pulling back along $\cat{X}_f: \cat{X}_{s} \to \cat{X}_{s'}$, for edges $f: s' \to s$ in $S^{op}$, and by the functors
\[ \bar{\alpha}: \spectra^{\smsh}_{\langle n' \rangle} \to \spectra^{\smsh}_{\langle n \rangle} \]
associated to a morphism $\alpha: \langle n' \rangle \to \langle n \rangle$ by the cocartesian fibration $\spectra^{\smsh} \to \Surj$.

The mapping space in the $\infty$-category $\Fun_{S^{op}}(\cat{X},\spectra)^{\smsh}$, over some morphism
\[ (f,\alpha): (s',\langle n' \rangle) \to (s,\langle n \rangle) \]
in $S^{op} \times \Surj$, is then given by
\[ \begin{split} \Hom_{\Fun_{S^{op}}(\cat{X},\spectra)^{\smsh}} &((F_1,\dots,F_{n'}),(G_1,\dots,G_{n}))_{(f,\alpha)} \\ &\homeq \Hom_{\Fun(\cat{X}_{s},\spectra^{\smsh}_{\langle n \rangle})}(\bar{\alpha}(F_1,\dots,F_{n'})\cat{X}_f,(G_1,\dots,G_{n})). \end{split} \]
Since $\Fun(\cat{X}_{s},\spectra^{\smsh}_{\langle n \rangle}) \homeq \Fun(\cat{X}_{s},\spectra)^{n}$, this formula yields condition (iii) of Definition~\ref{def:family} and gives us the multi-morphism spaces
\[ \Hom_{\Fun_{S^{op}}(\cat{X},\spectra)^{\smsh}}(F_1,\dots,F_{n'};G)_f \homeq \Hom_{\Fun(\cat{X}_{s},\spectra)}(\cat{X}_f^*(F_1) \smsh \dots \smsh \cat{X}_f^*(F_{n'});G). \]
Finally, since $\cat{X}_f^*$ is an exact functor and $\Fun(\cat{X}_{s},\spectra)$ is a stable $\infty$-category, we get the second part of Definition~\ref{def:stable-family}, so $p'$ is a stable $S^{op}$-family of $\infty$-operads, and hence also condition (c).
\end{proof}

Let us introduce the following terminology for the type of $\infty$-operad family described in Proposition~\ref{prop:objectwise-fam}.

\begin{definition} \label{def:cocart-fam}
Let $S$ be an $\infty$-category. A \emph{cocartesian $S$-family of symmetric monoidal $\infty$-categories} is an $S$-family of $\infty$-operads
\[ p: \cat{C}^{\otimes} \to S \times \Surj \]
for which the structure map $p$ is a cocartesian fibration.
\end{definition}

\begin{remark} \label{rem:cocart-fam}
A cocartesian $S$-family of symmetric monoidal $\infty$-categories $p$ corresponds to an $S$-indexed diagram of symmetric monoidal $\infty$-categories $\cat{C}_s^{\otimes} \to \Surj$ and symmetric monoidal functors $\cat{C}_f: \cat{C}_s \to \cat{C}_{s'}$ for each edge $f: s \to s'$ in $S$. Such a family is stable if and only if each $\cat{C}_s$ is stable and each $\cat{C}_f$ is exact. In this case the multi-morphism spectra are given by
\[ \Map_{\cat{C}^{\otimes}}(X_1,\dots,X_n;Y)_f \homeq \Map_{\cat{C}_{s'}}(\cat{C}_f(X_1) \otimes \dots \otimes \cat{C}_f(X_n), Y). \]
\end{remark}

\begin{definition} \label{def:FX}
Now suppose that the cartesian fibration $p: \cat{X} \to S^{op}$ encodes a diagram of pointed compactly-generated $\infty$-categories and reduced finitary functors. Let
\[ p'': \cat{F}_{\cat{X}}^{\smsh} \to S^{op} \times \Surj \]
be the restriction of the map $p'$ of Construction~\ref{cons:objectwise-fam} to the full subcategory $\cat{F}_{\cat{X}}^{\smsh} \subseteq \Fun_{S^{op}}(\cat{X},\spectra)^{\smsh}$ generated by the reduced finitary functors $\cat{X}_s \to \spectra$ for $s \in S$. Then $p''$ is also a stable cocartesian $S^{op}$-family of symmetric monoidal $\infty$-categories because each fibre is the stable symmetric monoidal $\infty$-category $\cat{F}_{\cat{X}_s}$, and the functors $\cat{X}_f^*$ restrict to exact symmetric monoidal functors $\cat{F}_{\cat{X}_{s'}} \to \cat{F}_{\cat{X}_{s}}$.
\end{definition}

Recall that since each $\cat{F}_{\cat{X}_s}$ is not small, we have to restrict further to compact objects before applying the Day convolution. These objects are not necessarily preserved by the pullback functors $\cat{X}_f^*$ so we have to enlarge the subcategories of compact objects to include all those functors obtained by pulling back a compact object.

\begin{definition} \label{def:FXomega}
Let $p: \cat{X} \to S^{op}$ and $p'': \cat{F}_{\cat{X}}^{\smsh} \to S^{op} \times \Surj$ be as in Definition~\ref{def:FX}, and let
\[ p_1: (\hat{\cat{F}}_{\cat{X}}^\omega)^{\smsh} \to S^{op} \times \Surj \]
be the restriction of $p''$ to the symmetric monoidal subcategories generated by those reduced finitary functors $\cat{X}_s \to \spectra$ that are of the form $\cat{X}_f^*(G)$ for some compact object $G \in \cat{F}^{\omega}_{\cat{X}_{s'}}$ and some $f: s \to s'$ in $S$.
\end{definition}

\begin{proposition}
The map $p_1: (\hat{\cat{F}}_{\cat{X}}^\omega)^{\smsh} \to S^{op} \times \Surj$ is a stable cocartesian $S^{op}$-family of symmetric monoidal $\infty$-categories for which each fibre $(\hat{\cat{F}}_{\cat{X}}^{\omega})_s$ is essentially small and generates $\cat{F}_{\cat{X}_s}$ under filtered colimits.
\end{proposition}

We now wish to apply to $p_1$ a functorial version of the Day convolution symmetric monoidal $\infty$-category of Construction~\ref{cons:day-infty}. The following proposition summarizes what we need.

\begin{proposition} \label{prop:day-infty-fam}
Let $p_1: \cat{C}^{\otimes} \to S^{op} \times \Surj$ be a stable cocartesian $S^{op}$-family of symmetric monoidal $\infty$-categories such that each fibre $\cat{C}_s$ is essentially small. Let $q: \cat{D}^{\otimes} \to \Surj$ be a stable symmetric monoidal $\infty$-category for which the monoidal product commutes with colimits in each variable. Then there is a stable cocartesian $S^{op}$-family of symmetric monoidal $\infty$-categories
\[ p_2: \Fun_{S^{op}}(\cat{C},\cat{D})^{\otimes} \to S^{op} \times \Surj \]
such that
\begin{enumerate}
  \item the fibre of $p_2$ over $s \in S$ is the Day convolution symmetric monoidal $\infty$-category
    \[ (p_2)_s: \Fun(\cat{C}_s,\cat{D})^{\otimes} \to \Surj \]
    of Construction~\ref{cons:day-infty};
  \item the multi-morphism spectra over $f: s' \to s$ in $S^{op}$ are given by
    \[ \begin{split} \Map_{\Fun_{S^{op}}(\cat{C},\cat{D})^{\otimes}}(A_1,\dots,A_n;B)_f &\homeq \Map_{\Fun(\cat{C}_{s},\cat{D})}(\cat{C}_f^!(A_1) \otimes \dots \otimes \cat{C}_f^!(A_n),B) \\ &\homeq \Map_{\Fun(\cat{C}_{s},\cat{D})}(\cat{C}_f^!(A_1 \otimes \dots \otimes A_n),B) \\ &\homeq \Map_{\Fun(\cat{C}_{s'},\cat{D})}(A_1 \otimes \dots \otimes A_n, B\cat{C}_f) \end{split} \]
    for $A_1,\dots,A_n: \cat{C}_{s'} \to \cat{D}$ and $B: \cat{C}_{s} \to \cat{D}$. Here $\cat{C}_f^!$ denotes left Kan extension along the symmetric monoidal functor $\cat{C}_f: \cat{C}_{s'} \to \cat{C}_{s}$ associated to $f$ by the family $p_1$.
\end{enumerate}
\end{proposition}
\begin{proof}
There are two parts to this construction: (i) functoriality of the Day convolution in the variable $\cat{C}$; (ii) formation of the left Kan extensions $\cat{C}_f^!$.

Let $\Catinf^{\otimes}$ be the nerve of the simplicial category whose objects are the (non-unital) symmetric monoidal $\infty$-categories, and whose mapping spaces are the maximal Kan complexes in the $\infty$-categories of symmetric monoidal functors. Let $\Catinf^{\otimes,\mathrm{lax}}$ be the corresponding construction with all maps of $\infty$-operads, not only the monoidal functors.

Lurie's explicit construction of Day convolution in~\cite[2.2.6.18]{lurie:2017} is functorial with respect to $\infty$-operad maps and determines a functor
\[ \Fun(-,\cat{D})^{\otimes}: (\Catinf^{\otimes})^{op} \to \Catinf^{\otimes,\mathrm{lax}}. \]
Nikolaus proves in~\cite[3.8]{nikolaus:2016} that when $F:\cat{C}^{\otimes}_1 \to \cat{C}^{\otimes}_2$ is symmetric monoidal, and under the given conditions on $\cat{D}$, the induced map
\[ F^*: \Fun(\cat{C}_2,\cat{D})^{\otimes} \to \Fun(\cat{C}_1,\cat{D})^{\otimes} \]
admits a symmetric monoidal left adjoint $F^!$. The functor $\Fun(-,\cat{D})^{\otimes}$ above therefore takes values in the subcategory
\[ \Catinf^{\otimes,\mathrm{lax},R} \subseteq \Catinf^{\otimes,\mathrm{lax}} \]
whose morphisms are the lax monoidal functors that admit a symmetric monoidal left adjoint.

Now let $\Catinf^{\otimes,L} \subseteq \Catinf^{\otimes}$ be the subcategory whose morphisms are the symmetric monoidal functors that admit a lax monoidal left adjoint. An argument similar to that of~\cite[5.5.3.4]{lurie:2009} implies there is an equivalence of $\infty$-categories
\[ (\Catinf^{\otimes,\mathrm{lax},R})^{op} \homeq \Catinf^{\otimes,L} \]
that sends a lax monoidal functor to its symmetric monoidal left adjoint. Combining this construction with the functor $\Fun(-,\cat{D})^{\otimes}$ above, we therefore have a map
\[ \Fun^!(-,\cat{D})^{\otimes}: \Catinf^{\otimes} \to \Catinf^{\otimes} \]
given on objects by Construction~\ref{cons:day-infty}, and which sends a symmetric monoidal functor $F$ to the symmetric monoidal left adjoint $F^!$.

Our given cocartesian fibration $p_1: \cat{C}^{\otimes} \to S^{op} \times \Surj$ corresponds to a functor $S^{op} \to \Catinf^{\otimes}$. Applying $\Fun^!(-,\cat{D})^{\otimes}$, we get another functor $S^{op} \to \Catinf^{\otimes}$ and hence another cocartesian family of symmetric monoidal $\infty$-categories $p_2: \Fun_{S^{op}}(\cat{C},\cat{D})^{\otimes} \to S^{op} \times \Surj$. 

Finally, each symmetric monoidal $\infty$-category $\Fun(\cat{C}_s,\cat{D})^{\otimes}$ is stable, and each left adjoint $\cat{C}_f^!$ is exact, so the family $p_2$ is stable. By construction the multi-morphism spectra are given by the desired formulas.
\end{proof}

\begin{definition} \label{def:p2}
Let $p: \cat{X} \to S^{op}$ be a cartesian fibration as in Definition~\ref{def:FX}, and let $p_1: \hat{\cat{F}}^{\omega}_{\cat{X}} \to S^{op} \times \Surj$ be as in Definition~\ref{def:FXomega}. Applying Proposition~\ref{prop:day-infty-fam} with this $p_1$ and with $q: \spectra^{\smsh} \to \Surj$ the usual (non-unital) symmetric monoidal $\infty$-category of spectra under the smash product, we obtain a stable cocartesian $S^{op}$-family of symmetric monoidal $\infty$-categories
\[ p_2: \Fun_{S^{op}}(\hat{\cat{F}}^{\omega}_{\cat{X}},\spectra)^{\otimes} \to S^{op} \times \Surj. \]
\end{definition}

The final step is to apply a fibrewise version of Construction~\ref{cons:symmon-op}. The details of this construction are in Appendix~\ref{app:dual-cocart}; here is what we need from it.

\begin{proposition} \label{prop:symmon-op-fam}
Let $p_2: \cat{M}^{\otimes} \to S^{op} \times \Surj$ be a stable cocartesian $S^{op}$-family of symmetric monoidal $\infty$-categories corresponding to a diagram of exact symmetric monoidal functors $\cat{M}_f: \cat{M}_{s'} \to \cat{M}_s$ between stable symmetric monoidal $\infty$-categories, for each $f: s' \to s$ in $S^{op}$. Then there is a stable $S$-family of symmetric monoidal $\infty$-categories
\[ p_3: \cat{M}^{op,\otimes} \to S \times \Surj \]
such that:
\begin{enumerate}
  \item for each $s \in S$, the fibre $\cat{M}^{op,\otimes}_s \to \Surj$ is equivalent to the opposite symmetric monoidal $\infty$-category of $\cat{M}^{\otimes}_s$, as in Construction~\ref{cons:symmon-op};
  \item for each $f: s \to s'$ in $S$, we have multi-morphism spectra
  \[ \Map_{\cat{M}^{op,\otimes}}(Y_1,\dots,Y_n;X)_f \homeq \Map_{\cat{M}_s}(\cat{M}_f(X),Y_1 \otimes \dots \otimes Y_n). \]
\end{enumerate}
\end{proposition}

\begin{definition} \label{def:p3}
Let $p: \cat{X} \to S^{op}$ be a cartesian fibration as in Definition~\ref{def:FX}, and let $p_2: \Fun_{S^{op}}(\hat{\cat{F}}^{\omega}_{\cat{X}},\spectra)^{\otimes} \to S^{op} \times \Surj$ be as in Definition~\ref{def:p2}. Applying Proposition~\ref{prop:symmon-op-fam}, we obtain a stable $S$-family of symmetric monoidal $\infty$-categories
\[ p_3: \Fun_{S^{op}}(\hat{\cat{F}}^{\omega}_{\cat{X}},\spectra)^{op,\otimes} \to S \times \Surj \]
whose fibre over $s \in S$ is the symmetric monoidal $\infty$-category
\[ \Fun(\hat{\cat{F}}^{\omega}_{\cat{X}_s},\spectra)^{op,\otimes} \to \Surj \]
and with multi-morphism spectra
\[ \Map_{\hat{\cat{F}}^{\omega}_{\cat{X}_s}}(\cat{X}_f^!(A),B_1 \otimes \dots \otimes B_n) \]
for $A: \hat{\cat{F}}^{\omega}_{\cat{X}_{s'}} \to \spectra$ and $B_1,\dots,B_n: \hat{\cat{F}}^{\omega}_{\cat{X}_s} \to \spectra$, where $\cat{X}_f^!$ is the left Kan extension along the functor $\cat{X}_f^*$.
\end{definition}

\begin{definition} \label{def:D}
With $p$ and $p_3$ as in Definition~\ref{def:p3}, let $\mathbb{D}^{\otimes}_{\cat{X}}$ be the full subcategory of $\Fun_{S^{op}}(\hat{\cat{F}}^{\omega}_{\cat{X}},\spectra)^{op,\otimes}$ whose objects over $s \in S$ are those of the $\infty$-operad $\mathbb{I}^{\otimes}_{\cat{X}_s}$. Let
\[ p_0: \mathbb{D}^{\otimes}_{\cat{X}} \to S \times \Surj \]
be the restriction of $p_3$ to this subcategory. Recall that the objects of $\mathbb{I}^{\otimes}_{\cat{X}_s}$ are the functors
\[ \der_1(-)(X) : \cat{F}_{\cat{X}_s} \to \spectra \]
for $X \in \spectra(\cat{X}_s)$, and that we usually denote this object just by $X$.
\end{definition}

\begin{theorem} \label{thm:D}
Let $p: \cat{X} \to S^{op}$ be a cartesian fibration that encodes a diagram of pointed compactly-generated $\infty$-categories and reduced finitary functors. Then the map
\[ p_0: \mathbb{D}^{\otimes}_{\cat{X}} \to S \times \Surj \]
is a stable $S$-family of $\infty$-operads whose fibre over $s$ is the stable $\infty$-operad $\mathbb{I}^{\otimes}_{\cat{X}_s}$, and for which the bimodule corresponding to a morphism $f:s \to s'$ in $S$ has multi-morphism spectra given by
\[ \Map_{\mathbb{D}^{\otimes}_{\cat{X}}}(X_1,\dots,X_n;Y)_f \homeq \der_n(\cat{X}_f)(X_1,\dots,X_n;Y) \]
for $X_1,\dots,X_n \in \spectra(\cat{X}_s)$ and $Y \in \spectra(\cat{X}_{s'})$. Here $\cat{X}_f: \cat{X}_s \to \cat{X}_{s'}$ is the functor associated to $f$ by the cartesian fibration $p$.
\end{theorem}
\begin{proof}
The fact that $p_0$ is a stable $S$-family of $\infty$-operads follows from the fact that $p_3$ is, and that each $\mathbb{I}_{\cat{X}_s}^{\otimes}$ is a stable $\infty$-operad. It remains to identify the multi-morphism spectra. These are a priori given by
\[ \Nat_{\hat{\cat{F}}^{\omega}_{\cat{X}_s}}(\der_1(\cat{X}_f^*(-))(Y),\der_1(-)(X_1) \otimes \dots \otimes \der_1(-)(X_n)) \]
which can be written as
\[ \Nat_{\hat{\cat{F}}^{\omega}_{\cat{X}_s}}(\der_1(-\cat{X}_f)(Y),\der_1(-)(X_1) \otimes \dots \otimes \der_1(-)(X_n)). \]
The argument of Proposition~\ref{prop:sp-op} tells us that both the Day convolution and the natural transformation object can be calculated over $\cat{F}_{\cat{X}_s}$ instead of $\hat{\cat{F}}^{\omega}_{\cat{X}_s}$. It then follows from Theorems~\ref{thm:dn} and \ref{thm:dF} that the multi-morphism spectra are given by
\[ \der_n(\cat{X}_f)(X_1,\dots,X_n;Y) \]
as claimed.
\end{proof}

\begin{definition} \label{def:D-bimod}
Let $F: \cat{C} \to \cat{D}$ be a reduced finitary functor between pointed compactly-generated $\infty$-categories. We write
\[ \mathbb{D}_F^{\otimes} \to \Delta^1 \times \Surj \]
for the $(\mathbb{I}_{\cat{D}}^{\otimes},\mathbb{I}_{\cat{C}}^{\otimes})$-bimodule given by applying Definition~\ref{def:D} in the case that $p$ is the cartesian fibration $\cat{X} \to (\Delta^1)^{op}$ corresponding to the functor $F$. The multi-morphism spectra of the bimodule $\mathbb{D}^{\otimes}_F$ are then simply the derivatives of $F$.
\end{definition}

\begin{proposition} \label{prop:corep-fam}
The stable $S$-family of $\infty$-operads of Definition~\ref{def:D} is corepresentable (i.e. a locally cocartesian fibration) via the functors
\[ \Delta_n(\cat{X}_f): \spectra(\cat{X}_s)^n \to \spectra(\cat{X}_{s'}) \]
for $f: s \to s'$ in $S$.
\end{proposition}
\begin{proof}
From Theorem~\ref{thm:D} and Definition~\ref{def:deriv} we have
\[ \begin{split} \Map_{\mathbb{D}^{\otimes}_{\cat{X}}}(X_1,\dots,X_n;Y)_f &\homeq \der_n(\cat{X}_f)(X_1,\dots,X_n;Y) \\ &\homeq \Map_{\spectra(\cat{X}_{s'})}(Y,\Delta_n(\cat{X}_f)(X_1,\dots,X_n)) \end{split} \]
which implies the claim.
\end{proof}

\begin{remark}
Proposition~\ref{prop:corep-fam} can be interpreted as providing for the existence of maps of the form
\[ \Delta_kG(\Delta_{n_1}F, \dots, \Delta_{n_k}F) \to \Delta_{n_1+\dots+n_k}(GF) \]
that make the construction $\Delta_*$ lax monoidal (up to higher coherent homotopies).
\end{remark}

We conclude with a conjectured chain rule that generalizes that of Arone and the author for the $\infty$-categories of based spaces and spectra~\cite{arone/ching:2011}.

\begin{conjecture} \label{conj:chain}
Let $F: \cat{C} \to \cat{D}$ and $G: \cat{D} \to \cat{E}$ be reduced finitary functors between pointed compactly-generated $\infty$-categories. Let $\cat{X} \to (\Delta^2)^{op}$ be the cartesian fibration that encodes the diagram consisting of the functors $F$, $G$ and $GF$, and let
\[ \mathbb{D}_{\cat{X}}^{\otimes} \to \Delta^2 \times \Surj \]
be the associated $\Delta^2$-family of $\infty$-operads of Theorem~\ref{thm:D}. Then the corresponding map, described in Remark~\ref{rem:chain}, of $(\mathbb{I}_{\cat{E}}^{\otimes},\mathbb{I}_{\cat{C}}^{\otimes})$-bimodules,
\[ \mathbb{D}_G^{\otimes} \circ_{\mathbb{I}_{\cat{D}}^{\otimes}} \mathbb{D}_F^{\otimes} \to \mathbb{D}_{GF}^{\otimes} \]
is an equivalence.
\end{conjecture}

One of the challenges in proving~\ref{conj:chain} is identifying the underlying property of a stable $\Delta^2$-family of $\infty$-operads which implies that this associated map of bimodules is an equivalence. This property should be some version of Lurie's `flatness' condition~\cite[B.3]{lurie:2017}. Notice that \cite[3.1.4.2]{lurie:2017} describes a consequence of that assumption in the unstable case. This result suggests that we describe the chain rule in terms of induced functors between $\infty$-categories of algebras. In particular, it seems that any stable $S$-family of $\infty$-operads should determine a locally cartesian fibration that classifies those induced functors. The Chain Rule would then follow from showing that this locally cartesian fibration is in fact cartesian in the case at hand.

We note further, however, that the proof of the Chain Rule in~\cite{arone/ching:2011} depends considerably on Koszul duality between the operad $\der_*I_\cat{C}$ and the cooperad formed by the derivatives of the functor $\Sigma^\infty\Omega^\infty$. Since that duality does not play an explicit role in the current constructions, we might expect to need a new argument here.

Finally, we should note that Lurie proves a version of Conjecture~\ref{conj:chain} in \cite[6.3.2]{lurie:2017} for \emph{coderivatives} instead of derivatives. We expect the operad theory developed here to be Koszul dual, in a suitable sense, to Lurie's, in which case we might expect to deduce~\ref{conj:chain} directly from Lurie's results.

\appendix

\section{The chain rule for spectrum-valued functors} \label{app:chain}

In the proof of Theorem~\ref{thm:dF} we needed a chain rule for composites of functors $G: \cat{C} \to \spectra$ and $F: \spectra \to \spectra$. The purpose of this section is to state and prove the needed result, which is a mild generalization of \cite[1.15]{ching:2010}.

\begin{theorem} \label{thm:chain}
Let $\cat{C}$ be a pointed compactly-generated $\infty$-category and let $G: \cat{C} \to \spectra$ and $F: \spectra \to \spectra$ be reduced functors. Assume that $F$ preserves filtered homotopy colimits. Then for $X_1,\dots,X_n \in \spectra(\cat{C})$ we have
\[ \begin{diagram} \dgARROWLENGTH=1em \node{\der_n(FG)(X_1,\dots,X_n)} \arrow{s,l}{\sim} \\ \node{\prod_{\mu \in \mathsf{P}(n)} \der_k(F) \smsh \der_{n_1}(G)(\{X_i\}_{i \in \mu_1}) \smsh \dots \smsh \der_{n_k}(G)(\{X_i\}_{i \in \mu_k})} \end{diagram} \]
where the product is over the set $\mathsf{P}(n)$ of unordered partitions $\mu$ of $\{1,\dots,n\}$ into $k$ pieces $\mu_1,\dots,\mu_k$, with $n_j = |\mu_j|$.
\end{theorem}
\begin{proof}
We follow the approach of \cite{ching:2010} very closely. Indeed, many of the results proved there carry over to this more general situation with no change. Specifically, we can construct, as in \cite[2.5]{ching:2010}, a map
\[ \Delta: FG \to \prod_{\lambda} [P_{k_1},\dots,P_{k_r}]\creff_r(F)(P_{\ell_1}G,\dots,P_{\ell_r}G) \]
where $\lambda$ varies over expressions of the form
\[ n = k_1\ell_1 + \dots + k_r\ell_r. \]
with $k_i$ and $\ell_i$ positive integers such that $\ell_1 < \dots < \ell_r$. We can also prove, as in \cite[4.2]{ching:2010}, that $\Delta$ induces an equivalence on $D_n$, and hence on \ord{$n$} derivatives. Moreover, we can show, as in the proof of \cite[2.6]{ching:2010}, that the \ord{$n$} derivative of the functor
\[ [P_{k_1},\dots,P_{k_r}]\creff_r(F)(P_{\ell_1}G,\dots,P_{\ell_r}G) \]
is equivalent to the \ord{$n$} derivative of the $n$-homogeneous functor
\begin{equation} \label{eq:der-hom} (\der_kF \smsh (D_{\ell_1}G)^{\smsh k_1} \smsh \dots \smsh (D_{\ell_r}G)^{\smsh k_r})_{h\Sigma_{k_1} \times \dots \times \Sigma_{k_r}} \end{equation}
where $k = k_1+\dots+k_r$. It now remains to calculate this \ord{$n$} derivative at an $n$-tuple $(X_1,\dots,X_n)$ in $\spectra(\cat{C})$.

Since all the functors involved here are homogeneous, and thus factor via $\Sigma^\infty_{\cat{C}}: \cat{C} \to \spectra(\cat{C})$, we can assume without loss of generality that $\cat{C}$ is stable. Using the equivalence
\[ D_{\ell}G(X) \homeq \der_{\ell}G(X,\dots,X)_{h\Sigma_{\ell}} \]
we can write the functor (\ref{eq:der-hom}) as mapping $X$ to
\[ (\der_kF \smsh \der_{\ell_1}G(X,\dots,X)^{\smsh k_1} \smsh \dots \smsh \der_{\ell_r}G(X,\dots,X)^{\smsh k_r})_{hH(\lambda)} \]
where $H(\lambda)$ denotes the subgroup $(\Sigma_{\ell_1} \wr \Sigma_{k_1}) \times \dots \times (\Sigma_{\ell_r} \wr \Sigma_{k_r})$ of $\Sigma_n$ formed from wreath products. It's convenient to rewrite the formula above as
\[ (\der_kF \smsh \der_{n_1}G(X,\dots,X) \smsh \dots \smsh \der_{n_k}G(X,\dots,X))_{hH(\lambda)} \]
where $n_1,\dots,n_k$ are the numbers $\ell_1,\dots,\ell_r$ with $\ell_i$ repeated $k_i$ times.

Now when $E: \cat{C}^n \to \spectra$ is a multilinear functor, the \ord{$n$} derivative of the functor $X \mapsto E(X,\dots,X)$ at $(X_1,\dots,X_n)$ can be written as
\[ \prod_{\sigma \in \Sigma_n} E(X_{\sigma(1)},\dots,X_{\sigma(n)}) \]
It follows from all these calculations that $\der_n(FG)(X_1,\dots,X_n)$ can be expressed as
\[ \prod_{\lambda} \left(\prod_{\sigma \in \Sigma_n} \der_kF \smsh \der_{n_1}G(X_{\sigma(1)},\dots,X_{\sigma(n_1)}) \smsh \dots \smsh \der_{n_k}G(X_{\sigma(n-n_k+1)},\dots,X_{\sigma(n)}) \right)_{hH(\lambda)}. \]
It remains to identify this expression with the formula stated in the Theorem, which we do this by showing that a choice of $\lambda$, together with a coset $[\sigma]$ of $H(\lambda)$ in $\Sigma_n$, uniquely corresponds to an unordered partition of $\{1,\dots,n\}$.

In one direction, we map the pair $(\lambda,[\sigma])$ to the partition whose pieces are the sets $(\sigma(1),\dots,\sigma(n_1)), \dots, (\sigma(n-n_k+1),\dots,\sigma(n))$. On the other hand, given an unordered partition $\mu$, let $k_j$ be the number of pieces of size $\ell_j$ (determining $\lambda$). If we put the pieces of $\mu$ in ascending size order, and concatenate them, we get an element $\sigma \in \Sigma_n$ which determines a coset of $H(\lambda)$. This is well-defined because changing the order of elements within each piece, or the order of pieces of the same size, only changes $\sigma$ by an element of $H(\lambda)$. It is a simple check that these two constructions are inverse, setting up the desired correspondence. Via this bijection, the expression given above for $\der_n(FG)(X_1,\dots,X_n)$ corresponds with the desired formula.
\end{proof}

\section{Fibrewise duals for cocartesian fibrations} \label{app:dual-cocart}

The goal of this section is to prove Proposition~\ref{prop:symmon-op-fam}. Recall that the opposite of a symmetric monoidal $\infty$-category is built from the following general construction of~\cite{barwick/glasman/nardin:2018}.

\begin{definition}
Let $q: Y \to T$ be a cocartesian fibration. Then there is a cartesian fibration
\[ q^{\vee}: Y^{\vee} \to T^{op} \]
that classifies the same diagram of $\infty$-categories as $q$.

The construction $q \mapsto q^{\vee}$ is functorial in the following way. Let $p: Y \to Y'$ be a map between cocartesian fibrations over $T$ that takes cocartesian edges to cocartesian edges. Then $p$ induces a map
\[ p^{\vee}: Y^{\vee} \to Y'^{\vee} \]
between cartesian fibrations over $T^{op}$ that takes cartesian edges to cartesian edges.
\end{definition}

\begin{definition} \label{def:fibre-dual}
Let $S$ be an $\infty$-category, and let $p: Y \to S^{op} \times T$ be a map of simplicial sets such that:
\begin{enumerate}
  \item for each $s \in S$, the fibre $p_s: Y_s \to T$ is a cocartesian fibration;
  \item each $p_s$-cocartesian lift of a morphism in $T$ is also $p$-cocartesian.
\end{enumerate}
It follows that $p$ is a map between cocartesian fibrations over $T$ that takes  cocartesian edges to cocartesian edges. Therefore $p$ induces a map
\[ p^{\vee}: Y^{\vee} \to (S^{op} \times T)^{\vee} \homeq S^{op} \times T^{op} \]
between cartesian fibrations over $T^{op}$ that takes cartesian edges to cartesian edges.

Taking opposites, we obtain another map
\[ p^{\vee,op} : Y^{\vee,op} \to S \times T \]
between cocartesian fibrations over $T$ that takes cocartesian edges to cocartesian edges.
\end{definition}

\begin{definition} \label{def:fibre-dual-2}
Now let $p: \cat{M}^{\otimes} \to S^{op} \times \Surj$ be a stable cocartesian $S^{op}$-family of symmetric monoidal $\infty$-categories. Then $p$ satisfies the conditions in Definition~\ref{def:fibre-dual}, and so determines a map
\[ p^{\vee,op}: \cat{M}^{op,\otimes} \to S \times \Surj \]
between cocartesian fibrations over $\Surj$ that takes cocartesian edges to cocartesian edges. By factoring $p^{\vee,op}$ as an acyclic cofibration followed by a fibration in the cocartesian model structure of (naturally-)marked simplicial sets over $\Surj$ (see~\cite[3.1.3.9]{lurie:2009}), we can assume that $p^{\vee,op}$ is a fibration in that model structure.
\end{definition}

We now check that $p^{\vee,op}$ satisfies the requirements of Proposition~\ref{prop:symmon-op-fam}.

\begin{proposition} \label{prop:fibre-dual-1}
The map $p^{\vee,op}: \cat{M}^{op,\otimes} \to S \times \Surj$ constructed in Definition~\ref{def:fibre-dual-2} is a stable $S$-family of symmetric monoidal $\infty$-categories. The fibre of $p^{\vee,op}$ over $s \in S$ is equivalent to the opposite symmetric monoidal $\infty$-category $(p_s)^{\vee,op}: \cat{M}_s^{op,\otimes} \to \Surj$ and the multi-morphism spectra are given, for $f: s \to s'$ in $S$, by
\[ \Map_{\cat{M}^{op,\otimes}}(X_1,\dots,X_n;Y)_f \homeq \Map_{\cat{M}_s}(\cat{M}_f(Y),X_1 \otimes \dots \otimes X_n) \]
where $\cat{M}_f: \cat{M}_{s'} \to \cat{M}_{s}$ is the symmetric monoidal functor associated to $f$ by the original cocartesian family $p$. 
\end{proposition}
\begin{proof}
We check that $p^{\vee,op}$ satisfies the conditions of Definition~\ref{def:family}. Since $p^{\vee,op}$ is a cocartesian fibration, it follows from~\cite[3.1.5.1]{lurie:2009} that it is also a categorical fibration.

In~\cite[1.4]{barwick/glasman/nardin:2018}, it is shown that $(-)^{\vee,op}$ is a self-equivalence of the $\infty$-categories of cocartesian fibrations over $\Surj$, so this construction preserves fibre sequences. Therefore the fibre $(p^{\vee,op})_s$ is equivalent to $(p_s)^{\vee,op}$ as claimed. This observation also verifies condition (i) of Definition~\ref{def:family}.

Since $\cat{M}^{op,\otimes} \to \Surj$ is a cocartesian fibration, each morphism $\alpha: \langle n \rangle \to \langle m \rangle$ in $\Surj$ has, for $(X_1,\dots,X_n) \in \cat{M}^{\otimes}_{\langle n \rangle}$, a cocartesian lift
\[ \bar{\alpha}: (X_1,\dots,X_n) \to (X'_1,\dots,X'_m) \]
where $X'_i = \bigotimes_{\alpha(i) = j} X_j$.

By construction, $p^{\vee,op}$ maps $\bar{\alpha}$ to an edge of $S \times \Surj$ that is also cocartesian over $\Surj$. It follows from (the dual of) \cite[2.4.1.3(3)]{lurie:2009} that each $\bar{\alpha}$ is also $p^{\vee,op}$-cocartesian, which implies (ii).

For (iii), consider objects $(X_1,\dots,X_n) \in \cat{M}_s$ and $(Y_1,\dots,Y_m) \in \cat{M}_{s'}$, a morphism $f: s \to s'$ in $S$, and a morphism $\alpha: \langle n \rangle \to \langle m \rangle$ in $\Surj$.

We then have
\[ \Hom_{\cat{M}^{op,\otimes}}((X_1,\dots,X_n),(Y_1,\dots,Y_m))_f \homeq \Hom_{(\cat{M}^{\otimes})^{\vee}}((Y_1,\dots,Y_m),(X_1,\dots,X_n))_{f^{op}}. \]
Since the cartesian fibration $(\cat{M}^{\otimes})^{\vee} \to \Surj^{op}$ encodes the same diagram of $\infty$-categories and functors as the cocartesian fibration $\cat{M}^{\otimes} \to \Surj$, the part of the above mapping space over $\alpha$ is equivalent to
\[ \Hom_{\cat{M}^{\otimes}_{\langle m \rangle}} ((Y_1,\dots,Y_m),(X'_1,\dots,X'_m))_{f^{op}} \]
and hence to
\[ \prod_{i = 1}^{m} \Hom_{\cat{M}_s}(\cat{M}_f(Y_i),\bigotimes_{\alpha(j) = i}X_j). \]
Using this decomposition, we verify part (iii) of \ref{def:family}, complete the proof that the family $\cat{M}^{op,\otimes}$ is stable, and obtain the desired description of its multi-morphism spectra.
\end{proof}

\providecommand{\bysame}{\leavevmode\hbox to3em{\hrulefill}\thinspace}
\providecommand{\MR}{\relax\ifhmode\unskip\space\fi MR }
\providecommand{\MRhref}[2]{%
	\href{http://www.ams.org/mathscinet-getitem?mr=#1}{#2}
}
\providecommand{\href}[2]{#2}

\end{document}